\documentstyle{amsppt}

\voffset=-24pt
\hoffset=2.5mm
\magnification=\magstep1
\parindent=12pt
\baselineskip=14pt
\hfuzz=5pt
\mathsurround 2pt
\parskip 6pt
\overfullrule=0pt

% blackboard bold font
\define\A{{\Bbb A}}
\define\C{{\Bbb C}}
\define\R{{\Bbb R}}
\define\Q{{\Bbb Q}}
\define\Z{{\Bbb Z}}
% greek
\define\a{\alpha}
\redefine\b{\beta}

\define\e{\epsilon}
\predefine\l{\ll}
\redefine\l{\lambda}
\redefine\o{\omega}
\define\ph{\varphi}

\define\s{\sigma}
\redefine\P{\Phi}
\predefine\Sec{\S}
\redefine\L{\Lambda}
% other

\define\back{\backslash}

\define\lra{\longrightarrow}
\redefine\tt{\otimes}

\define\liminv#1{\underset{\underset{#1}\to\leftarrow}\to\lim\,}
\define\limdir#1{\underset{\underset{#1}\to\rightarrow}\to\lim\,}

\define\isoarrow{\ {\overset{\sim}\to{\longrightarrow}}\ }

\define\ord#1{{\underset{#1}\to{\text{ ord }}\,}}

\parindent=12pt
\baselineskip=14pt
\define\F{\Bbb F}
\define\K{\Cal K}

\define\fp{{\F_p}}
\define\fps{{\F_{p^2}}}
\define\zps{{\Z_{p^2}}}
\define\zpsx{{\Z_{p^2}^\times}}
\define\qps{{\Q_{p^2}}}
\define\fn{\frak n}
\define\tfn{\fn_\tl}
\define\fm{\frak m}
\define\fr{\frak r}
\define\ol{\Cal O_L}
\define\otl{\Cal O_\tl}
\define\om{\Cal O_M}

%\define\ker{\text{ker}}
\define\nass{\noalign{\smallskip}}

\redefine\H{H}

\redefine\j{\text{\bf j}}
\define\y{\text{\bf y}}
\define\und#1{\underline{#1}}
\redefine\O{\Omega}
\parskip=12pt

\define\End{\text{\rm End}}
\define\Hom{\text{\rm Hom}}
\define\Sch{\text{\rm Sch}}
\define\Lie{\text{\rm Lie}}

\redefine\L{\Cal L}  %%\Bbb L might be better...
\define\tl{{\tilde{L}}}
\define\im{\text{im}}
\predefine\volold{\vol}
\redefine\vol{\text{\rm vol}}
\define\pr{\text{\rm pr}}
\define\diag{\text{\rm diag}}

%\bye

%\input siegbib

%%%%%%%
%%
%%    reference macros for
%%    Cycles on Siegel 3-folds and derivatives of Eisenstein series
%%
%%%%%%%

%\documentstyle{amsppt}
%\input siegtop
%\input siegbib

\define\boutotcarayol{{\bf 1}} %
\define\deligne{{\bf 2}}
\define\delignetwo{{\bf 3}}               %
\define\eichler{{\bf 4}}
\define\fulton{{\bf 5}}
\define\genestier{{\bf 6}}   %letter to R.
\define\grosskeating{{\bf 7}}
\define\hashimoto{{\bf 8}}
\define\hashibuki{{\bf 9}}
\define\ibukatoort{{\bf 10}}
\define\kaiser{{\bf 11}}
\define\katoort{{\bf 12}}
\define\kitaokanote{{\bf 13}}
%\define\kitaokaone{{\bf }} % best ref for the def of rep densities, e.g.,
book?
\define\kitaokatwo{{\bf 14}}   %Proc Japan Acad.
%\define\kitaokathree{{\bf }}  %Advanced Studies article
\define\kitaokabook{{\bf 15}}
\define\kottwitztam{{\bf 16}} % Tamagawa numbers paper
\define\kottwitz{{\bf 17}}   % JAMS
\define\cycles{{\bf 18}}     % Duke paper
\define\ked{{\bf 19}}         % Annals paper
%\define\kmone{{\bf }}        % Math Ann. I
%\define\kmtwo{{\bf }}         % Math Ann II
\define\kmihes{{\bf 20}}
\define\krhb{{\bf 21}}         % Hilbert-Blumenthal case
\define\krdrin{{\bf 22}}       % Drinfeld case
\define\moretbailly{{\bf 23}}
\define\oort{{\bf 24}}
\define\ribetone{{\bf 25}}
\define\ribettwo{{\bf 26}}
\define\serre{{\bf 27}}
\define\shih{{\bf 28}}      %
\define\shimura{{\bf 29}}
\define\stamm{{\bf 30}}
\define\sweet{{\bf 31}}       % rdps for metaplectic case
\define\tonghai{{\bf 32}}     % rep. densities papers

%\bye

\phantom{.}
\vskip 1in
%\vskip0.9 in

\centerline{\bf Cycles on Siegel 3-folds and derivatives of
Eisenstein series}
\smallskip
\centerline{by}
\centerline{Stephen S. Kudla\footnote{NSF grant number
DMS-9622987}}
\centerline{and}
\centerline{Michael Rapoport}
%\vfill\eject
\bigskip
\subheading{Introduction}
\par\noindent
This is the first of two companion papers in which we try to
generalize the results of one of us, \cite{\ked}, to higher
dimensions. In that paper, it was established that the nonsingular
Fourier coefficients of the derivative at 0 of certain incoherent
Eisenstein series on the metaplectic group in 4 variables are
closely related to the value of the height pairing of a pair of
arithmetic cycles on a Shimura curve. Both sides of the identity to
be proved turned out to be a sum of terms enumerated by the places
of $\Q$. In the present papers we are only concerned with the
situation at a nonarchimedean prime.
\hfill\break
It is a hope,
already expressed in \cite{\ked}, that a similar relation holds in
general between the value of the derivative at 0 of certain
incoherent Eisenstein series on the metaplectic group in $2n$
variables and the height pairing of suitable arithmetic cycles on
Shimura varieties associated to orthogonal groups of signature
$(n-1, 2)$. This would constitute an arithmetic analogue of the
result of the first author \cite{\cycles} which relates the {\it
value}\/ at 1/2 of certain {\it coherent}\/ Eisenstein series with
the intersection pairing on suitable classical cycles on these
Shimura varieties.
\hfill\break
A first difficulty with this general
program is that models over the integers of the Shimura varieties
associated to orthogonal groups are not well understood. For low
values of $n$ there are, however, exceptional isomorphisms which
relate the groups in question to symplectic groups, and the Shimura
varieties associated to them have integral models which one can
investigate. In the present paper we are concerned with the
exceptional isomorphism which relates the orthogonal group of
signature $(3,2)$ with the symplectic group in 4 variables. In the
companion paper \cite{\krhb} we are concerned with the Shimura variety
associated to an orthogonal group of signature $(2,2)$ which is
related to certain Hilbert-Blumenthal surfaces.

Let us now be more specific about the contents of this paper. Let
$B$ be an indefinite quaternion algebra over $\Q$. Let $C=M_2(B)$
and put
$$V=\{ x\in C \mid\ x'=x,\ {\roman{tr}}^0(x)=0\}\ \ ,\tag0.1$$
where $x\mapsto x'= {}^t x^{\iota}$ is the involution on $C$
induced by the main involution on $B$. Then $(V,q)$, with $q$ defined
by $x^2=q(x)\cdot 1$, is a quadratic space of signature $(3,2)$ and
the group $G=GSpin(V)$ of $V$ can be identified
with a twisted form of the group of symplectic similitudes in 4
variables.
Let ${\Cal D}$ be
the space of oriented negative 2-planes in $V(\R)$ and let $K$ be a compact
open
subgroup of $G(\A_f)$. Then, the Shimura variety $\text{\rm Sh}(G,\Cal D)_K$,
whose complex points are given by
$${\roman{Sh}}(G, {\Cal D})_K(\C)= G(\Q)\setminus
[{\Cal D}\times G({\Bbb A}_f)/K], \tag0.2$$
is a (twisted) version of the Siegel 3-fold over $\Q$. For example, the case
of the split quaternion algebra $B=M_2(\Q)$ yields the usual Siegel modular
variety of genus 2.
\hfill\break
The exceptional isomorphism of $G=GSpin(V)$ with a form of $GSp_4$ plays a
fundamental role
throughout the paper. In particular, we use it to construct a good
integral model of $\text{\rm Sh}(G,\Cal D)_K$.
More precisely, we fix a prime $p>2$ such that $B$ is unramified at $p$
and take $K$ of the form $K=K^p.K_p$, where $K^p\subset G({\Bbb
A}_f^p)$ is sufficiently small and where $K_p$ is the natural
maximal compact open subgroup of $G(\Q_p)$. Then we use the modular
interpretation of ${\roman{Sh}}(G, {\Cal D})_K$ to construct a
smooth model ${\Cal M}$ over ${\roman{Spec}}\, \Z_{(p)}$, as a
parameter space of certain abelian varieties with additional
structure.

Algebraic cycles on ${\roman{Sh}}(G, {\Cal D})_K$ were defined
analytically in \cite{\cycles} as follows. For $x\in V^n$ let
$q(x)=\frac12\big((x_i,x_j)\big)\in {\roman{Sym}}_n(\Q)$, be the
matrix of inner products of the compnents of $x$ for the symmetric
bilinear form $(\ ,\ )$ associated to $q$. Assume that $d=q(x)$ is
positive-definite (hence $n\le 3$), and let ${\Cal D}_x$ be the
subspace of oriented negative 2-planes orthogonal to all entries of
$x$. Let $G_x$ be the pointwise stabilizer of $x$. Then
${\roman{Sh}}(G_x, {\Cal D}_x)$ is a sub-Shimura variety of
${\roman{Sh}}(G, {\Cal D})$, and thus defines a cycle of
codimension $n$ in ${\roman{Sh}}(G, {\Cal D})_K$. These cycles are
a special case of the totally geodesic cycles in locally symmetric
spaces studied in \cite{\kmihes} and elsewhere. A slight
generalization of the previous construction yields a cycle
$Z(d,\omega;K)$ of ${\roman{Sh}}(G, {\Cal D})_K$ which is
associated to any positive definite $d\in {\roman{Sym}}_n(\Q)$ and
any $K$-invariant compact open subset $\omega$ of $V({\Bbb
A}_f)^n$.
\hfill\break
The next step is to give a modular definition of these cycles.
First, for one of the abelian varieties parametrized by ${\Cal M}$,
we define the notion of a {\it special endomorphism}
(Definition~2.1). The space of such endomorphisms is a finitely
generated free $\Z$-module equipped with a quadratic form $q$. The
cycle $Z(d,\o;K)$ ($=Z(d,\o)$ if $K$ is fixed) is then obtained by
imposing an $n$-tuple $\j$ of special endomorphisms such that
$q(\j)=d$, and satisfying an additional compatibility with respect
to $\o$. If $\omega$ satisfies an integrality condition at $p$,
this definition can be used to extend the cycle $Z(d,\o)$ to a
cycle ${\Cal Z}(d,\omega)$ for the integral model ${\Cal M}$ of the
Shimura variety. Here, by a cycle on ${\Cal M}$, we mean a scheme
which maps by a finite unramified morphism to ${\Cal M}$. At this
point we meet a very important problem: in contrast to ${\Cal M}$,
the cycles ${\Cal Z}(d, \omega)$ will no longer be smooth, in
general. In fact, they often are not flat over $\Z_{(p)}$ and may
even have embedded components. Our justification for our choice of
this integral extension of the classical cycles is that their
definition is very simple, has a nice inductive structure with
respect to intersection, and that we are able to prove something
about them. Before stating these results, we note that, while the
arithmetic cycles ${\Cal Z}(d,\omega)$ can be defined for any $d\in
{\roman{Sym}}_n(\Q)$, any $n$,  they are nonempty only when $d$ is
positive semidefinite and with coefficients in $\Z_{(p)}$.

We fix positive integers $n_1,\ldots, n_r$ with $n_1+\ldots
+n_r=4$ and,  for each $i$, we choose a positive definite $d_i\in
{\roman{Sym}}_{n_i}(\Z_{(p)})$ and a $K$-invariant open compact
subset $\omega_i\in V({\Bbb A})^{n_i}$. The resulting cycles
${\Cal Z}(d_i, \omega_i)$ on ${\Cal M}$ have generic fibres of
codimension $n_i$. We form the fibre product
$${\Cal Z}={\Cal Z}(d_1, \omega_1)\times_{\Cal M}\ldots\times_{\Cal M}
{\Cal Z}(d_r, \omega_r)\ \ .\tag0.3$$ To each point $\xi$ of ${\Cal
Z}$, we then associate its {\it fundamental matrix}\/ $T_{\xi}\in
{\roman{Sym}}_4(\Z_{(p)})_{\ge 0}$, defined by $T_\xi = q(\j)$ where
$\j=(\j_1,\dots,\j_r)$ is the $4$-tuple of special endomorphisms imposed at
a point of the fiber product. Note that the diagonal blocks of $T$ are
$(d_1,\dots,d_r)$. The function $\xi\mapsto T_{\xi}$ is
locally constant for the Zariski topology on ${\Cal Z}$ and induces
a disjoint sum decomposition in which the summands are again special cycles of
a definite kind,
$${\Cal Z}= \coprod_{\matrix
T\in {\roman{Sym}_4}(\Z_{(p)})_{\ge 0}\\
{\roman{diag}}(T)=(d_1,\ldots, d_r)\\
\endmatrix}
{\Cal Z}(T, \omega)\ \ .\tag0.4$$ Here $\omega=\omega_1\times
\ldots\times
\omega_r$.
The summand on the right corresponding to $T$ is the set of points
$\xi$ where $T_{\xi}=T$. This decomposition illustrates the
inductive nature of the special cycles mentioned above.
\hfill\break
The decomposition (0.4) bears some formal similarity to the
partitioning into isogeny classes that occurs in the approach of
Langlands-Kottwitz to the calculation of the zeta function of a
Shimura variety. In that approach the stable conjugacy class of the
Frobenius endomorphism is the most basic invariant of an isogeny
class. In our context this role is played by the fundamental
matrix. One of our discoveries is that the fundamental matrix and
more specifically its divisibility by $p$ governs the intersection
behaviour of the special cycles. Thus if ${\roman{det}}(T_{\xi})\ne
0$, i.e.\ $T=T_{\xi}$ is positive definite, then the point $\xi$
lies in characteristic $p$ and is not the specialization of a point
of ${\Cal Z}$ in characteristic 0. In this case,  the connected
component $\Cal Z(T,\o)$ of ${\Cal Z}$ containing $\xi$ consists
entirely of supersingular points of ${\Cal M}$. Contrary to what
one might expect, however, the condition $\det(T_\xi)\ne 0$ is not
sufficient to ensure that $\xi$ is an isolated point of
intersection. One of our main results is the characterization of
when this is the case.
\proclaim{Theorem 0.1}
Let $\xi\in {\Cal Z}$ with ${\roman{det}}(T_{\xi})\ne 0$. Then
$\xi$ is an isolated intersection point if and only if $T_{\xi}$
represents 1 over $\Z_p$. In this case the underlying abelian
variety is isomorphic to a power of a supersingular elliptic curve.
\endproclaim

When $T=T_{\xi}$ does not represent 1 over $\Z_p$ (but still is
positive definite), then the connected component $\Cal Z(T,\o)$ of ${\Cal Z}$
containing $\xi$ is a union of projective lines and, in fact, one can
enumerate these lines. It turns out that the more divisible
$T_{\xi}$ is by $p$, the more components there will be. A more
thorough analysis of the set of irreducible components can be found
in \cite{\krhb}. We point out that this phenomenon of excess
intersection does not occur in the case of Shimura curves at a
place of good reduction \cite{\ked},  but it does occur at a place of
bad reduction \cite{\krdrin}.
\hfill\break
With the
previous notation let us put
$$<{\Cal Z}(d_1,\omega_1),\ldots,
{\Cal Z}(d_r,\omega_r)>_p^{\roman{proper}} =
\sum_{\matrix
\xi\in{\Cal Z},\\
\xi\ \text{isolated}\\
\endmatrix}
e(\xi)\ \ ,\tag0.5$$ where each isolated intersection point $\xi$
appears with multiplicity $e(\xi)= {\roman{lg}}\, {\Cal O}_{{\Cal
Z}, \xi}$, the length of the local ring of
$\Cal Z$ at $\xi$.

We next come to the relation with Eisenstein series, for which we refer the
reader to section~8 or the first part of \cite{\ked} for more details.
Let $W$ be a symplectic space over $\Q$ of dimension 8
and let
$$W=W_1+\ldots +W_r\tag0.6$$
be a decomposition of $W$ into symplectic spaces $W_i$
of dimension $2n_i$. Let
$$i:Mp_{1,{\Bbb A}}\times \ldots\times Mp_{r,{\Bbb A}}\to Mp_{\Bbb A}
\tag0.7$$
be the corresponding embedding of metaplectic groups. Let $\Phi(s)$
be the standard section of the induced representation $I(s,\chi_V)$
of $Mp_{\Bbb A}$ which is of the form $\Phi(s)
=\Phi_\infty(s)\otimes \Phi_f(s)$.  Here the finite part is
associated to the Schwartz function ${\roman{char}}\, \omega\in
S(V({\Bbb A}_f)^4)$ under the natural map $S(V({\Bbb A}_f)^4)
\rightarrow I_f(0,\chi)$ defined via the Weil representation.
Similarly, the component $\P_\infty(s)$ at $\infty$ is associated
to the Gaussian for the $5$-dimensional quadratic space $V'(\R)$
over $\R$ of signature $(5,0)$ under the map
$S(V'(\R)^4)\rightarrow I_\infty(0,\chi)$. Thus the section $\P(s)$
is determined by
$$\o =\o_1\times \dots\times \o_r,\tag0.8$$
and is incoherent in the sense of \cite{\ked}. In particular, for $h\in
Mp_\A(W)$,
the corresponding Eisenstein series $E(h,s,\P)$ vanishes at the center of
symmetry $s=0$. For any $(h_1,\ldots, h_r)\in Mp_{1,\R}\times
\ldots \times Mp_{r,\R}$ we put
$$F_{d_1,\ldots, d_r}(h_1,\ldots,h_r,\Phi)_p^{\roman{proper}} =
\sum_{\matrix
T\in{\roman{Sym}}_4(\Z_{(p)})_{>0}\\ T\ \text{represented by}\
V({\Bbb A}_f^p),\\
\text{but not by}\ V(\Q_p)\\
T\ \text{represents 1 over}\ \Z_p\\
\endmatrix}
E'_T(i(h_1,\ldots, h_r), 0, \Phi)\tag0.9$$
On the right, we sum over
certain Fourier coefficients of the derivative at 0 of the
Eisenstein series for $Mp_{\Bbb A}$. Our second main result is the
following identity (Corollary~9.4).
\proclaim{Theorem 0.2}
We have
$$\align
F_{d_1, \ldots, d_r} (h_1, \ldots, h_r, \Phi)_p^{\roman{proper}} =\
& c\, W_{d_1}^{5/2}(h_1)\ldots W_{d_r}^{5/2}(h_r) \cdot
{\roman{log}}\, p\cdot
\\
\nass
{}&\cdot\vol(\pr(K))\cdot <{\Cal Z}(d_1, \omega_1),\ldots, {\Cal
Z}(d_r,
\omega_r)>_p^{\roman{proper}}\ ,\tag0.10
\endalign$$
where $c=\frac12 \vol(SO(V'(\R)))$.
\endproclaim
\noindent Unexplained notation may be found in the body of the text. The
identity is proved by unravelling both sides of (0.10), where, for
the right side, we use the decomposition (0.4) and Theorem 0.1. The
identity then reduces to the statement that, for $T\in
{\roman{Sym}}_4(\Z_{(p)})_{>0}$ such that $T$ is not represented by
$V(\Q_p)$ and where $T$ represents 1 over $\Z_p$, we have
$$\left[ ({\roman{log}}\, p)^{-1}\cdot
\frac{W'_{T,p}(e,0,\Phi_p)}{W_{T,p}(e,0,\Phi'_p)}\right]
\left[ {\roman{vol}}(K)^{-1}\cdot I_{T,f}(\varphi_f^{(p)})\right]
=<{\Cal Z}(T,\omega)>_p^{\roman{proper}}\tag0.11$$
Here, in the first factor on the left, there appears a quotient of
the {\it derivative}\/ at 0 of a certain Whittaker function for the
quadratic space $V(\Q_p)$ by the {\it value}\/ at 0 of a Whittaker
function for a twist $V'(\Q_p)$, and, in the second factor, a
Fourier coefficient of a theta integral. It turns out that the
first factor equals the multiplicity $e(\xi)$ of any point $\xi\in
{\Cal Z}(T,\omega)$ (which is constant), while the second factor is
equal to the number of points in ${\Cal Z}(T,\omega)$. For the
multiplicity $e(\xi)$, the calculation can be reduced to a problem
on one-dimensional formal groups of height 2 which has been solved
by Gross and Keating \cite{\grosskeating}. For the calculation of
the Whittaker functions we use the results of Kitaoka
\cite{\kitaokatwo} on local representation densities. It should be
pointed out that we are using here the length of the local ring
$\Cal O_{\Cal Z,\xi}$ as the multiplicity of a point $\xi$, whereas
the sophisticated definition would also involve Tor-terms. It is a
fundamental question whether these correction terms vanish. This
question we have to leave open.
\hfill\break
In summary, we may say that Theorem 0.2 is proved by explicitly
computing both sides of (0.10) and comparing them. It would of
course be highly desirable to find a more direct connection between
the analytic side and the algebro-geometric side of this identity.

We now give an overview of the structure of this paper. In section
1, we introduce the Shimura variety and formulate the moduli problem
solved by ${\Cal M}$. Our special cycles are introduced in section
2. We define the fundamental matrix in section 3 and isolate there
the part of ${\Cal Z}$ lying purely in characteristic $p$. It is
clear from the above description that to proceed further we need a
thorough understanding of the supersingular locus of ${\Cal
M}\times_{\roman{Spec}\, \Z_{(p)}}{\roman{Spec}}\, {\Bbb F}_p$.
This is essentially due to Moret-Bailly \cite{\moretbailly} and Oort
\cite{\oort}.
In section 4, we give a presentation of their results in terms of
Dieudonn\'e theory, better suited for our needs. A similar
presentation was independently given by Kaiser \cite{\kaiser} for a
different purpose. The heart of the paper is section 5. In it we
determine the space of special endomorphisms of certain Dieudonn\'e
modules and deduce the characterization of isolated intersection
points (Theorems~5.11, 5.12 and~5.14).
Here again the exceptional isomorphism plays a vital role.
In section 6, we explain the reduction of the calculation of
$e(\xi)$ to the result of Gross and Keating, and, in section 7, we
explain how to count the number of isolated points. Section 8 is a
review of the Fourier coefficients of Siegel Eisenstein Series. In
section 9, we bring everything together and prove the identity (0.10)
above. In section 10, we review some results of Kitaoka and show how
they can be used to prove the formulas on Whittaker functions
needed in section 9. Finally there is an appendix containing some
facts on Clifford algebras in our special situation.
\hfill\break
In conclusion we wish to thank A.\ Genestier for very useful
discussions on our special cycles which helped us to correct some
misconceptions we had about them. We also thank Th.\ Zink for
helpful remarks. We thank the NSF and the DFG for their support. S.
K. would like to express his appreciation for the hospitality of
the Univ.~Wuppertal and the Univ.~of Cologne during January 1995
and May and June of 1997 respectively. Finally, M.\ R.\ is very
grateful to the Math Department of the University of Maryland for
inviting him and making his stay in Washington a memorable
pleasure.

\vfill\eject
\subheading{\Sec1. The Shimura variety}

In this section, we review the construction of the  Siegel 3-folds
associated to indefinite quaternion algebras over $\Q$, and the
corresponding moduli problem. The use of the Clifford algebra is
modeled on \cite{\shih}. We refer to the appendix for some facts on
those Clifford algebras that will be relevant for our purposes.

Let $B$ be an indefinite quaternion algebra over $\Q$, let
$C=M_2(B)$, with involution $x'={}^tx^\iota$, and let
$$V=\{\ x\in C\ \mid \ x'=x \text{ and } {\roman{tr}}(x)=0\ \}.\tag1.1$$
We define a quadratic form $q$ on $V$ by setting $x^2=q(x)\cdot
1_2\in M_2(B)$, cf.\ Appendix, A.3. Since $B$ is indefinite, the
signature of $(V,q)$ is $(3,2)$, cf.\ Appendix, A.6,  and the Witt
index of $V$ over $\Q$ is $2$ if $B=M_2(\Q)$ and $1$ if $B$ is a
division algebra, cf.\ Appendix, A.3. Let $C(V)$ be the Clifford
algebra of the quadratic space $(V,q)$. Since, for $x\in V\subset
C$, $x^2=q(x)$, there is a natural algebra homomorphism $C(V)\lra
C$ extending the inclusion of $V$ into $C$. The restriction of this
map to the even Clifford algebra $C^+(V)$ induces an isomorphism
$$C^+(V)\simeq C.\tag1.2$$

Let
$$G=GSpin(V)=\{g\in C^\times\mid gg'=\nu(g)\ \},\tag1.3$$
cf.\ Appendix, A.3, so that $G$ is a twisted form over $\Q$ of
$GSp_4$, cf.\ Appendix, A.2. The group $G$ acts on $V\subset C$ by
conjugation and this action yields an exact sequence
$$1\lra Z\lra G\lra SO(V)\lra 1,\tag1.4$$
where $Z$ is the center of $G$.

Let ${\Cal D}$ be the space of oriented negative $2$-planes in
$V(\R)$. This space has two connected components and the group
$G(\R)$ acts transitively on it, via its action on $V(\R)$. For an
oriented $2$-plane $z\in {\Cal D}$, let $z_1$, $z_2\in z$ be a
properly oriented basis such that the restriction of the quadratic
form $q$ from $V(\R)$ to $z$ has matrix $-1_2$ for the basis $z_1$,
$z_2$. Let $j_z = z_1z_2\in C(\R)$. Viewing $j_z$ as the image of
the element $z_1z_2\in C(V(\R))$, the Clifford algebra of $V(\R)$,
and recalling the commutative diagram of section A.3 of the
Appendix, we see that $j_z'=-j_z$ and that $j_z^2=-z_1^2z_2^2=-1$.
Hence, $j_zj_z'=1$ and so, $j_z\in G(\R)$. There is an isomorphism
of algebras over $\R$,
$$\C\isoarrow C^+(z)\qquad i\mapsto z_1z_2,\tag1.5$$
where $C^+(z)$ is the even Clifford algebra of the real $2$-plane $z$.
The composition of this map with the map
$$C^+(z)\subset C^+(V(\R)) \isoarrow C(\R)=M_2(B(\R))\tag1.6$$
induces a morphism, defined over $\R$, $h_z:\Bbb S \lra G$, where
$\Bbb S = R_{\C/\R}\Bbb G_m$, as usual. Note that $h_z(i) = j_z$.
The space ${\Cal D}$ can thus be viewed as the space of conjugacy
classes of such maps under the action of the group $G(\R)$. The
data $(G,{\Cal D})$ or $(G,h_z)$ defines a Shimura variety
$Sh(G,{\Cal D})$, \cite{\deligne}, \cite{\delignetwo}, whose
canonical model is defined over $\Q$. Note that ${\Cal D}$is
isomorphic to two copies of the Siegel space of genus $2$, and, if
$B=M_2(\Q)$, $Sh(G,{\Cal D})$ is just the Siegel modular variety of
genus $2$.

Since $G$ satisfies the Hasse principle, the Shimura variety
represents a certain moduli problem over $(Sch/\Q)$,
\cite{\kottwitz}. To define this we must introduce more notation.

Fix a maximal order $\Cal O_B$ in $B$ such that $\Cal O_B^\iota =
\Cal O_B$, and let $\Cal O_C=M_2(\Cal O_B)$. Let $D(B)$ be the
product of the primes $p$ at which $B_p$ is division, and, as in
\cite{\boutotcarayol}, choose $\tau\in B^\times$ such that
$\tau^\iota=-\tau$, $\tau^2=-D(B)$, and $\tau \Cal O_B
\tau^{-1}=\Cal O_B$.  By section~A.7 of the Appendix, the map
$x\mapsto x^*=\tau x^\iota \tau^{-1}$ is a positive involution of
$B$ preserving $\Cal O_B$. Also, for
$$\a=\pmatrix \tau&{}\\{}&\tau\endpmatrix\in M_2(B),\tag1.7$$
$\a'=-\a$ and $x^*=\a x' \a^{-1} = \a^{-1}x'\a$ is a positive involution of
$C$,
preserving $\Cal O_C$.

Let $U=\Cal O_C$, viewed as a module for $\Cal O_C$ under both left and right
multiplication. Define an alternating form:
$$<\ ,\ >\ :U\times U\lra \Z\tag1.8$$
by
$$<x,y>\ = {\roman{tr}}(y'\a^{-1} x).\tag1.9$$
Then
$$<cx,y>\ = {\roman{tr}}(y'\a^{-1}cx) = {\roman{tr}}(y'\a^{-1}c\a\a^{-1} x) =\
<x,c^*y>,\tag1.10$$
and
$$<xc,y>\ = {\roman{tr}}(y'\a^{-1}xc) = {\roman{tr}}(cy'\a^{-1} x) =\
<x,yc'>.\tag1.11$$
Thus, if $g\in G$,
$$<xg,yg> = \nu(g) <x,y>,\tag1.12$$
and, in particular, for $z\in {\Cal D}$,
$$<x j_z,y j_z>\ =\ <x,y>.\tag1.13$$

We fix a compact open subgroup $K\subset G(\A_f)$. The functor
$M_K$ associates to $S\in (\Sch/\Q)$ the set of quadruples,
$(A,\iota,\lambda,\bar\eta)$, up to isomorphism, where
\roster
\item"{(i)}" $A$ is an abelian scheme over $S$, up to isogeny,
\item"{(ii)}" $\iota:C\lra End^0(A)$ is a homomorphism such that
$${\roman{det}}(\iota(c);\Lie(A)) = N^o(c)^2,$$
where $N^o(c)$ is the reduced norm on $C$.
\item"{(iii)}" $\lambda$ is a $\Q$-class of polarizations on $A$
which induce the involution $*$ on $C$:
$$\lambda\circ\widehat{\iota(c)}\circ\lambda^{-1} = \iota(c^*).$$
\item"{(iv)}" $\bar\eta$ is a $K$-class of isomorphisms
$$\eta:\hat{V}(A) \isoarrow U\tt \A_f$$
which are $C$-linear (for the left module structure on $U$) and
respect the symplectic forms on both sides up to a constant in $\A_f^\times$.
Here
$$\hat{V}(A) = \prod_{\ell}T_\ell(A)\tt\Q.$$
\endroster
For the precise meaning of the datum (iv) we refer to
\cite{\kottwitz}, p.\ 390. In particular, if $S={\roman{Spec}}\, k$
is the spectrum of a field, the $K$-class $\overline{\eta}$ is
supposed to be stable under the action of the Galois group
${\roman{Gal}}(\overline{k}/k)$ where $\overline{k}$ is the
algebraic closure used to form the Tate module of $A$.

Note that the abelian scheme $A$ will have relative dimension $8$
over $S$.

\proclaim{Proposition 1.1}
For $K$ neat this moduli problem is representable by a smooth quasi-projective
scheme $M_K$ over $\Q$ and
$$M_K(\C) \simeq Sh(G,{\Cal D})(\C).$$
\endproclaim
\demo{Proof} For the representability, see \cite{\kottwitz}. We prove the last
assertion in  detail, since the conventions involved will be used
later.

For $\tau\in B^\times$, as above, let
$$\tau_0= D(B)^{-\frac12} \tau\in B^\times(\R),\tag1.14$$
so that $\tau_0^2=-1$.
Choose $\beta\in B^\times$ such that
$$\beta\tau=-\tau\beta,\qquad\text{ and }\qquad \beta^\iota=-\beta.\tag1.15$$
Since $B$ is indefinite, $\beta^2>0,$ and we can set
$$\beta_0=(\beta^2)^{-\frac12}\beta\in B^\times(\R),\tag1.16$$
so that $\beta_0^2=1$. The vectors
$$\pmatrix {}&\beta_0\\ -\beta_0&{}\endpmatrix, \qquad\text{ and } \qquad
\pmatrix {}&\tau_0\beta_0\\-\tau_0\beta_0&{}\endpmatrix \in V(\R)\tag1.17$$
form a standard basis of an oriented negative $2$-plane $z_0\in {\Cal D}$,
and
$$j_{z_0}=\pmatrix {}&\beta_0\\ -\beta_0&{}\endpmatrix
\pmatrix {}&\tau_0\beta_0\\-\tau_0\beta_0&{}\endpmatrix = \pmatrix
\tau_0&{}\\{}&\tau_0
\endpmatrix
=D(B)^{-\frac12} \a =:\a_0.\tag1.18$$

\proclaim{Lemma 1.2} For any $z\in {\Cal D}$,
$$<xj_z,y>\ =\ < y j_z,x>,$$
and, for $x\in U(\R)$, $x\ne 0$,
$$<x j_z,x>\ >0,$$
if $z$ lies in the same connected component of ${\Cal D}$ as $z_0$, and
$$<x j_z,x>\ <0,$$
if $z$ and $z_0$ lie in opposite components.
\endproclaim
\demo{Proof} For the first assertion:
$$<x j_z,y>\ =\ -<x,y j_z>\ =\ <y j_z,x>.\tag1.19$$
For the second, write $z= g z_0$ for $g\in G(\R)$, so that
$$j_z= g j_{z_0} g^{-1}= g \a_0 g^{-1}.\tag1.20$$
Then, we have
$$\align
<x j_z,x>&=\ <x g \a_0 g^{-1}, x>\\
{}&=\nu(g)^{-1} <x g \a_0,x g > \\
{}&= \nu(g)^{-1}{\roman{tr}}((xg)'\a^{-1} xg\a_0)\tag1.21\\
{}&= \nu(g)^{-1} D(B)^{-\frac12}\, {\roman{tr}}(\a (xg)' \a^{-1} (xg))\\
{}&= \nu(g)^{-1} D(B)^{-\frac12}\,{\roman{tr}}((xg)^*(xg)).
\endalign$$
Since $x\mapsto x^*$ is a positive involution, this gives the claim.
\qed\enddemo

Let ${\Cal D}^+$ be the connected component of ${\Cal D}$ containing
$z_0$ and ${\Cal D}^-$ the connected
component of ${\Cal D}$ not containing $z_0$. Then,
for any $z\in {\Cal D}^\pm$, we obtain a (principally) polarized abelian
variety over $\C$,
$$A_z=(U(\R), j_z, U(\Z),\pm<\ ,\ >)\tag1.22$$
with $\dim A_z=8$ and with an action
$$\iota: \Cal O_C\hookrightarrow \End(A_z).\tag1.23$$
Note that
$$\hat{V}(A_z) = U(\hat{\Bbb Z})\tt\Q=U(\A_f).\tag1.24$$
If
$$\gamma\in \Gamma =\{\ g\in G(\Q)^+\mid U(\Z) g = U(\Z)\ \},\tag1.25$$
then right multiplication by $\gamma^{-1}$ induces an isomorphism
$$A_z\isoarrow A_{\gamma z}.\tag1.26$$
Thus $\Gamma\back {\Cal D}^+$ parametrizes such principally
polarized abelian varieties, up to isomorphism.

More generally, to $(z,g)\in {\Cal D}\times G(\A_f)$, we associate the
collection $(A,\iota,\lambda,\bar\eta)$ defined by:
\roster
\item"{$\bullet$}" $(A,\iota) = (A_z,\iota)$, where $A_z$ is taken up to
isogeny.
\item"{$\bullet$}" $\lambda$ is the $\Q$-class of polarizations determined by
$<\ ,\ >$.
\item"{$\bullet$}" $\bar\eta$ is the $K$-class containing the isomorphism:
$$\hat{V}(A_z) = U(\A_f) \ \overset{r(g)}\to{\isoarrow}\ U(\A_f).$$
\endroster
Note that, if $\gamma\in G(\Q)$ and $k\in K$, then $(\gamma z, \gamma g k)$
defines a collection isomorphic to that defined by $(z,g)$, via the
element of $\Hom^0(A_z,A_{\gamma z})$ given on $U(\R)$ by right
multiplication by $\gamma^{-1}$. The map
$$G(\Q)(z,g)K \mapsto (A,\iota,\lambda,\bar\eta)/ \sim\tag1.27$$
yields the isomorphism
$$G(\Q)\back {\Cal D}\times G(\A_f) /K \isoarrow M_K(\C).\qquad\qed\tag1.28$$
\enddemo

\define\uzp{U_{\Z_p}}

We now turn to the construction of a $p$-integral model. Fix a
prime $p$ such that $p\nmid D(B)$, so that $C\tt \Q_p\simeq
M_4(\Q_p)$. Let $\Cal O_C$ be the maximal order chosen above,
and note that the maximal order $\Cal O_C\tt\Z_p$ in $C\tt\Q_p$
is the stabilizer of the lattice $U_{\Z_p}=U\tt\Z_p$ in $U\tt\Q_p$
under both right and left multiplication. The choice of $\tau$ made
before (1.7) ensures that $<\ ,\ >$ defines a perfect pairing
$$<\ ,\ >:\uzp\times \uzp \lra \Z_p.\tag1.29$$
Let $K_p$ be the stabilizer of $\uzp$ in $G(\Q_p)$, acting on  $U_{\Q_p}$ via
right multiplication. Let $K^p\subset G(\A_f^p)$ be compact open, and
take $K=K_p\cdot K^p$.

We now want to formulate a moduli problem over $(\Sch/\Z_{(p)})$
which extends the previous one. The functor $\Cal M_{K^p}$ associates
to $S\in (\Sch/\Z_{(p)})$ the set of isomorphism classes of quadruples
$(A,\iota,\lambda,\bar\eta^p)$ where
{\parindent=12pt
\roster
\item"{(i)}" $A$ is an abelian scheme over $S$, up to prime to $p$
isogeny
\item"{(ii)}" $\iota: \Cal O_C\tt \Z_{(p)}\lra \End(A)\tt\Z_{(p)}$ is a
homomorphism such that, for $c\in \Cal O_C$,
$$\det(\iota(c);\Lie(A)) = N^o(c)^2,$$
where $N^o$ is the reduced norm on $C$.
\item"{(iii)}" $\lambda$ is a $\Z_{(p)}^\times$\snug-class of isomorphisms
$A\lra \hat A$ such that $n\lambda$, for a suitable natural number $n$,
is induced by an ample line bundle on $A$.
\item"{(iv)}" $\bar\eta^p$ is a $K^p$\snug-class of $\Cal O_C$-linear
isomorphisms (in the sense of Kottwitz)
$$\eta^p:\hat V^p(A) \isoarrow U\tt\A_f^p,$$
which respects the symplectic form on both sides up to a constant
in $(\A_f^p)^\times$. Here
$$\hat{V}^p(A) = \prod_{\ell\ne p}T_\ell(A)\tt\Q.$$
\endroster}

In the determinant condition above, the equality is meant as an identity of
polynomial
functions. In the case at hand, it simply says $\dim A=8$.
\proclaim{Proposition 1.3}
For $K^p$ neat the above moduli problem is representable by a smooth
quasiprojective scheme ${\Cal M}_{K^p}$ over ${\roman{Spec}}\ {\Bbb Z}_{(p)}$.
Its
generic fibre can be canonically identified with $M_K$,
$${\Cal M}_{K^p}\times_{{\roman{Spec}}\ {\Bbb Z}_{(p)}} {\roman{Spec}}\ {\Bbb
Q}
=M_K\ \ .$$
\endproclaim
Let us briefly explain the last identification on geometric points. Let
$S$ be the spectrum of an algebraically closed field of characteristic $0$.
Let us consider $(A,\iota,\lambda,\bar\eta^p)\in \Cal M_{K^p}(S)$. Then the
$p$-adic Tate module $T_p(A)$ is equipped with a perfect symplectic form,
unique up to
scaling by $\Z_p^\times$ and hence there is an $\Cal O_C\tt\Z_p$\snug-linear
isomorphism
$$\eta_p:T_p(A) \isoarrow \uzp,$$
which respects the symplectic forms up to $\Z_p^\times$. The set of
such $\eta_p$'s form a single orbit for $K_p$, which acts via right
multiplication in $\uzp$. Hence, from
$(A,\iota,\lambda,\bar\eta^p)$, we obtain an object
$(A\tt\Q,\iota\tt\Q,\lambda\tt\Q,\bar\eta^p\cdot\bar \eta_p)$ of
$M_K(S)$. Passage in the other direction is similar. For example,
in the isogeny class $A$ and for $\eta\in\bar\eta$, there is an
abelian variety $B$, unique up to prime to $p$ isogeny, such that
$\eta_p(T_p(B)) = \uzp$.

The above proposition tells us that, when $K=K^p\cdot K_p$, as
above, then $\Cal M_{K^p}$ provides us with a smooth model of
$Sh(G,D)_K$ over $\Z_{(p)}$. From now on, we will use the same
notation for both moduli problems, if this does not cause
confusion.

\comment
%%%%%%%%%%%%%% Shifted material %%%%%%%%%%%%
%% This material should be returned to section~1.8 of the Appendix
%% OR kept in this section!!?
%%%%%%%%%%%%%%

For $\tau\in B^\times$, as above, let
$$\tau_0= D(B)^{-\frac12} \tau\in B^\times(\R),$$
so that $\tau_0^2=-1$.
Choose $\eta\in B^\times$ such that
$$\eta\tau=-\tau\eta,\qquad\text{ and }\qquad \eta^\iota=-\eta.$$
Since $B$ is indefinite, $\eta^2>0,$ and we can set
$$\eta_0=(\eta^2)^{-\frac12}\eta\in B^\times(\R),$$
so that $\eta_0^2=1$. The vectors
$$\pmatrix {}&\eta_0\\ -\eta_0&{}\endpmatrix, \qquad\text{ and } \qquad
\pmatrix {}&\tau_0\eta_0\\-\tau_0\eta_0&{}\endpmatrix \in V(\R)$$
form a standard basis of an oriented negative $2$-plane $z_0\in D$,
and
$$j_{z_0}=\pmatrix {}&\eta_0\\ -\eta_0&{}\endpmatrix
\pmatrix {}&\tau_0\eta_0\\-\tau_0\eta_0&{}\endpmatrix = \pmatrix
\tau_0&{}\\{}&\tau_0
\endpmatrix
=D(B)^{-\frac12} \a =:\a_0.$$

\proclaim{Lemma} For any $z\in D$,
$$<xj_z,y>\ =\ < y j_z,x>,$$
and, for $x\in U(\R)$, $x\ne 0$,
$$<x j_z,x>\ >0,$$
if $z$ lies in the same connected component of $D$ as $z_0$, and
$$<x j_z,x>\ <0,$$
if $z$ and $z_0$ lie in opposite components.
\endproclaim
\demo{Proof} For the first assertion:
$$<x j_z,y>\ =\ -<x,y j_z>\ =\ <y j_z,x>.$$
For the second, write $z= g z_0$ for $g\in G(\R)$, so that
$$j_z= g j_{z_0} g^{-1}= g \a_0 g^{-1}.$$
Then, we have
$$\align
<x j_z,x>&=\ <x g \a_0 g^{-1}, x>\\
{}&=\nu(g)^{-1} <x g \a_0,x g > \\
{}&= \nu(g)^{-1}{\roman{tr}}((xg)'\a^{-1} xg\a_0)\\
{}&= \nu(g)^{-1} D(B)^{-\frac12}\, {\roman{tr}}(\a (xg)' \a^{-1} (xg))\\
{}&= \nu(g)^{-1} D(B)^{-\frac12}\,{\roman{tr}}((xg)^*(xg)).
\endalign$$
Since $x\mapsto x^*$ is a positive involution, this gives the claim.
\qed\enddemo

Let $D^+$ be the connected component of $D$ containing $z_0$. Then,
for any $z\in D^+$, we obtain a (principally) polarized abelian variety
$$A_z=(U(\R), j_z, U(\Z), <\ ,\ >)$$
with $\dim A_z=8$ and with an action
$$\iota: \Cal O_C\hookrightarrow \End(A_z).$$
If
$$\gamma\in \Gamma =\{\ g\in G(\Q)\mid U(\Z) g = U(\Z)\ \},$$
then, right multiplication by $\gamma^{-1}$ induces an isomorphism
$$A_z\isoarrow A_{\gamma z}.$$

\endcomment

\vfill\eject
\subheading{\Sec2. Special cycles}

In this section we give a modular definition of the special cycles in
$Sh(G,{\Cal D})$, which were defined analytically in \cite{\cycles}. We then
explain the relation between the two definitions.

Recall that the quadratic form on the space $V\subset C=M_2(B)$
was defined by $x^2=q(x)\cdot 1_2$. Let
$$(x,y) = q(x+y)-q(x)-q(y)\tag2.1$$
be the corresponding bilinear form, so that $q(x)=\frac12 (x,x)$.
If $x=(x_1,x_2,\dots,x_n)\in V^n({\Bbb Q})$, we let
$$q(x) = \frac12 ((x_i,x_j))_{i,j}\in {\roman{Sym}}_n(\Q).\tag2.2$$
This defines a quadratic map $q:V^n\lra {\roman{Sym}}_n$.

Fix a positive integer $n$. For $d\in {\roman{Sym}}_n(\Q)$ a
symmetric rational matrix, let
$$\O_d = \{\ x\in V^n\ \mid\ q(x) = d\ \}\tag2.3$$
be the corresponding hyperboloid. The group $G$ acts diagonally on
$V^n$ and preserves $\O_d$.

Cycles in $Sh(G,{\Cal D})$ were defined analytically in
\cite{\cycles} as follows. For $x\in \O_d(\Q)$, let $<x>\ \subset
V$ be the $\Q$-subspace spanned by the components of $x$, and let
$V_x
=\  <x>^\perp$ be its orthogonal complement. Let ${\Cal D}_x$
denote the space of oriented negative $2$-planes in $V_x(\R)$, and
let $G_x$ be the pointwise stabilizer of $< x>$ in $G$. Note that
$G_x
\simeq GSpin(V_x)$, and that ${\Cal D}_x\subset {\Cal D}$.
Moreover, for $z\in {\Cal D}_x$, the homomorphism $h_z$ factors
through $G_x(\R)$. Thus there is a natural morphism of Shimura
varieties, rational over $\Q$,
$$Sh(G_x,{\Cal D}_x) \lra Sh(G,{\Cal D}).\tag2.4$$
If the space $<x>$ is not positive-definite, then ${\Cal
D}_x=\emptyset$. If $<x>$ is positive-definite of dimension $r$
then $d$ is positive semi-definite of rank $r$,
${\roman{sig}}(V_x)= (3-r, 2)$ and ${\Cal D}_x$ has codimension $r$
in ${\Cal D}$. Hence the previous construction is only interesting
when $d$ is positive semi-definite and even only when $d$ is
positive definite with $n\leq 3$.

For a fixed compact open subgroup $K\subset G(\A_f)$ and for $h\in
G(\A_f)$, there is a cycle, namely the image of the map
$$Z(x,h;K): G_x(\Q)\back {\Cal D}_x\times G_x(\A_f)/ (G_x(\A_f)\cap hKh^{-1})
\lra G(\Q)\back {\Cal D}\times G(\A_f)/K\tag2.5$$
given by $(z,g)\mapsto (z,gh)$. This map is finite and generically
injective, hence the cycle image is taken with multiplicity 1. This
cycle of codimension $r={\roman{rk}}(d)$ on $Sh(G,{\Cal D})_K$ is
rational over $\Q$.

Assume that $\O_d(\Q)\ne \phi$
and fix $x_0\in \O_d(\Q)$.
Let $\ph\in S(V(\A_f)^n)^K$ be a Schwartz function which is $K$-invariant,
and write
$$\text{supp}(\ph)\cap \O_d(\A_f) = \coprod_{r} K h_r^{-1}x_0\tag2.6$$
for elements $h_r\in G(\A_f)$. Then define the {\bf weighted cycle}:
$$Z(d,\ph;K) = \sum_r \ph(h_r^{-1}x_0)\cdot Z(x_0,h_r;K).\tag2.7$$
This cycle is independent of the choice of $x_0$ and  of the orbit
representatives $h_r$. It is a (weighted linear combination of)
cycle(s) of codimension $r={\roman{rk}}(d)$ on $Sh(G,{\Cal D})_K$
and is rational over $\Q$.

\par\noindent
If $\varphi$ is the characteristic function of a $K$-invariant
compact open subset $\omega$ of $V({\Bbb A}_f)^n$, then
$Z(d,\omega;K)= Z(d,\varphi;K)$ can be considered as a disjoint
union of maps (2.5), or as the union of the images of these maps.

We introduce the following definition, which will play a key
role throughout the paper.
\proclaim{Definition 2.1}
Let $(A,\iota,\lambda,\overline{\eta})\in M_K(S)$ . A {\bf special
endomorphism of} $(A,\iota, \lambda, \overline{\eta})$ is an
element $j\in {\roman{End}}_S^0(A,\iota)$ which satisfies
$$j^*=j\ \ \text{and}\ \ {\roman{tr}}^0(j)=0\ \ .\tag2.8$$
\endproclaim

Here $*$ denotes the Rosati involution of $\lambda$. Also note that
${\roman{End}}^0(A,\iota)$ is a finite-dimensional semisimple
${\Bbb Q}$-algebra, so that the reduced trace appearing here makes
sense. Indeed, this is well-known when $S$ is the spectrum of a
field. The case when $S$ is irreducible follows by reduction to its
generic point, and the general case follows by considering the
irreducible components of $S$.

\proclaim{Lemma 2.2}
Let $j$ be a special endomorphism of $(A,\iota,
\lambda,\overline{\eta})\in M_K(S)$, where $S$ is connected. Then
$$j^2=q(j)\cdot {\roman{id}}\ \ ,\tag2.9$$
with $q(j)\in {\Bbb Q}$.
\endproclaim

\demo{Proof}
Again we may reduce first to the case where $S$ is irreducible and
then to the case when $S$ is the spectrum of a field. However for
$\eta\in\overline{\eta}$ let $x=\eta^{*}(j)\in
{\roman{End}}_C(U({\Bbb A}_f))= C({\Bbb A}_f)$. Under the last
identification the adjoint involution $*$ w.r.t.\ $<\ ,\ >$
corresponds to the involution $'$ on $C({\Bbb A})$, cf.\ (1.11).
Hence $x$ lies in $V({\Bbb A}_f)$ and the assertion follows, cf.\
appendix A.3.\qed
\enddemo

The previous Lemma justifies the following definitions. Let $S$ be
a connected scheme and $\xi= (A,\iota, \lambda, \overline{\eta})\in
M_K(S)$ be an $S$-valued point of $M_K$. Let
$$C_{\xi}^0={\roman{End}}^0_S(A,\iota)^{\roman{op}}\tag2.10$$
and
$$V_{\xi}^0=\{ x\in C_{\xi}^0\mid x^*=x\ \text{and}\
{\roman{tr}}^0(x)= 0\}\ \ .\tag2.11$$ Then $V_{\xi}^0$ is the
finite-dimensional ${\Bbb Q}$-vector space of special endomorphisms
with quadratic form
$$q_{\xi}:V_{\xi}^0\longrightarrow {\Bbb Q}\tag2.12$$
given by $x^2=q_{\xi}(x)\cdot {\roman{id}}_A$. By the universal
property of the Clifford algebra of $(V_{\xi}^0, q_{\xi})$ there is
a natural homomorphism
$$C(V_{\xi}^0, q_{\xi})\longrightarrow C_{\xi}^0\ \ .\tag2.13$$
This structure is compatible with specialization. If $S'\subset S$
is a connected closed subscheme, let $\xi'\in M_K(S')$ be the
restriction of $\xi$. Then we have a homomorphism of ${\Bbb
Q}$-algebras
$$C_{\xi}^0 ={\roman{End}}_S^0(A,\iota)^{\roman{op}}\hookrightarrow
{\roman{End}}_{S'}(A,\iota)^{\roman{op}} =C_{\xi'}^0\tag2.14$$
inducing a map
$$V_{\xi}^0\hookrightarrow V_{\xi'}^0$$
of quadratic spaces.

Let us spell out these concepts in the classical case.
\proclaim{Lemma 2.3}Let $\xi\in M_K({\Bbb C})$ with parameter $(z,g)$ in
$Sh(G,{\Cal D})_K$.
Let $A_z^{\roman{top}}= U({\Bbb R})/ U({\Bbb Z})$ be the real torus
underlying $A_z$.
\par\noindent
(i)
$$C(\Q)\isoarrow \End^0(A^{\text{top}}_z,\iota)^{\text{op}},
\qquad y\mapsto r(y),$$
where $r(y)$ denotes the action of $y\in C(\Q)$ on $U(\R)\supset U(\Q)$
by right multiplication. Moreover,
$r(y)^* =r(y')$.\hfill\break
(ii) $$C^0_\xi \simeq Cent_{C(\Q)}(j_z),\qquad\text{ and }\qquad V^0_\xi
 \simeq \{\  x\in V(\Q)\mid x j_z = j_z x\ \}.$$
(iii) $$Cent_{C(\R)}(j_z) \cap V(\R) = z^\perp.$$
In particular,
$$V_\xi^0 = V(\Q)\cap z^\perp,$$
and so $0\le \dim_\Q V_\xi^0\le 3$.
\endproclaim

\demo{Proof} The first two assertions are obvious by (1.11).
To prove the last assertion let $z_1, z_2\in z$ be a properly
oriented basis such that the restriction of the quadratic form $q$
to $z$ has matrix $-1_2$ in terms of this basis. Let $v\in V(\Bbb
R)$ with $(v,z_i)=a_i$, $i=1,2$. Then
$$v\cdot j_z= v\cdot (z_1\cdot z_2)= z_1z_2v -a_2z_1+a_1z_2=j_zv-
a_2z_1+a_1z_2.$$
Hence $v\in {\roman{Cent}}_{C({\Bbb R})}(j_z)$ iff $a_1=a_2=0$,
i.e.\ iff $v\in z^\bot$.\qed
\enddemo

Let us return to the abstract situation.

\proclaim{Lemma 2.4}
Let $\xi=(A, \iota, \lambda, \overline{\eta})\in M_K(S)$ be a point
with values in a connected scheme $S$. The quadratic space
$V_{\xi}^0$ is positive-definite.
\endproclaim
\demo{Proof}
We may assume that $S$ is the spectrum of a field. The assertion
follows from the positivity of the Rosati involution, since
$$q(x)\cdot {\roman{id}}_A =x^2=x\cdot x^*\ \ ,\ \ x\in V_{\xi}^0\ \ .
\qed$$
\enddemo

\par
We next give a modular definition of the cycles introduced above.
We take here the point of view that a cycle is given by a finite
unramified morphism into the ambient scheme. Let $K\subset G(\A_f)$
be a compact open subgroup, and let $\omega\subset V({\Bbb A}_f)^n$
be a $K$-invariant compact open subset. Consider the functor on
$(\Sch/\Q)$ which associates to a scheme $S$ the set of isomorphism
classes of 5-tuples $(A,\iota,\lambda,\bar\eta;\j)$ where
$(A,\iota,\lambda,\bar\eta)\in M_K(S)$. Here the additional element
$\j = (j_1,\dots,j_n)\in \text{End}^0(A,\iota)^n$ is an
$n$\snug-tuple of special endomorphisms of $A$, satisfying the
following conditions. {\parindent=12pt
\roster
\item"{(2.15)}" For some (and hence for all) $\eta\in\bar\eta$, the
element $\eta^*(\j)\in \text{End}_C(U(\A_f))^n$ lies in $\omega$.
\item"{(2.16)}" $q(\j)=d$.
\endroster}
\noindent
Let us explain the condition (2.15). As in the proof of Lemma 2.2
above, for any $\eta\in\overline{\eta}$
$$x=\eta^*(\j)\in V({\Bbb A}_f)^n\subset C({\Bbb A}_f)^n=
{\roman{End}}_C(U({\Bbb A}_f))^n\ \ .$$ The condition imposes that
$x\in\omega$. If $\eta$ is changed to $r(k)\circ
\eta$, with $k\in K$ and $r(k)\in \text{End}_C(U(\A_f))$ the
endomorphism defined by right multiplication by $k$, then
$$(r(k)\circ \eta)^*(\j) = r(k)\circ \eta^*(\j)\circ r(k)^{-1}.\tag2.17$$
The condition (2.15)  asserts that $\eta^*(\j) = r(x)$ for some
$x\in\omega$. If this is the case, then
$$(r(k)\circ \eta)^*(\j) = r(k)\circ \eta^*(\j)\circ r(k)^{-1} =
r(k^{-1}xk),\tag2.18$$ and $k^{-1}xk\in \omega$. Thus the condition
(2.15) depends only on $\bar\eta$.
\hfill\break
To interpret condition (2.16) we may assume $S$ to be connected.
Let $(\ ,\ )$ be the bilinear form on the space of special
endomorphisms of $(A,\iota,\lambda,\overline{\eta})$ associated to
the quadratic form $q$ of lemma 2.2. Then
$q(\j)=\frac{1}{2}((j_i,j_j))_{i,j}\in {\roman{Sym}}_n({\Bbb Q})$
is defined as in (2.2). The condition (2.16) requires that
$q(\j)=d$.
\proclaim{Proposition 2.5}
The above functor has a coarse moduli scheme ${\Cal Z}(d,\omega)$.
If $K$ is neat, then ${\Cal Z}(d,\omega)$ is a fine moduli scheme
and the forgetful morphism
$${\Cal Z}(d,\omega)\longrightarrow M_K\tag2.19$$
is finite and unramified. Furthermore ${\Cal Z}(d,\omega)({\Bbb C})
=Z(d,\omega, K)$.
\endproclaim
\demo{Proof}
The first statement follows easily from the second. Let us assume
that $K$ is neat. The relative representability of the forgetful
morphism by a morphism of finite type follows in a standard way
from Grothendieck's theory of Hilbert schemes since $M_K$ may be
considered as a moduli scheme of polarized abelian varieties with
additional structure. By the N\'eron universal property the
valuative criterion for properness is satisfied. Since the matrix
$d$ gives the squares $j_i^2$ of the special endomorphisms, the
morphism is quasi-finite and hence finite. The unramifiedness
follows from the rigidity theorem for abelian varieties.
\hfill\break
The last statement is to be interpreted as an equality between the
image of (2.5) and ${\Cal Z}(d,\omega)({\Bbb C})$, and follows
easily from Lemma 2.3 above.
\enddemo

We now assume that $p\nmid 2D(B)$ and that $K= K^p\cdot K_p$ with
$K^p$ neat, as in Proposition 1.3, and we formulate a $p$-integral
version of the previous moduli problem.

Before doing this let us point out that for a point
$\xi=(A,\iota,\lambda, \overline{\eta}^p)\in {\Cal M}_{K^p}(S)$ of
the $p$-integral version of our moduli problem with values in a
connected scheme $S$ we may transpose the concepts above. Hence we
introduce the ${\Bbb Z}_{(p)}$-algebra
$$C_{\xi}={\roman{End}}_S(A,\iota)^{\roman{op}}
\otimes {\Bbb Z}_{(p)}\tag2.20$$
and
$$V_{\xi}=\{ x\in C_{\xi};\ x^*=x\ \text{and}
\ {\roman{tr}}^0(x)= 0\}\ \ .\tag2.21$$
The latter is a ${\Bbb Z}_{(p)}$-module with a ${\Bbb
Z}_{(p)}$-valued positive definite quadratic form. The elements of
$V_{\xi}$ will again be called the special endomorphisms of
$(A,\iota,\lambda, \overline{\eta}^p)$.

Let now again $d\in {\roman{Sym}}_n({\Bbb Q})$. Let
$\omega^p\subset V(\A_f^p)^n$ be a $K^p$-invariant open compact
subset. Then a point of the corresponding moduli problem ${\Cal
Z}(d,\omega^p)$ on a $\Z_{(p)}$-scheme $S$ is an isomorphism class
of 5-tuples $(A,\iota,\lambda,\bar\eta^p;\j)$ where
$(A,\iota,\lambda,\bar\eta^p)$ is an object of $\Cal M_{K^p}(S)$
and where $\j\in \big(\End(A,\iota)\tt\Z_{(p)}\big)^n$ is an
$n$-tuple of special endomorphisms which satisfies (2.16) above
and, in addition,
$$(\eta^p)^*(\j) \in \omega^p.\tag{2.22}$$
These conditions are to be interpreted in the same way as
(2.15)-(2.16) above.
\par
To clarify the relation between the $p$-integral version $\Cal
Z(d,\omega^p)$ and the previous ${\Cal Z}(d,\omega)$, let
$$\omega_p = V(\Z_p)^n\tag2.23$$
where $V(\Z_p) = V(\Q_p)\cap (\Cal O_C\tt\Z_p)$, the intersection
taking place inside of $C\tt\Q_p$. Let
$$\omega =\omega_p\times \omega^p\ \ ,\tag2.24$$
a $K$-invariant open compact subset of $V({\Bbb A}_f)^n$.
\proclaim{Proposition 2.6}
If $K^p$ is neat, the functor ${\Cal Z}(d,\omega^p)$ is
representable by a scheme which maps by a finite unramified
morphism to ${\Cal M}_{K^p}$. Furthermore, there is an
identification
$${\Cal Z}(d,\omega^p)\times_{{\roman{Spec}}\ {\Bbb Z}_{(p)}}
{\roman{Spec}}\, {\Bbb Q} ={\Cal Z}(d,\omega)\ \ .$$
\endproclaim
\noindent
{\bf Remark 2.7:} By Lemma 2.4 the scheme ${\Cal Z}(d,\omega^p)$ is
empty unless $d$ is positive semi-definite. Similarly ${\Cal
Z}(d,\omega^p)=\emptyset$, unless $d\in {\roman{Sym}}_n({\Bbb
Z}_{(p)})$. Note that it may well happen that ${\Cal
Z}(d,\omega^p)$ is non-empty but where both sides of the equality
in Proposition 2.6 are empty. In fact we will later consider cases
in which $d\in {\roman{Sym}}_4({\Bbb Z}_{(p)})$ is positive
definite so that ${\Cal Z}(d,\omega)=\emptyset$ and when ${\Cal
Z}(d,\omega^p)\ne \emptyset$.

>From now on, since we will be interested in the arithmetic
situation, we will simplify our notation by denoting $\omega$ what
is denoted by $\omega^p$ above, i.e.,
$$\omega\subset V({\Bbb A}_f^p)^n\tag2.25$$
is a $K^p$-invariant open compact subset.

\vfill\eject
\subheading{\Sec3. The intersection problem}

We continue to fix $p\nmid 2D(B)$ and a neat open compact subgroup
$K^p\subset G({\Bbb A}_f^p)$ as at the end of section 2. Then
${\Cal M}={\Cal M}_{K^p}$ is a regular noetherian scheme of
dimension 4. We wish to consider the intersection of the cycles
introduced in a modular way in the previous section. Let us set up
our problem in a more precise way.

We fix integers $n_1,\ldots, n_r$ with $1\leq n_i\leq 4$ and with
$n_1+\dots +n_r=4$. For each $i$, we choose $d_i\in
{\roman{Sym}}_{n_i}({\Bbb Q})$ positive definite, and a
$K^p$-invariant open compact subset $\omega_i\subset V({\Bbb
A}_f^p)^{n_i}$. Let
$${\Cal Z} = {\Cal Z}(d_1,\omega_1)\times_{\Cal M}\dots
\times_{\Cal M} {\Cal Z}(d_r, \omega_r)\ \
\tag3.1$$
be the fiber product of the corresponding special cycles.

By what has been said in section 2, since the codimensions of the
generic fibres of our special cycles add up to the arithmetic
dimension of ${\Cal M}_{K^p}$, one might expect that ${\Cal Z}$
consists of finitely many points of characteristic $p$. We will see
that this is in fact quite false, but we will be able to determine
that part of ${\Cal Z}$ which lies purely in characteristic $p$ and
also determine the isolated points of ${\Cal Z}$.

Let $\xi$ be a point of ${\Cal Z}$, with corresponding point
$(A_{\xi}, \iota, \lambda, \overline{\eta}^p)\in {\Cal M}$. We
denote by $C_{\xi}$ and $(V_{\xi}, q_{\xi})$ the ${\Bbb
Z}_{(p)}$-algebra and the quadratic ${\Bbb Z}_{(p)}$-module
associated to $(A_{\xi}, \iota, \lambda, \overline{\eta}^p)$, cf.\
(2.20). The projections ${\Cal Z}\to {\Cal Z}(d_i, \omega_i)$
define $n_i$-tuples of special endomorphisms
$$\j_i\in V_{\xi}^{n_i}\ \ ,\ \ i=1,\ldots, r\ \ .\tag3.2$$
Let
$$T_\xi =\frac12
\pmatrix (\j_1,\j_1)_\xi&
\ldots&
(\j_1,{\bold j}_r)_\xi\\
\vdots&&\vdots\\
({\bold j}_r,\j_1)_\xi&
\ldots&
({\bold j}_r, {\bold j}_r)_\xi\endpmatrix
\in Sym_4(\Z_{(p)}),\tag3.3$$
where $(\ ,\ )_\xi$ is the bilinear form associated to $q_\xi$. Here, as
always, $p\ne2$.
The matrix $T_{\xi}$ is called the {\bf fundamental matrix
associated to the intersection point} $\xi$ of the special cycles
${\Cal Z}(d_1,\omega_1),
\ldots, {\Cal Z}(d_r,\omega_r)$. We note that the blocks on
the diagonal of $T_{\xi}$ are $d_1,\ldots, d_r$. By the results of
section 2, the function $\xi\mapsto T_{\xi}$ is constant on each
connected component of ${\Cal Z}$. Therefore, for $T\in
{\roman{Sym}}_4({\Bbb Z}_{(p)})$ we may introduce
$$\align
{\Cal Z}_T=\ & ({\Cal Z}(d_1,\omega_1)\cap\dots\cap {\Cal
Z}(d_r,\omega_r))_T=\\
\ &
\text{union of the connected components of}\ {\Cal Z}\tag3.4\\
&
\text{consisting of the points $\xi$ with $T_{\xi}=T$}
\endalign$$

We note here the hereditary nature of our construction, given by
$${\Cal Z}(T,\omega_1\times\dots\times \omega_r)=
({\Cal Z}(d_1,\omega_1)\times_{\Cal M}\dots
\times_{\Cal M} {\Cal Z}(d_r,\omega_r))_T\ \ ,\tag3.5$$
valid provided that the blocks on the diagonal of $T$ are
$d_1,\ldots, d_r$. We may therefore write
$$\align
{\Cal Z} &
={\Cal Z}(d_1, \omega_1)\times_{\Cal M}\ldots\times_{\Cal M} {\Cal Z}
(d_r, \omega_r)\\
{}&=\coprod\limits_T {\Cal Z}_T\\
{}&=\coprod\limits_{\matrix
T\in {\roman{Sym}}_4({\Bbb Z}_{(p)})_{\ge 0}\\
{\roman{diag}}(T)=(d_1,\ldots, d_r)
\endmatrix} {\Cal Z}(T,\omega)\ \ .\tag3.6
\endalign$$
Here $\omega=\omega_1\times \ldots\times \omega_r$.

We shall see that the fundamental matrix governs the intersection
behaviour of our special cycles. We first note the following
result.

\proclaim{Proposition 3.1}
Let $\xi\in {\Cal Z}={\Cal Z}(d_1,\omega_1) \times_{\Cal
M}\ldots\times_{\Cal M} {\Cal Z}(d_r,
\omega_r)$ where $\omega_i\subset V({\Bbb A}_f^p)^{n_i}$ and
$d_i\in {\roman{Sym}}_{n_i}({\Bbb Q})$ positive-definite with
$n_1+\ldots +n_r=4$. Suppose that ${\roman{det}}(T_{\xi})\ne 0$.
Then $\xi$ lies in the special fiber of $\Cal Z$, and $\xi$ does
not lie in the closure of any point of ${\Cal Z}$ in the generic
fiber.
\endproclaim

\demo{Proof}
By Lemma 2.4 the assumption on $T_{\xi}$ means that $T_{\xi}\in
{\roman{Sym}}_4({\Bbb Z}_{(p)})$ is positive definite. However, for
a point of ${\Cal M}$ in characteristic zero, the space of special
endomorphisms is contained in a 3-dimensional positive definite
quadratic space.
\qed\enddemo

Next suppose that $\xi\in \Cal M(\bar\F_p)$. In this case, the
standard Honda-Tate results yield information about the
possibilities for $C_\xi^0$. We have $\iota:M_2(B)=C\hookrightarrow
\End^0(A_\xi)$, so that, up to isogeny $A_\xi\simeq A\times A$
where $\dim A = 4$ and there is an embedding $B\hookrightarrow
\End^0(A)$.
\proclaim{Lemma 3.2} Suppose that $p\nmid D(B)$. Then there are no simple
abelian varieties $A_0$ over $\bar\F_p$ with $\dim A_0=2$ or $4$ and with
$B\hookrightarrow \End^0(A_0)$.
\endproclaim
\demo{Proof} If $A_0$ is simple over $\F_q$, then
$E=\End^0(A_0)$ is a central simple algebra over $F = \Q(\pi_{A_0})$, and
$$2 \dim A_0 = [E:F]^{\frac12}\cdot [F:\Q].\tag3.7$$
Here $\pi_{A_0}$ denotes as usual the Frobenius endomorphism. If
$\dim A_0\ge 2$ and $A_0$ remains simple over $\bar\F_p$, then $F$
is a CM field. Suppose that $\dim A_0=2$, so that $[F:\Q]
= 2$ or $4$. The second case is excluded, since then $E=F$ is
commutative. In the first case, $E$ is a division quaternion
algebra over $F$ ramified only at places over $p$. Thus $p$ splits
in $F$ and $inv_v(E)=inv_{\bar{v}}(E)=\frac12$ for $v\mid p$. But
the embedding $B\hookrightarrow E$ yields an isomorphism $B\tt_\Q
F\simeq E$. This is possible only if $p\mid D(B)$ and $F$ splits
$B$ at all other primes.
%\hfill\break

If $\dim A=4$, then $[F:\Q]=2$, $4$ or $8$, and the last case is again
excluded since $E=F$. In the case $[F:\Q]=4$, $E$ is a quaternion algebra
over $F$, ramified only at primes lying over $p$, and $B\tt_\Q F\simeq E$.
This cannot occur if $p\nmid D(B)$. Finally, if $[F:\Q]=2$, then $p$
splits in $F$ and $E$
is a division algebra over $F$ of dimension $16$ with invariants
$\frac14$ and $\frac34$ at the primes over $p$. There is no homomorphism
from a quaternion algebra $B\tt_\Q F$ into such an algebra.
\qed\enddemo

Returning to $A$, and assuming that $p\nmid D(B)$, we see that $A$ cannot
be simple and that any simple factor of $A$ of dimension $1$ or $2$
must occur with multiplicity at least $2$. Thus we have various possibilities
for
$A$, up to isogeny:
\roster
\item"{(3.8.i)}" $A\simeq A_2\times A_2$, with $\dim A_2=2$ simple and
$\End^0(A_2)\simeq F$ for a CM field $F$ with $[F:\Q]=4$ which
splits $B$, i.e., such that $B\tt_\Q F\simeq M_2(F)$. Then,
$\End^0(A) \simeq M_2(F)$,  $C_\xi^0=\End^0(A,\iota) \simeq F$, and
$V_\xi^0=\Q$.
\item"{(3.8.ii)}"  $A\simeq A_2\times A_2$, with $\dim A_2=2$ simple and
$\End^0(A_2)\simeq E$, where $E$ is a quaternion algebra over a CM
field $F$ with $[F:\Q]=2$. More precisely, $p$ splits in $F$ and
$E\simeq \H_p\tt_\Q F$, where $\H_p$ is the quaternion algebra over
$\Q$ ramified at $\infty$ and $p$. Let $B'$ be the quaternion
algebra over $\Q$ whose invariants agree with those of $B$ except
at $\infty$ and $p$. Then $\End^0(A)\simeq M_2(E)$,
$\End^0(A,\iota)\simeq B'\tt_\Q F$, and $V_\xi^0=\{x\in B'\mid
{\roman{tr}}(x)=0\}$. Here note that $B\tt_\Q B'\simeq M_2(\H_p)$
and hence that $(B\tt_\Q F)\tt_F (B'\tt_\Q F) \simeq M_2(E)$.
\item"{(3.8.iii)}" $A\simeq A_0^2\times A_1^2$ where $A_0$ and $A_1$
are non-isogenous ordinary elliptic curves. Then $\End^0(A)\simeq
M_2(F_0)\times M_2(F_1)$ for imaginary quadratic fields $F_0$ and
$F_1$, which split $B$. Then, $\End^0(A,\iota)\simeq F_0\times
F_1$, and $V_\xi^0=\Q$.
\item"{(3.8.iv)}" $A\simeq A_0^4$, for an ordinary elliptic curve $A_0$. Then
$\End^0(A)\simeq M_4(F_0)$ where the imaginary quadratic field
$F_0$ splits $B$, $\End^0(A,\iota)\simeq M_2(F_0)$
and $V_\xi^0 \simeq \{x\in M_2(F_0)\mid
{}^t\bar x =x,\ {\roman{tr}}(x)=0\}$.
\item"{(3.8.v)}" $A\simeq A_0^2\times A_1^2$, where $A_0$ is a supersingular
elliptic curve and $A_1$ is an ordinary elliptic curve.
Then $\End^0(A)\simeq M_2(\H_p)\times M_2(F_1)$, $\End^0(A,\iota)\simeq
B'\times F_1$.
Since the Rosati involution acts on $\End^0(A,\iota)$
by $(b,a)\mapsto (b^\iota,\bar a)$, the
conditions $x^*=x$ and ${\roman{tr}}(x)=0$ force $V_\xi^0\simeq \Q$.
\item"{(3.8.vi)}" $A\simeq A_0^4$, for a supersingular elliptic curve $A_0$.
Then
$\End^0(A)\simeq M_4(\H_p)$, $\End^0(A,\iota)\simeq M_2(B')$, and
$$V_\xi^0=\{x\in M_2(B')\mid x' = {}^tx^\iota = x, \ {\roman{tr}}(x)=0\}= V'.$$
\endroster

For the last identification we are using the proposition in section
A.4 of the Appendix. Indeed, the Rosati involution on
${\roman{End}}^0(A_{\xi})\simeq M_8(H_p)$ is of main type by A.4
and its restriction to $M_2(B)$ of neben type. Hence the involution
on $M_2(B')$ is indeed of main type.

Note that $\dim V_\xi^0\le 3$, with the exception of the
supersingular case (3.8.vi). As a consequence, we have the following:
\proclaim{Proposition 3.3} Let $T\in Sym_4({\Bbb Z}_{(p)})$
and $\omega\subset V({\Bbb A}_f^p)^4$ with corresponding special
cycle ${\Cal Z}(T,\omega)$. If ${\roman{det}}(T)\ne 0$, then the
point set underlying ${\Cal Z}(T,\omega)$ maps to the supersingular
locus of ${\Cal M}\times_{{\roman{Spec}\, {\Bbb Z}_{(p)}
}}{\roman{Spec}}\, {\Bbb F}_p$. In particular, ${\Cal Z}(T,\omega)$
is proper over ${\roman{Spec}}\, {\Bbb Z}_{(p)}$ with support in
the special fibre.
\endproclaim
\demo{Proof}
Indeed the previous results imply that this is true for closed points.
\enddemo

\proclaim{Corollary 3.4}
For $i=1,\ldots, r$, let $d_i\in {\roman{Sym}}_{n_i}({\Bbb Q})$ be
positive definite with $n_1+\ldots +n_r= 4$, and let
$\omega_i\subset V({\Bbb A}_f^p)^{n_i}$ with corresponding special
cycles ${\Cal Z}(d_i,
\omega_i)$. For $T\in {\roman{Sym}}_4({\Bbb Z}_{(p)})$ with diagonal
blocks $d_1,\ldots, d_r$, let ${\Cal Z}_T$ be the union of the
connected components of ${\Cal Z}(d_1, \omega_1)\times_{\Cal
M}\ldots \times_{\Cal M}{\Cal Z}(d_r, \omega_r)$ where the
fundamental matrix has value $T$. If ${\roman{det}}(T)\ne 0$ then
the point set underlying ${\Cal Z}_T$ lies over the supersingular
locus of ${\Cal M}\times_{\roman{Spec}\, {\Bbb
Z}_{(p)}}{\roman{Spec}}\, {\Bbb F}_p$.
\qed
\endproclaim

Having answered these very crude questions on the intersection
behaviour of our special cycles, we are led to ask more precise
questions. Again for $i=1,\ldots, r$ let
$d_i\in{\roman{Sym}}_{n_i}({\Bbb Q})$ be positive-definite with
$n_1+\ldots +n_r=4$ and let $\omega_i\subset V({\Bbb A}_f^p)^{n_i}$
with corresponding cycles ${\Cal Z}(d_1,\omega_1), \ldots, {\Cal
Z}(d_r,
\omega_r)$. We then ask:

a) Under which conditions do the cycles ${\Cal
Z}(d_1,\omega_1),\ldots, {\Cal Z}(d_r,\omega_r)$ intersect
properly? More precisely, can one parametrize the isolated points
of ${\Cal Z}={\Cal Z}(d_1,\omega_1)\times_{\Cal M}\ldots
\times_{\Cal M} {\Cal Z}(d_r, \omega_r)$ and calculate at such an isolated
point $y$,
$$e(y)={\roman{lg}}_{{\Cal O}_{{\Cal Z},y}} ({\Cal O}_{{\Cal Z}, y})
\ \ ?\tag3.9$$

b) Let $Y$ be a connected component of ${\Cal Z}={\Cal
Z}(d_1,\omega_1)\times_{\Cal M}\ldots \times_{\Cal M}{\Cal
Z}(d_r,\omega_r)$ lying over the supersingular locus of ${\Cal
M}\times_{\roman{Spec}\, {\Bbb Z}_{(p)}} {\roman{Spec}}\, {\Bbb
F}_p$. The intersection number along $Y$ is
$$\chi(Y, {\Cal O}_{{\Cal Z}_1}\otimes_{\Cal O_{\Cal M}}^{\Bbb L} \ldots
\otimes_{\Cal O_{\Cal M}}^{\Bbb L}{\Cal O}_{{\Cal Z}_r})\ \ ,\tag3.10$$
cf.\ \cite{\krdrin}, \cite{\serre}. An important question to answer
is when the derived tensor product here can be replaced by an
ordinary tensor product, i.e.\ by ${\Cal O}_{\Cal Z}$. In the case
when $Y$ is an isolated point this would mean that the length in
(3.9) is in fact the intersection number of ${\Cal Z}_1,\ldots,
{\Cal Z}_r$ at $y$. In particular one may ask, when does the
intersection number along $Y$ depend only on the scheme $Y$?
Related to this question is the problem of the singularities of the
schemes ${\Cal Z}(d,\omega)$: under which conditions are they
Cohen-Macaulay, or even locally complete intersections? In general
they are neither \cite{\krhb}.
\par\noindent

Our next task will be to investigate the structure of the
supersingular locus ${\Cal M}^{ss}\subset {\Cal
M}\times_{{\roman{Spec}}\ {\Bbb Z}_{(p)}} {\roman{Spec}}\  {\Bbb
F}_p$.

\vfill\eject
\subheading{\Sec4. Structure of the supersingular locus}

As mentioned in the introduction, the results of this section are a
presentation of results of Moret-Bailly \cite{\moretbailly} and Oort
\cite{\oort}.
A similar presentation was independently given by Kaiser \cite{\kaiser}.
\par
We put $\F=\bar\F_p$, and let $W=W(\F)$ be the ring of Witt vectors of $\F$
and $\K= W\tt_{\Z_p}\Q_p$ its quotient field. Also write $W[F,V]$ for the
Cartier ring of $\F$.

Throughout this section, we assume that $p\nmid
D(B)$, and we fix an isomorphism $\Cal O_C\tt_\Z\Z_p\simeq
M_4(\Z_p)$.

Suppose that $\xi=(A,\iota,\lambda,\bar\eta^p)\in \Cal
M^{\text{ss}}(\F)$, and let $A(p)$ be the p-divisible (formal)
group of $A$. The action of $\Cal O_C\tt_\Z\Z_p\simeq M_4(\Z_p)$ on
$A(p)$ then induces a decomposition $A(p)\simeq A_0(p)^4$, where
$A_0(p)$ is a p-divisible formal group of dimension $2$ and height
$4$. Let $L_0$ be the (contravariant) Dieudonn\'e module of
$A_0(p)$ and let $\L= L_0\tt_W\K$ be the associated isocrystal.
This does not depend on the choice of $\xi\in \Cal
M^{\text{ss}}(\F)$, up to isomorphism.
\hfill\break
More precisely,
we fix a base point $\xi_o = (A_o, \iota_o, \lambda_o,
\eta_o^p)\in {\Cal M}^{ss}({\Bbb F})$ and let ${\Cal
L}=L_0\otimes_W{\Cal K}$ be the isocrystal associated to it. The
isocrystal $\L$ has a polarization $<\ ,\ >$, is isoclinic with
slope $\frac12$, and has $\dim_{\K}\L=4$. Then $F$ is $\s$-linear,
$V=pF^{-1}$ is $\s^{-1}$-linear, and
$$<Fx,y>\ =\ <x,Vy>^\s.\tag4.1$$
If $\xi = (A,\iota,\lambda,\overline{\eta}^p)\in {\Cal
M}^{ss}({\Bbb F})$ is another point, then the choice of an isogeny
between $\xi$ and $\xi_o$ defines a $W$-lattice $L\subset {\Cal
L}$.

For a $W$-lattice $L\subset
\L$ of rank $4$, set
$$L^\perp=\{\ x\in \L\ \mid\  <x,L>\,\subset W\ \}.\tag4.2$$
\proclaim{Definition 4.1}  a) A $W$-lattice $L$ in $\L$ is
{\bf special} if and only if $L=c\cdot L^\perp$, for some $c\in \K^\times$.
Otherwise, $L$ is {\bf non-special}.
\hfill\break
b) A $W$-lattice $L$ in ${\Cal L}$ is {\bf admissible} if
$$L\supset FL\supset pL\ \ .$$
\endproclaim
We note that, if $L$ is an admissible lattice, then, since ${\Cal L}$ is
isoclinic of ${\roman{slope}}\ 1/2$, we have
$${\roman{dim}}_{\Bbb F}\  L/FL=2\ \ .$$
We define a set of lattices as follows:
$$X=\{ L\subset {\Cal L};\ L\ \text{admissible and special}\}\ \ .
\tag4.3$$
If $L\in X$ then $FL\in X$. This follows from $(FL)^\bot =V^{-1}\cdot
L^\bot$, cf.\ (4.1).
\par
The conditions in our moduli problem imply that the lattice
$L\subset {\Cal L}$ associated to $\xi\in {\Cal M}^{ss}({\Bbb F})$
and an isogeny between $\xi$ and $\xi_o$ actually lies in $X$. Note
that each admissible lattice is the Dieudonn\'e module of a
$p$-divisible formal group of dimension 2 and height 4 over ${\Bbb
F}$.
\par
Let $L\subset {\Cal L}$ be a lattice. Then $L$ is special or
non-special depending on whether the generalized index
$\left\lbrack L^\bot :L\right\rbrack\in {\Bbb Z}$ is divisible by 4
or not.  In any case this index is always even. If $L$ is special
and $L'\subset L$ is an inclusion of index 1, then $L'$ is
non-special and conversely.
\par
If $L$ is special, we can replace $L$ by $\a\cdot L$, for $\a\in
\K^\times$ to obtain a lattice with $L=L^\perp$ or $L=p L^\perp$. If $L$ is
non-special, we can scale to obtain a lattice with either
$$L\mathop{\subset}\limits_{\ne} L^\perp\mathop{\subset}\limits_{\ne} p^{-1}L,
\ \text{or}\ L^\bot\mathop{\subset}\limits_{\ne} L
\mathop{\subset}\limits_{\ne} p^{-1}L^\bot,\tag4.4$$
with all indices equal to 2. We will call lattices scaled in this
way {\bf standard}.

Recall from (3.8.vi) that
$\End^0(A_{\xi_o},\iota)^{\text{op}}=:C'\simeq M_2(B')$, where $B'$
is the definite quaternion algebra over $\Q$ with the same local
invariants as $B$ at all primes $\ell\ne p$. As before, let $V'=\{\
x\in M_2(B')\mid x'=x\text{ and } {\roman{tr}}(x)=0\
\}$. Let
$$G'=\{\ g\in C^{\prime,\times}\mid g V' g^{-1}=V'\ \text{ and } \ gg'=\nu(g)\
\}.\tag4.5$$
Note that the action of $G'(\Q_p)$ on $A_{\xi_o}(p)$ up to isogeny
passes to $\L$. In fact,
$$G'(\Q_p)\simeq \{\ g\in GL(\L)\mid \  <gx,gy>\,= \nu(g)<x,y>,\ \
Fg=gF \ \}.\tag4.6$$
Here $\nu(g)\in {\Cal K}^\times$.
\hfill\break
The action of $G'(\Q_p)$ preserves the set of lattices $X$. Fix an
isomorphism $B(\A_f^p)\simeq B'(\A_f^p)$ and, hence, an isomorphism
$G(\A_f^p)\simeq G'(\A_f^p)$. Then, the usual analysis identifies
$G'({\Bbb Q})$ with the group of self-isogenies of $\xi_o$ and
yields an isomorphism
$$\Cal M^{\text{ss}}(\F) \simeq G'(\Q)\back\bigg( X\times
G'(\A_f^p)/K^p\bigg).\tag4.7$$

We will now describe the lattices in $X$ in more detail.
%First, we would like to determine the orbits of $G'(\qp)$ in $X$.

\proclaim{Definition 4.2} For $L\in X$, let
$$a(L) = \dim_\F L/(FL+VL)\ \ .$$
\endproclaim

\par
Since $a(L)={\roman{dim}}_{\Bbb F} {\roman{Hom}}_{W\lbrack
F,V\rbrack}(L, {\Bbb F})$, we see that $a(L)$ is the $a$-number, \cite{\oort},
of
the $p$-divisible group $A_0(p)$ associated to $L$, i.e.
$$a(L)= {\roman{Hom}}_{\Bbb F}(\alpha_p, A_0(p))\ \ .$$
Since
$$\matrix {}&{}&L&{}&{}\\
\nass
{}&{}&\uparrow&{}&{}\\
\nass
{}&{}&FL+VL&{}&{}\\
{}&\nearrow&{}&\nwarrow&{}\\
FL&{}&{}&{}&VL\\
{}&\nwarrow&{}&\nearrow&{}\\
{}&{}&FL\cap VL&{}&{}\\
\nass
{}&{}&\uparrow&{}&{}\\
\nass
{}&{}&pL&{}&{}\endmatrix,\tag4.8$$
we have
$$a(L)=\cases 2 &\text{ if $FL=VL$, }\\
1&\text{ if $[L:FL+VL]=1$.}\endcases\tag4.9$$
Let
$$X_0=\{\ L\in X\mid a(L)=2\ \}.\tag4.10$$
Such lattices will be called {\bf superspecial}.

In addition to the superspecial lattices,
the following type of lattice will play a key role in the
description of the structure of $X$.
\proclaim{Definition 4.3} A lattice $\tl \subset {\Cal L}$ is {\bf
distinguished}
if $\tilde L$ is admissible and $F\tl=c\tl^\perp$ for some $c\in {\Cal
K}^\times$.
\endproclaim
We denote by $\tilde X$ the set of distinguished lattices. Obviously, if
$\tilde L\in\tilde X$ is distinguished, the index $\lbrack\tilde L^\bot
:\tilde L\rbrack$ is congruent to 2 mod 4. Thus $\tilde L$ is non-special.
Note that if $\tilde L$ is distinguished, then $F\tilde L=V\tilde L$.
Indeed, by (4.1) for any lattice $\tilde L$ we have $(F\tilde L)^\bot
=V^{-1}\tilde L^\bot$. Hence if $F\tilde L=c\cdot \tilde L^\bot$ we get
$$\tilde L=c (F\tilde L)^\bot =cV^{-1}\tilde L^\bot =c
V^{-1}c^{-1}F\tilde L= V^{-1}F\tilde L\ \ ,$$
i.e.\ $V\tilde L=F\tilde L$, as claimed. Similarly one sees that if $\tilde
L\in\tilde X$, then $F\tilde L\in\tilde X$.
\par
We note that if, in the identity defining a distinguished lattice
$\tilde L$, the order of $c$ is odd, then $\tilde L$ may be scaled
to be standard in the sense of the first alternative of (4.4)
above. If the order of $c$ is even, then $\tilde L$ can be scaled
to be standard in the sense of the second alternative of (4.4),
and hence,  $F\tilde L$ can be scaled to be standard in the sense of the
first alternative of (4.4).
\par
For any $\tl\in \tilde{X}$ and for any
$\F$-line $\ell\subset \tl/F\tilde L$, let $L =L(\ell)$
be the inverse image of $\ell$ in $\tl$. Thus
$$\matrix \tl&\supset&L&\supset& F\tilde L\\
\nass
\downarrow&{}&\downarrow&{}&\downarrow\\
\tl/F\tilde L&\supset&\ell&\supset&0\endmatrix\tag4.11
$$
\proclaim{Lemma 4.4} For $\ell\subset \tl/F\tilde L$,
 $L=L(\ell)\in X$.
\endproclaim

\demo{Proof} First, since $F\tl = V\tl$,
we have
$$FL \subset F\tl  \subset L,\tag4.12$$
and
$$FL\supset  FV\tl = p\tl \supset pL.\tag4.13$$
Hence $L$ is admissible.
\hfill\break
Next, we have
$$\tl\supset L\supset F\tilde L
,\tag4.14$$
where all inclusions have index $1$. But $\tilde L$ is non-special and hence
$L$ is special.
\qed\enddemo

The above proof in fact shows the following. Suppose that $\tilde
L\in\tilde X$ with $F\tilde L=p\cdot \tilde L^\bot$. Then
$L(\ell)^\bot
=pL(\ell)$. If $F\tilde L^\bot=p\tilde L$, then $L(\ell)^\bot =L(\ell)$.

Thus to any distinguished $\tilde L$ we have associated a projective line
${\Bbb P}(\tilde L/F\tilde L)$ and a family of admissible special lattices
parametrized by the ${\Bbb F}$-points of this projective line. These
projective lines have a natural ${\Bbb F}_{p^2}$-structure which we now
describe.
\par
For any $W$-lattice $L$ in $\L$, we have
$$FL=VL \iff F^2L=FVL=pL\iff p^{-1}F^2 L=L.\tag4.15$$
\proclaim{Lemma 4.5} Suppose that $p^{-1}F^2 L=L$, and let
$$L_0=\{\ x\in L \mid p^{-1}F^2x=x\ \}.$$
Then $L_0$ is a $\zps$-module and
$$L_0\tt_{\zps}W  \simeq L.\qquad\qed$$
\endproclaim

If $\tl\in \tilde{X}$ is distinguished, then $\tl$ is preserved by the
$\s^2$-linear endomorphism $p^{-1}F^2$, and we have
$\tl\simeq \tl_0\tt_{\zps}W$. Moreover, $F\tl$ is also preserved by
$p^{-1}F^2$, and $(F\tl)_0 = F(\tl_0)$. Thus, the two dimensional
$\F$-vector space $\tl/F\tl$
has a natural $\fps$-structure:
$$\tl/F\tl \simeq \tl_0/F\tl_0\tt_{\fps}\F.\tag4.16$$
We may then view any line $\ell$ as an element of $\Bbb P(\tl_0/F\tl_0)(\F)$.
We denote by ${\Bbb P}_{\tilde L}$ the projective line ${\Bbb P}(\tilde
L_0/F\tilde L_0)$ over ${\Bbb
F}_{p^2}$.
\proclaim{Lemma 4.6} Under the isomorphism
$$\tl/F\tl \simeq \tl_0/F\tl_0\tt_{\fps}\F,$$
the automorphism induced by $p^{-1}F^2$ on $\tl/F\tl$
coincides with $1\tt \s^2$ on \hfill\break $\tl_0/F\tl_0\tt_{\fps}\F$.
Hence,
$$p^{-1}F^2(L(\ell)) = L(\s^2(\ell)),$$
where $\ell$ is identified with a point in ${\Bbb P}_{\tilde L}({\Bbb F})$.
\qquad\qed\endproclaim

\proclaim{Corollary 4.7} A lattice $L(\ell)$ associated to a
distinguished $\tl$ is superspecial, i.e.,
has $a(L(\ell))=2$, if and only if $\ell\in
\Bbb P_{\tilde L}(\fps)$.
\endproclaim

\proclaim{Proposition 4.8}
Suppose that $L\in X$ with $a(L)=1$, and let
$$\tilde L=F^{-1}(FL+VL)\ \ .$$
Then $\tilde L$ is distinguished and $L=L(\ell)$ for a unique line
$\ell\in{\Bbb P}_{\tilde L}({\Bbb F})\setminus {\Bbb P}_{\tilde L}({\Bbb
F}_{p^2})$.
\endproclaim
\demo{Proof}
Let $L^\bot =c\cdot L$. Then
$$
(F\tilde L)^\bot
=(FL)^\bot\cap (VL)^\bot
=V^{-1}L^\bot\cap F^{-1}L^\bot
=p^{-1}c\cdot (FL\cap VL).
$$
On the other hand, $F^2\tilde L=F^2L+pL$. Let $S=L/pL$, and let $f$ and $v$
be the $\sigma$-linear resp.\ $\sigma^{-1}$-linear endomorphisms of $S$ induced
by $F$ and $V$.
Since $FV=VF=p$, we have $fv=v f=0$ and so ${\roman{ker}}(f)=
{\roman{im}}(v)$ and ${\roman{ker}}(v)= {\roman{im}}(f)$ are 2-dimensional
subspaces of $S$. However, for any $L\in X$ there is some $j\geq 2$ with
$F^jL\subset pL$ and hence $f$ is nilpotent. If $f^2=0$, then $F^2L=pL$
since both lattices have index 4 in $L$ and this would imply $a(L)=2$,
contrary to our assumption. Therefore, since ${\roman{im}}(f)$ is
2-dimensional we must have that ${\roman{im}}(f^2)$ is one-dimensional and
${\roman{im}}(f^2)= {\roman{im}}(f)\cap {\roman{im}}(v)$. Hence
$$F^2\tilde L=F^2L+pL= FL\cap VL\ \ .$$
It follows that $F(F\tilde L)= p\cdot c^{-1}(F\tilde L)^\bot$. On the other
hand, $F\tilde L$ is admissible, since
$$pF\tilde L= p(FL+VL)\subset pL\subset F^2\tilde L= FL\cap VL\subset
FL\subset F\tilde L=FL+VL$$
where all inclusions are of index 1. It follows that $F\tilde L\in\tilde X$
and hence also $\tilde L\in\tilde X$. Finally $L=L(\ell)$ for the line
$$\ell= L/F\tilde L\subset \tilde L/F\tilde L\ \ .\qquad\qed$$
\enddemo
\par
We summarize the above construction in the following theorem.
\proclaim{Theorem 4.9}
There is a natural $G'({\Bbb Q}_p)$-equivariant map
$$\coprod_{\tilde L\in\tilde X} {\Bbb P}_{\tilde L}({\Bbb F})\longrightarrow
X$$
which induces a bijection
$$\coprod_{\tilde L}({\Bbb P}_{\tilde L}({\Bbb F})\setminus {\Bbb P}_{\tilde
L}({\Bbb F}_{p^2}))\overset\sim\to\longrightarrow X\setminus X_0\ \ .$$
The map associates to $(\tilde L, \ell)$, where $\ell\subset\tilde L
/F\tilde L$ is a line, the element $L=L(\ell)\in X$.
\par
The action of $g\in G'({\Bbb Q}_p)$ on the index set of the left hand side
is lifted in the obvious way to the whole set appearing on the left hand side.
\endproclaim
\noindent
{\bf Remark 4.10.}
It can be shown that the map above is in fact a morphism, i.e.,\ is the map
on ${\Bbb F}$-points induced by a morphism of schemes over
${\roman{Spec}}\ {\Bbb F}_p$,
$$\coprod_{\tilde L\in\tilde X} {\Bbb P}_{\tilde L} \longrightarrow {\Cal
M}^{ss}\ \ .$$ This can be shown by the method of Oort,
\cite{\oort}, or using Cartier theory, as in Stamm, \cite{\stamm}.
Using either of these methods one can construct a morphism of
schemes over ${\roman{Spec}}\, {\Bbb F}_p$,
$$G'({\Bbb Q})\setminus \left\lbrack (\coprod\limits_{\tilde L\in \tilde X}
{\Bbb P}_{\tilde L})
\times G' ({\Bbb A}_f^p)/K^p\right\rbrack \longrightarrow {\Cal M}^{ss}$$
which turns out to be the normalization of the curve ${\Cal M}^{ss}$.
\par
The `distinguished curves' cross at the superspecial points. To describe this,
it will be useful to have a normal form for superspecial lattices.
\proclaim{Lemma 4.11} Fix
$\delta\in {\Bbb Z}^\times_{p^2}$ with $\delta^\s=-\delta$. Let
$L\in X_0$ be  superspecial and standard.
\hfill\break
(i) Suppose that $L=L^\perp$. Then there is a basis $e_1$, $e_2$,
$e_3$, $e_4$ for $L$ over $W$ such that $e_3=Fe_1$, $e_4=Fe_2$,
$Fe_3=pe_1$, $Fe_4=pe_2$ and such that the matrix for the
polarization is
$$(<e_i,e_j>)_{i,j} = \delta\cdot\pmatrix 0&1_2\\-1_2&0\endpmatrix,$$
(ii) If $L=pL^\perp$, then $L=FL'$ where $L'\in X_0$ with $L'=(L')^\perp$.
\endproclaim

\proclaim{Proposition 4.12} Suppose that $L\in X_0$ is superspecial and
standard. \hfill\break
(i) If $L^\perp = L$, consider lattices $\tl$ such that
$L\mathop{\supset}\limits_{\ne} \tl\mathop{\supset}\limits_{\ne} FL$
and such that $F\tilde L=p\tilde L^\bot$. Such $\tl$'s are
distinguished; there are $p+1$ of them,
and they can be described explicitly as follows.
Let $e_1,\dots,e_4$ be a standard basis as in Lemma~4.11. Then
the distinguished $\tl$'s have the form
$$\tl = W(e_1+\mu e_2) + FL$$
where $\mu\in {\Bbb Z}^\times_{p^2}$ such that $\mu\mu^\s\equiv
-1\mod p$.\hfill\break (ii) If $L=pL^\perp$, then the distinguished
$\tl$'s containing $L$ with index $1$ are those associated, as in
(i),  to $L'=F^{-1}L$.
\endproclaim
\demo{Proof of Lemma~4.11} Since $p^{-1}F^2$ is a $\s^2$-linear
automorphism of $L$, we can write $L=L_0\otimes_{{\Bbb Z}_{p^2}}W$ for the rank
$4$ lattice
$L_0$ of fixed points of $p^{-1}F^2$. Let $S_0=L_0/pL_0$, a $4$-dimensional
symplectic vector space over $\fps$, and note that $FL_0/pL_0$ is an
isotropic 2-plane in $S_0$, which is paired with the quotient $L_0/FL_0$.
We can then choose $e_1$ and $e_2\in L_0$ whose images form a basis for
$L_0/FL_0$ and such that $<e_1,e_2>=0$, after modification by elements of
$FL_0$,
if necesssary. The elements $e_1$, $e_2$, $e_3:=Fe_1$ and $e_4:=Fe_2$ then
give a $W$-basis for $L$, and $Fe_3=F^2e_1=pe_1$, and $Fe_4=F^2e_2=pe_2$, as
required, since $e_1$ and $e_2\in L_0$. The matrix for the polarization is then
$$\pmatrix 0&A\\-{}^tA&0\endpmatrix\qquad\text{ where } A= <\und{e},F\und{e}>,
\text{ with } \und{e} =\pmatrix e_1\\e_2\endpmatrix.\tag4.17$$
Note that $\det(A)\in \zpsx$, and that
$$-{}^tA^\s= <F\und{e},\und{e}>^\s = <\und{e},V\und{e}> =
<\und{e},F\und{e}>=A,\tag4.18$$ since $V=pF^{-1}$ and so, on $L_0$,
$V= F^2\cdot F^{-1} = F$. If we change the vector $\und{e}$ to
$a\cdot \und{e}$, for $a\in GL_2(\zps)$, then $A$ changes to $a A
{}^ta^\s$. Since $\det(A)\in\zpsx$ and since the norm map $N:\zpsx
\lra \Z_p^\times$ is surjective, it is easy to check that, for a
suitable choice of $a$ we can obtain $a A{}^ta^\s = \delta\cdot
1_2$.
\qed\enddemo
\demo{Proof of Proposition~4.12} Let us prove (i).
Using the standard basis of Lemma~4.11, we have $L=[e_1,e_2,e_3,e_4]$ (the
square
brackets indicate the $W$-span) and $FL=[pe_1,pe_2,e_3,e_4]$. Any lattice $\tl$
with
$L\supset\tl\supset FL$ and with $[L:\tl]=1$ has the form
$$\tl = W\cdot(ae_1+be_2) + FL,\tag4.19$$
where at least one of $a$ and $b\in W$ is a unit. If $a$ is a unit, we can
write
$\tl = [ e_1+\mu e_2, pe_2,e_3,e_4]$. Then
$$F\tl = [e_3+ \mu^\s e_4, pe_4,pe_1,pe_2]\quad\text{ and }\quad
p\tl^\perp = [ e_4-\mu e_3, pe_4, pe_1,pe_2].\tag4.20$$
Comparing, we see that $\mu$ must be a unit and that $\mu\mu^\s \equiv -1\mod
p$,
as claimed. It is easy to check that the case in which $a$ is not a unit
yields no solutions. The assertion (ii) is trivial.
\qed\enddemo

\proclaim{Corollary 4.13}
The map appearing in Theorem 4.9.\ is surjective. Any lattice in
$X_0$ has $p+1$ preimages which all lie on distinct lines. In fact,
the preimages of $L\in X_0$ correspond to the distinguished
lattices $F^{-1}\tilde L$ where $\tilde L$ ranges over the lattices
associated to $L$ in (i) of Proposition~4.12 ( resp. to
distinguished lattices $\tilde L$ associated to $L$ in (ii) of
Proposition~4.12). Finally, the images of two distinct lines ${\Bbb
P}_{\tilde L}$ and ${\Bbb P}_{\tilde L'}$ have at most one lattice
in common which then lies in $X_0$.
\endproclaim

\demo{Proof} The last assertion follows since, if $L,L'\in X$, $L\ne L'$, both
lie on
${\Bbb P}_{\tilde L}$, then $\tilde L=L+L'$.
\qed\enddemo
\par
The next result gives a standard basis for a distinguished lattice.
\proclaim{Lemma 4.14} Let $\tl$ be a distinguished lattice which is standard.
\hfill\break
(i) If $F\tilde L=p\tilde L^\bot$, then there
exists a $W$-basis $e_1,\dots,e_4$ of $\tl$ such that
$e_3=Fe_1$, $e_4=Fe_2$, $Fe_3=pe_1$, $Fe_4=pe_2$, and such that
the polarization has matrix
$$(<e_i,e_j>)_{i,j} = \delta\pmatrix {}&1&{}&{}\\-1&{}&{}&{}\\
{}&{}&{}&-p\\{}&{}&p&{}\endpmatrix.$$
(ii) If $F\tilde L^{\bot}= p\tilde L$, then $\tilde L=F\tilde L'$ where $\tilde
L'\in\tilde X$ with $F\tilde
L'=p\cdot \tilde L^{'\bot}$.
\endproclaim
\demo{Proof of (i)} Let $\tl_0$ be the fixed points of
$p^{-1}F^2$ on $\tl$. Since $F\tl = p\tl^\perp$, $<\ ,\ >$ induces
a nondegenerate symplectic form on the two dimensional
$\fps$-vector space $\tl_0/F\tl_0$. Choose $e_1$ and $e_2\in\tl_0$
whose images in $\tl_0/F\tl_0$ are a basis for this space and such
that $<e_1,e_2>=\delta$. Let $e_3=Fe_1$ and $e_4=Fe_2$, so that, as
in Lemma 4.11, $Fe_3=F^2e_1=pe_1$ and $Fe_4=F^2e_2=pe_2$. The
polarization then has matrix
$$\pmatrix \delta J&A\\-{}^tA&-p\delta J\endpmatrix\tag4.21$$
where $J=\pmatrix {}&1\\-1&{}\endpmatrix$ and $A=<\und{e},F\und{e}> =
-{}^tA^\s$,
as in the proof of Lemma~4.11. In the present case, however,
$A\equiv 0\mod p$. A Hensel's Lemma argument shows that we can
replace $\und{e}$ by $a\und{e}+bF\und{e}$ with $a\in GL_2(\zps)$
and $b\in M_2(\zps)$ to achieve $A=0$, while preserving the
condition $<\und{e},\und{e}>=\delta J$.
\qed\enddemo

Recall that $G'(\Q_p)$, given by (4.6) above, acts on the set of
admissible lattices. For any lattice $L$, $(gL)^\perp =
\nu(g)^{-1}g(L^\perp)$. If $L\in X$ is a special lattice, with
$L=c\cdot L^\perp$, then $gL = \nu(g)c\cdot (gL)^\perp$, so that
$gL$ is again special. Moreover, $a(gL)=a(L)$ so that the subset of
superspecial lattice is preserved. Also, if $\tl$ is distinguished,
and if $g\in G'(\Q_p)$, then $g\tl$ is again distinguished. Since
the valuation of $\nu(g)$ is an arbitrary integer,  for any
$L\in X$ (resp.\ $\tilde L\in\tilde X$) there is $g\in G'({\Bbb
Q}_p)$ such that $gL$ (resp.\ $g\tilde L$) is standard with
$(gL)^\bot=gL$ (resp.\ $F(g\tilde L)=p\cdot (g\tilde L)^\bot$). By
Lemmas~4.11 and 4.14, we have:
\proclaim{Corollary 4.15} $G'(\Q_p)$ acts transitively on the
set of superspecial lattices and on the set of distinguished lattices.
\endproclaim

We would finally like to compute the stablizers in $G'(\Q_p)$
of the superspecial and distinguished lattices.

Let $B'$ be as above, and, identifying $\qps$ with a subfield of
$B'_p$,  write $B'_p=\qps + \Pi \qps$ for an element
$\Pi\in B_p^{\prime,\times}$ with $\Pi^2=p$ and such that $\Pi a = a^\s \Pi$,
for $a\in \qps$.
Let $\L_0$ be the fixed set for the automorphism $p^{-1}F^2$ of $\L$, and
let $\Pi$ operate on $\L_0$ by $F$. By construction, $\Pi^2=p$, and
so $\L_0$ is naturally a left vector space over $B'_p$ of dimension $2$.
\proclaim{Lemma 4.16} Let
$$\End_{\Cal K}(\L,F):=\{\a\in \End_{\Cal K}(\L) \mid F\a=\a F\ \}.$$
Then,
$$\End_{\Cal K}(\L,F) =\End_{\qps}(\L_0,F) = \End_{B'_p}(\L_0).$$
\endproclaim

The polarization on $\L$ induces a $\qps$-bilinear symplectic form on
$\L_0$, which still satisfies $<Fx,y>=<x,Vy>^\s$.
\proclaim{Lemma 4.17} Let $U$ be a left $B'_p$-vector space with a
$B'_p$-Hermitian form $(\ ,\ ):U\times U\rightarrow B'_p$. Thus
$(bx,cy)=b(x,y) c^\iota$ and $(y,x)=(x,y)^\iota$, where $b\mapsto b^\iota$ is
the
main involution on $B'_p$. Write
$$(x,y) = (x,y)_0\delta + (x,y)_1\delta \Pi,$$
where $(x,y)_0$ and $(x,y)_1\in \qps$. Then,
$$(\ ,\ )_1:U\times U\lra \qps$$
is a symplectic $\qps$-bilinear form on the $\qps$-vector space $U$ such that
$$(\Pi x,y)_1 = (x,\Pi y)_1^\s,\tag*$$
and
$$(x,y)_0 = -(x,\Pi y)_1.$$
The map $(\ ,\ )\mapsto (\ ,\ )_1$ yields a bijection
between the space of $B'_p$-Hermitian forms on $U$ and the space
of symplectic forms satisfying (*). Moreover,
$$\align
G'(\Q_p) &= \{\ g\in GL(\L)\mid <gx,gy>=\nu(g)<x,y>\text{ and } Fg=gF\}\\
{}&\simeq\{\ g\in GL_{B'}(\L_0)\mid (gx,gy) =\nu(g) (x,y)\ \}.
\endalign$$
\endproclaim
\demo{Proof} We just check
the behavior of $\Pi$. We have
$$\align
(\Pi x, y) &= \Pi(x,y)_0\delta +\Pi(x,y)_1\delta\Pi \tag4.22\\
{}&= -p(x,y)_1^\s\delta - (x,y)_0^\s\delta\Pi\\
\noalign{and}
(x,\Pi y)&= -(x,y)_0\delta\Pi -p(x,y)_1\delta\ \ . \tag4.23\endalign$$
Thus
$$(\Pi x,y)_1 = - (x,y)_0^\s = (x,\Pi y)_1^\s, \tag4.24$$
as required. It is at this point that the factor of $\delta$ is
required in the formulas.
\qed\enddemo

Let $\Cal O'=\Cal O_{B'_p} =\zps + \Pi \zps$ be the maximal order in
$B'_p$. If $L$ is an admissible lattice such that $p^{-1}F^2L=L$,
then the fixed point set $L_0$ of $p^{-1}F^2$ is naturally an $\Cal O'$-lattice
in the $B'_p$-vector space $\L_0$, and $\dim_\fps L_0/\Pi L_0 =2$.
Conversely, given any $\Cal O'$-lattice $\Lambda$ with this last property, we
set $F=\s\tt\Pi$ on $L=L(\Lambda):=W\tt_\zps \Lambda$. Then since $\Pi^2=p$ on
$\Lambda$, we have $L\supset FL\supset pL$ and $\dim_\F L/FL=2$, i.e.,
$L$ is admissible, and $p^{-1}F^2L=L$.
The following is easily checked, using the formulas of Lemma~4.17.
\proclaim{Lemma 4.18} (i) Suppose that $L\in X_0$ is superspecial with
$L=L^\perp$, and let $e_1,\dots,e_4$ be a standard basis as in
Lemma~4.11. Then $e_1'=\delta^{-1}e_1$ and $e_2'=\delta^{-1}e_2$ is
an $\Cal O'$-basis for $L_0$, and the matrix for the
$B'_p$-Hermitian form on $\L_0$ is $((e'_i,e'_j))_{i,j} =
1_2$.\hfill\break (ii) Suppose that $\tl\in \tilde{X}$ is
distinguished and that $F\tilde L=p\tilde L^\bot$, and let
$e_1,\dots,e_4$ be a standard basis as in Lemma~4.14. Then
$e_1'=\delta^{-1}e_1$ and $e_2'=-\delta^{-1}e_2$ form an $\Cal
O'$-basis for $\tl_0$ and
$$((e'_i,e'_j))_{i,j} =  \pmatrix {}&\Pi\\-\Pi&{}\endpmatrix.$$
\endproclaim
Thus, in classical language, cf.\ \cite{\shimura},
\cite{\hashimoto}, the superspecial lattices correspond to local
components of the principal genus of quaternion Hermitian lattices,
while the distinguished lattices correspond to local components of
a non-principal genus of such lattices.
%\cite{\shimura},
%\cite{\hashimoto}, etc..
\par
In less classical language we may describe our results in terms of the
Bruhat-Tits building of the adjoint group $G'_{\roman{ad}}$ over ${\Bbb
Q}_p$, comp.\ also \cite{\kaiser}. The building ${\Cal
B}(G'_{\roman{ad}}, {\Bbb Q}_p)$ is a tree and may be identified with the
fixed points
$${\Cal B}(G'_{\roman{ad}}, {\Bbb Q}_p)= {\Cal B}(G'_{\roman{ad}}, {\Cal
K})^F\ \ .$$ The special vertices in ${\Cal B}(G'_{\roman{ad}},
{\Bbb Q}_p)$ correspond to the equivalence classes of lattices
$L\subset {\Cal L}$ which are $F$-invariant. Here two lattices
$L_1$ and $L_2$ are equivalent if $L_1$ is homothetic to $L_2$ or
to $L_2^\bot$. Hence the special vertices are in one-to-one
correspondence with the distinguished lattices $\tilde L$ which are
standard and with $F\tilde L=p\tilde L^\bot$. The non-special
vertices in ${\Cal B}(G'_{ad}, {\Bbb Q}_p)$ correspond to the edges
in ${\Cal B}(G'_{\roman{ad}}, {\Cal K})$ whose vertices are
interchanged by $F$. Equivalently, they correspond to pairs $\{
L,FL\}$ of lattices in $X_0$ which are standard. We thus obtain
bijections
$$
\align
\tilde X
&
\leftrightarrow
{\Bbb Z}\times \{\text{special vertices in}\ {\Cal B}(G'_{\roman{ad}}, {\Bbb
Q}_p)\}\\
\intertext{and}
X_0
&
\leftrightarrow {\Bbb Z}\times \{\text{non-special vertices in}\ {\Cal
B}(G'_{\roman{ad}}, {\Bbb Q}_p)\}\ .
\endalign
$$
These bijections are $G'({\Bbb Q}_p)$-equivariant, where $g\in
G'({\Bbb Q}_p)$ acts on the ${\Bbb Z}$-component on the right via
$n\mapsto n+{\roman{ord}}(\nu(g))$. The action of $F$ on the left
corresponds to the translation $n\mapsto n+1$ on the first factor
and the trivial action on the second factor on the right.
Furthermore, a lattice $L\in X_0$ and $\tilde L\in \tilde X$ are
incident (i.e.\ $L\in {\Bbb P}_{\tilde L})$ if and only if the
corresponding vertices of ${\Cal B}(G'_{\roman{ad}}, {\Bbb Q}_p)$
lie on one and the same edge.

In these terms the stabilizer $K^{\text{d}}$ of a
distinguished lattice $\tilde L\in\tilde X$ is a special maximal
compact subgroup of $G'({\Bbb Q}_p)$, and the stabilizer $K^{\text{ss}}$ of
a superspecial lattice $L\in X_0$ is a
non-special maximal compact subgroup of $G'({\Bbb Q}_p)$.
\par
{\bf Remark 4.19} We return, for a moment, to the global situation,
and recall that $\tilde{X}$ is the set of distinguished lattices in $\L$. As
observed in
Remark~4.10, our
calculations `show' that the supersingular locus $\Cal M^{\text{ss}}$
is a union of rational curves and that the irreducible components
are in bijection with the set
$$G'(\Q)\back \bigg(\tilde{X} \times G(\A_f^p)/K^p\bigg)
\simeq G'(\Q)\back \bigg( G'(\Q_p)/K^{\text{d}}_p\times G(\A_f^p)/K^p\bigg),
\tag4.25$$
where $K^{\text{d}}_p$ is the stabilizer in $G'(\Q_p)$ of a fixed
distinguished lattice $\tl\in \tilde{X}$.
These curves cross, $p+1$ at a time, at the superspecial points,
and there are $p^2+1$ such crossing points on each component.
The set of all crossing points is in bijection with the
set
$$G'(\Q)\back\bigg(X_0\times G(K^p_f)/K^p\bigg)
\simeq G'(\Q)\back \bigg( G'(\Q_p)/K^{\text{ss}}_p\times G(\A_f^p)/K^p\bigg),
\tag4.26$$
where $K^{\text{ss}}_p$
is the stabilizer in $G'(\Q_p)$ of a fixed superspecial lattice $L\in X_0$.

We finally observe two consequences of our description of $\Cal M^{\text{ss}}$.

Fix a factorization $D(B)=D_1 D_2$, and let $K=\prod_\ell K_\ell$ be the
compact open
subgroup of $G(\A_f)$ with local factors
$$K_\ell = \cases K_\ell^{\text{ss}} &\text{ if $\ell\mid D_1$,}\\
K_\ell^{\text{d}} &\text{ if $\ell\mid D_2$,}\\
K_\ell^0 &\text{ if $\ell\nmid D(B)$.}
\endcases$$
Here, we have fixed a maximal order $R$ in $B$, and for $\ell\nmid
D(B)$, we fix an isomorphism $M_2(B_\ell)\simeq M_4(\Q_\ell)$ such
that $M_2(R_\ell)\simeq M_4(\Z_\ell)$.  Then let
$K_\ell^0=G(\Q_\ell)\cap M_4(\Z_\ell)$. Thus, for $\ell\mid D_1$
(resp. $\ell\mid D_2$), $K_\ell$ is the stabilizer of a Hermitian
$\Cal O_{B_\ell}$-lattice of principal (resp. non-principal) type,
and, for $\ell\nmid D(B)$, $K_\ell$ is a hyperspecial maximal
compact subgroup of $G(\Q_\ell)$. Note that, in contrast to the
general assumptions above, $K$ is not neat. Still, for a fixed
prime $p\nmid D(B)$, one can consider the coarse moduli space $\Cal
M_K$ (the quotient by a finite group of one of the schemes
considered above) and its points over $\F$.  Let $B^{(p)}$ denote
the definite quaternion algebra with $D(B^{(p)}) = D(B)p$. Then, by
(4.25), the components of the supersingular locus in the fiber of
$\Cal M_K$ at $p$ correspond to the classes of maximal Hermitian
lattices in the genus of type $(D_1,pD_2)$ for $B^{(p)}$. An
explicit formula for this number $H(D_1,pD_2)$ was found by
Hashimoto and Ibukiyama \cite{\hashibuki}. In the case $D(B)=1$, so
that $B=M_2(\Q)$, the abelian varieties parameterized by $\Cal
M_K(\F)$ have the form $A\simeq A_0^4$, where $A_0$ is a
principally polarized abelian surface. Thus, in this case, $\Cal
M_K\simeq A_{2,1}$, and the description of the supersingular locus
reduces to some of the information given by Katsura and Oort
\cite{\katoort}, Theorem~5.7, and Ibukiyama, Katsura and Oort
\cite{\ibukatoort}. In particular, the number of irreducible
components of the supersingular locus is $H(1,p)$.

As another example, fix a square free positive integer $D$ and
distinct primes $p_1$ and $p_2$ relatively prime to $D$. Consider
indefinite quaternion algebras $B_1$ and $B_2$ over $\Q$ with
discriminants $D(B_1) = D p_1$ and $D(B_2)=D p_2$. Let $G_1$ and
$G_2$ be the associated groups, via (1.3). As in (4.5), let $G_1'$
be the twist of $G_1$ at $p_2$ and let $G_2'$ be the twist of $G_2$
at $p_1$. These groups are both associated to the definite
quaternion algebra $B_1^{(p_2)}\simeq B_2^{(p_1)}$, and are
isomorphic. Fix an isomorphism $G_1'\simeq G_2'$ and compatible
isomorphisms
$$G_1(\A_f^{p_1p_2})\simeq G_1'(\A_f^{p_1p_2}) \simeq G_2'(\A_f^{p_1p_2})
\simeq G_2(\A_f^{p_1p_2}),$$
and let $K^{p_1p_2}=K_1^{p_1p_2}=K_2^{p_1p_2}$ be a
sufficiently small compact open subgroup.
Also let
$$\align
K_{1,p_1} &= K_{p_1}^{*_1},\qquad\text{for $*_1 = \text{\rm d}$ or $\text{\rm
ss}$,}\\
K_{1,p_2} &= K_{p_2}^0,\\
K_{2,p_1} &= K_{p_1}^0,\\
K_{2,p_2} &= K_{p_2}^{*_2},\qquad\text{for $*_1 = \text{\rm d}$ or $\text{\rm
ss}$,}
\endalign$$
where the notation is as above.
Let
$$\align
K_1^{*_1} &= K^{p_1p_2}K_{1,p_1}K_{1,p_2},\\
K_2^{*_2} &= K^{p_1p_2}K_{2,p_1}K_{2,p_2}.\endalign$$
Let $\Cal M_1^{*_1}=\Cal M_{K_1^{*_1}}$ and $\Cal M_2^{*_2}=\Cal M_{K_2^{*_2}}$
be the corresponding
moduli schemes, defined over
$\Z_{(p_2)}$ and $\Z_{(p_1)}$ respectively.

Then, using (4.25) and (4.26), there are (non-canonical but
equivariant) bijections between various sets of irreducible
components or crossing points as follows:
$$\align
\text{Components}\bigg((\Cal M_2^{\text{\rm d}}\times
\F_{p_1})^{\text{s.s.}}\bigg)
&\simeq
\text{Components}\bigg( (\Cal M_1^{\text{\rm d}}\times
\F_{p_2})^{\text{s.s.}}\bigg),\\
\nass\nass
\text{Components}\bigg((\Cal M_2^{\text{\rm ss}}\times
\F_{p_1})^{\text{s.s.}}\bigg)
&\simeq
\text{Crossing points}\bigg( (\Cal M_1^{\text{\rm d}}\times
\F_{p_2})^{\text{s.s.}}\bigg),\\
\nass\nass
\text{Crossing points}\bigg((\Cal M_2^{\text{\rm d}}\times
\F_{p_1})^{\text{s.s.}}\bigg)
&\simeq
\text{Components}\bigg( (\Cal M_1^{\text{\rm ss}}\times
\F_{p_2})^{\text{s.s.}}\bigg),\\
\nass
%\noalign{\noindent and}
\nass
\text{Crossing points}\bigg((\Cal M_2^{\text{\rm ss}}\times
\F_{p_1})^{\text{s.s.}}\bigg)
&\simeq
\text{Crossing points}\bigg( (\Cal M_1^{\text{\rm ss}}\times
\F_{p_2})^{\text{s.s.}}\bigg).
\endalign$$
Here we have written $(\Cal M_1^{\text{\rm ss}}\times \F_{p_2})^{\text{s.s.}}$
for the supersingular locus of the fiber over $p_2$ of $\Cal M_1^{*_1}$, where
$*_1=
\text{\rm ss}$, for example.
These results are in the spirit of
those of Ribet \cite{\ribetone}, \cite{\ribettwo}, who considers components and
their
crossing points for the fibers of Shimura curves and modular curves at
primes of bad reduction.

\vfill\eject
\subheading{\Sec5. Endomorphism algebras and points of proper intersection}

In this section, we consider the points of intersection of the special cycles
in the supersingular locus, using the information obtained in section 4 about
the
structure of this locus. In particular, in the decomposition
$$
{\Cal Z}(d_1, \omega_1)\times_{\Cal M}\ldots\times_{\Cal M} {\Cal Z}
(d_r, \omega_r)=\coprod\limits_{\matrix
T\in {\roman{Sym}}_4({\Bbb Z}_{(p)})_{\ge 0}\\
{\roman{diag}}(T)=(d_1,\ldots, d_r)
\endmatrix} {\Cal Z}(T,\omega)$$
of (3.6), we fix a matrix $T$ and we obtain a criterion, in terms of
$T$, for ${\Cal Z}(T,\omega)$ to consist of isolated points.
We also show that, even when $\det(T)\ne0$, there can be components
of the supersingular locus in the image of ${\Cal Z}(T,\omega)$ in $\Cal
M^{\text{ss}}$.

We retain the notation of sections 2--4, and we begin by obtaining information
about the endomorphism rings of various types of admissible lattices.

For an admissible lattice $L$, let ${\Cal O}_L=\End_W(L,F)$ be
the $\Z_p$-algebra of $W$-linear endomorphisms of $L$ which
commute with $F$. Note that $\End_W(L,F)$ is an order in the
$\Q_p$-algebra $\End_{\Cal K}(\L,F)=C'_p= C'\tt_\Q\Q_p\simeq
M_2(B'_p)$. Also, observe that $\End_W(L,F)=\End_W(F^jL,F)$ for
any $j\in \Z$. If $L=c\cdot L^\perp$ is special, we have
$$(F^jL)= p^jc\cdot (F^jL)^\perp.\tag5.1$$
Thus, to determine ${\Cal O}_L$ for $L\in X_0$ we may assume
$L=L^\perp$.

By Lemma~4.18, we immediately have the following.
\proclaim{Lemma 5.1} For any superspecial lattice $L\in X_0$,
(resp. any distinguished lattice $\tl\in \tilde{X}$)
$\End_W(L,F)$ (resp. $\End_W(\tl,F)$) is a maximal order in $C'_p$.
\endproclaim
In either case, this order is isomorphic to $M_2(\Cal O')$, where
$\Cal O' = \zps+\Pi\zps$, as in section 4. The map $M_2(\Cal
O')\lra M_2(\fps)$ given by reduction modulo $\Pi$ can be described
as follows. Consider the case of $L\in X_0$. As in section 4, let
$L_0$ be the fixed points of $p^{-1}F^2$ on $L$. Then define
$${\roman{red}}_L:\End_W(L,F) \lra \End_\fps(L_0/FL_0)\simeq M_2(\fps)\tag5.2$$
as the composition
$$\End_W(L,F) \isoarrow \End_\zps(L_0,F) \lra \End_\fps(L_0/FL_0).\tag5.3$$
This map is surjective. The surjective reduction map for $\tilde
L\in\tilde X$
$${\roman{red}}_{\tl}:\End_W(\tl,F)\lra \End_\fps(\tl_0/F\tl_0)\simeq
M_2(\fps)\tag5.4$$
is defined analogously. Note that $\tl/F\tl \simeq \tl_0/F\tl_0\tt_\fps \F$,
and that the endomorphism $\bar{\a}$ induced on $\tl/F\tl$ by $\a\in
\End_W(\tl,F)$
is ${\roman{red}}_\tl(\a)\tt 1$.

Next, suppose that $L\in X\setminus X_0$,
and let $\tl$ be the unique distinguished lattice associated to $L$ by
Proposition~4.8. Recall that $F\tl =  FL+VL$. In
particular, for every element $\a\in \End_W(L,F)$, $\a F\tl \subset F\tl$,
so that $\a\tl\subset\tl$, and
there is a natural homomorphism which is injective,
$$\End_W(L,F)\hookrightarrow \End_W(\tl,F).\tag5.5$$
On the other hand, there is a unique line $\ell\subset \tl/F\tl$
such that $L =L(\ell)$ is the inverse image of $\ell$ in $\tl$.
\proclaim{Lemma 5.2} Let $L\in X\setminus X_0$. With the notations introduced
above,
$$\End_W(L,F) = \{\ \a\in \End_W(\tl,F)\mid \bar{\a}(\ell)\subset \ell\ \}.$$
Here $\bar{\a}$ is the endomorphism of $\tl/F\tl$ induced by $\a$.
In fact, there are two possibilities.
\hfill\break
(i) If $\ell\in \Bbb P_{\tl}(\F) - \Bbb P_{\tl}(\Bbb F_{p^4})$, then
$$\End_W(L,F) = ({\roman{red}}_\tl)^{-1}(\fps\cdot 1_2).$$
(ii) If $\ell\in \Bbb P_{\tl}(\Bbb F_{p^4}) - \Bbb P_{\tl}(\fps)$, then
$$\End_W(L,F) = ({\roman{red}}_\tl)^{-1}(\Bbb F_{p^4}),$$
for some embedding $\Bbb F_{p^4}\hookrightarrow M_2(\fps)$.
\endproclaim
\demo{Proof} As remarked above, the automorphism of
$\tl/F\tl =\tilde L_0/F\tilde L_0\otimes_{{\Bbb F}_{p^2}}{\Bbb
F}$ induced by $p^{-1}F^2$ is just $1\tt\s^2$. Since, for any
$\a\in
\End_W(\tl,F)$, $\bar{\a}$ commutes with this automorphism,
$\bar{\a}(\ell)\subset \ell$ implies that
$\bar{\a}(\s^2(\ell))\subset\sigma^2(\ell)$. Since a non-scalar
endomorphism can have at most two eigenlines,
$\bar{\a}(\ell)\subset \ell$ and $\s^4(\ell)\ne\ell$ implies that
$\bar{\a} = a\cdot 1_2$, for $a\in \fps$. If $\s^4(\ell)=\ell$ but
$\s^2(\ell)\ne \ell$, and if $\bar{\a}$ is not a scalar
endomorphism, then $\ell$ and $\s^2(\ell)$ are the distinct
eigenlines of $\bar{\a}$. Then $\fps[\bar{\a}] \simeq \Bbb
F_{p^4}$, and any endomorphism $\bar{\beta}$, with $\beta\in
\End_W(L,F)$ must lie in this subfield of $M_2(\fps)$.
\qed\enddemo

Note that the lattices in (ii) of lemma 5.2 are characterized
intrinsically by the condition that $F^4L=p^2L$ but $F^2L\ne pL$.
We let $X_{(ii)}$ be the set of lattices appearing in (ii) and
$X_{(i)}=X\setminus X_{(ii)}\setminus X_0$ the set appearing in
(i). Recall that $\ol=\End_W(L,F)$ and $\otl=\End_W(\tl,F)$. Then
$$\align
{\roman{red}}_L(\ol)={\roman{red}}_L(\End_W(L,F)) &\simeq M_2(\fps)
\qquad\, \text{if}\  L\in X_0\\
{\roman{red}}_\tl(\ol)={\roman{red}}_\tl(\End_W(L,F)) &\simeq
\Bbb F_{p^4}  \qquad\qquad\ \text{if}\  L\in X_{(ii)}
\tag5.6\\
{\roman{red}}_\tl(\ol)={\roman{red}}_\tl(\End_W(L,F)) &\simeq \fps
\qquad\qquad\ \text{if}\  L\in X_{(i)}.
\endalign$$
In particular, the endomorphism algebras of all the $L$'s with $L\in X_{(i)}$
and with a given associated $\tl$ coincide. Any
endomorphism of one such $L$ preserves
all lattices $L'\in X$ in the image of ${\Bbb P}_{\tilde L}$.

Recall that $C'_p= \End_{\Cal K}(\L,F)$, and let
$$V'_p  = \{ \ x\in C'_p\mid x^\ast=x, \text{ and } {\roman{tr}}^0(x)=0\
\}.\tag5.7$$
Note that ${\Cal O}_L$ and ${\Cal O}_{\tilde L}$ are invariant
under the involution $\ast$ on $C'_p$. Indeed, let $L^\bot =cL$
and $x\in {\Cal O}_L$. Then
$$x^\ast (L^\bot)\subset L^\bot\ \ ,\ \ \text{i.e.}\ x^\ast(L)\subset L\ \
.$$
Similarly, if $F\tilde L=c\tilde L^\bot$ and $x\in {\Cal
O}_{\tilde L}$, then $x^\ast(\tilde L^\bot)\subset \tilde
L^\bot$, i.e.\ $x^\ast(F\tilde L)\subset F\tilde L$, i.e.,
$x^\ast(\tilde L)\subset \tilde L$, since $x^\ast$ commutes with
$F$.
\par
For $L\in X$  and for $\tl\in \tilde{X}$, let
$$N_L= \End_W(L,F)\cap V'_p \qquad\text{ and }\qquad N_\tl = \End_W(\tl,F)\cap
V'_p.\tag5.8$$
These are ${\Bbb Z}_p$-lattices in $V'_p$ on which the quadratic
form given by squaring, $x^2 = q(x)\cdot id$ is valued in $\Z_p$.
\par
We now describe the reduction maps for distinguished and for superspecial
lattices. We start with the case of distinguished lattices.

\proclaim{Lemma 5.3}
Let $\tilde L\in \tilde X$, and put ${\frak n}_{\tilde L}=
{\roman {red}}_{\tilde L}(N_{\tilde L})$. Then ${\frak n}_{\tilde L}$ is equal
to
$$\{ x=a\cdot 1_2;\ a\in{\Bbb F}_{p^2},\ a^{\sigma}=-a\}$$
and the ${\Bbb F}_p$-valued quadratic form $q$ on ${\frak
n}_{\tilde L}$ is given by $x^2=q(x)\cdot 1_2$, i.e.\
$q(x)=-a\cdot a^{\sigma}$. In particular, $q$ does not represent
1 and hence the Clifford algebra $C({\frak n}_{\tilde L})$ is
isomorphic to ${\Bbb F}_{p^2}$. The following diagram is
commutative
$$\matrix
N_{\tilde L}
&
\overset q\to\longrightarrow
&
{\Bbb Z}_p\\
\llap{$\scriptstyle {\roman{red}}_{\tilde L}$}\big\downarrow
&&
\big\downarrow\\
{\frak n}_{\tilde L}
&
\overset q\to\longrightarrow
&
{\Bbb F}_p
&
.
\endmatrix$$
\endproclaim

\demo{Proof}
Replacing $\tilde L$ by $F^j\tilde L$ we may assume that $\tilde L$
is standard with $F\tilde L=p\cdot \tilde L^{\bot}$. The symplectic
form $<\ ,\ >$ on ${\Cal L}$ induces a nondegenerate
alternating pairing
$$<\ ,\ > : \tilde L/F\tilde L\times \tilde L/F\tilde
L\longrightarrow {\Bbb F}\ \ .\tag5.9$$
This pairing descends to a non-degenerate alternating ${\Bbb
F}_{p^2}$-bilinear pairing on $\tilde L_0/F\tilde L_0$ with values in ${\Bbb
F}_{p^2}$. The induced involution on ${\roman{End}}_{{\Bbb F}_{p^2}}(\tilde
L_0 /F\tilde L_0)$ is compatible with the reduction map,
$${\roman{red}}_{\tilde L}(x^\ast)= {\roman{red}}_{\tilde L}(x)^\ast\ \ ,\ \
x\in {\Cal O}_L\ \ .$$ Now any endomorphism $\overline x$ of the
2-dimensional symplectic vector space $\tilde L_0/F\tilde L_0$
over ${\Bbb F}_{p^2}$ with $\overline x^{\ast}=\overline x$ is a
scalar. Hence for $x\in N_{\tilde L}$ we get
$${\roman{red}}(x)= a\cdot 1_2\ \ ,\ \ a\in {\Bbb F}_{p^2}\ \ .$$
But $x\in N_{\tilde L}$ acts on $\tilde L_0/p\tilde L_0$ preserving the
subspace $F\tilde L_0/p\tilde L_0$. Since $x$ commutes with $F$, it acts on
the subspace as $a^{\sigma}\cdot 1_2$. The condition ${\roman{tr}}^0(x)=0$
implies therefore that $a=-a^{\sigma}$. Therefore we have proved that
${\frak n}_{\tilde L}$ is contained in the subspace above. It is easy to see
that we have in fact an equality. The remaining assertions are obvious.
\qed\enddemo

Next we consider the case of superspecial lattices.
\proclaim{Lemma 5.4} Let $L\in X_0$ and put ${\frak n}_L={\roman{red}}_L(N_L)$.
\hfill\break
(i) $\fn_L$ is isomorphic to
$$\align
&\{\ x\in M_2(\fps)\mid {}^tx^\s =x\text{ and } tr(x)=0\ \}\\
{}=&\{\ x=\pmatrix a&b\\b^\s&-a\endpmatrix\mid a\in \fp,\ b\in \fps\ \}.
\endalign$$
The $\fp$-valued quadratic form $q$ on $\fn_L$ is given by
$x^2=q(x)\cdot 1_2$, i.e., $q(x) = -(a^2+bb^\s)$. \hfill\break
(ii) Let $C(\fn_L)$ be the Clifford algebra of
the three dimensional quadratic space $\fn_L$. Then
the natural map $C(\fn_L)\isoarrow M_2(\fps)$
is an isomorphism.\hfill\break
(iii) The following diagram is commutative.
$$\matrix
N_{L} &
\overset q\to\longrightarrow
&
{\Bbb Z}_p\\
\llap{$\scriptstyle {\roman{red}}_{L}$}\big\downarrow
&&
\big\downarrow\\
{\frak n}_{L} &
\overset q\to\longrightarrow
&
{\Bbb F}_p
&
.
\endmatrix$$
\endproclaim
\demo{Proof}
Replacing $L$ by $F^jL$ we may assume $L^\bot =L$. On $L_0/FL_0$ we have
the non-degenerate anti-hermitian form
$$(\ ,\ ): L_0/FL_0\times L_0/FL_0\longrightarrow {\Bbb F}_{p^2}\tag5.10$$
induced by the formula
$$(v,w)=\langle \tilde v, F\tilde w\rangle\mod p\ \ ,\tag5.11$$
where $\tilde v$ and $\tilde w$ are representatives of $v$ and $w$ in $L_0$.
We may find a basis of $L_0/FL_0$ such that the induced involution on
$M_2({\Bbb F}_{p^2})$ is given by $x\mapsto {}^tx^{\sigma}$. Now the lemma
is proved in a way similar to Lemma 5.3.\ above.
\qed\enddemo
\par
We now return to the points of intersection of the special cycles
in the supersingular locus. Let $T\in{\roman{Sym}}_4({\Bbb
Z}_{(p)})$ with ${\roman{det}}\,T\ne 0$ and $\omega\subset V({\Bbb
A}_f^p)^4$ with corresponding special cycle ${\Cal Z}(T,\omega)$.
Let $\xi\in {\Cal Z}(T,\omega)$
correspond to the collection $(A_\xi,\iota,\lambda,\bar\eta^p;\j)$.
By Corollary~4.3, the point corresponding to the collection
$(A_\xi,\iota,\lambda,\bar\eta^p)$ lies in $\Cal M^{\text{ss}}(\F)$.
Thus, $\End^0(A_\xi,\iota)^{\text{op}} =C^0_\xi = C'\simeq
M_2(B')$, where $B'$ is the definite quaternion algebra over $\Q$
with discriminant $D(B)p$. The last isomorphism here can be chosen
so that the Rosati involution corresponds to the involution
$x\mapsto x'={}^tx^\iota$ of $M_2(B')$. Then, as in (3.8.vi),
$$V^0_\xi=V' \simeq \{x\in M_2(B')\mid x'=x\text{ and } {\roman{tr}}(x)=0\
\}.\tag5.12$$
The components $j_1,\ldots, j_4$ of $\bold j$ lie in $V'$, and therefore, in
particular,
we must have $T>0$ if ${\Cal Z}(T,\omega)$ is to be non-empty.
\par
Let $L$ be the contravariant Dieudonn\'e module of the formal group
$A_0(p)$, where we write $A_\xi(p)\simeq A_0(p)^4$, as in section
4. By choosing an isogeny of $\xi$ with the chosen base point
$\xi_o$ we obtain, as in section 4, an identification $\L=
L\tt_W\Cal K$ of its isocrystal with that of the base point. Then
$L\in X$, and there is a natural algebra homomorphism
$$C_\xi\tt_{\Z}\Z_p=\End(A_\xi,\iota)^{\text{op}}\tt_{\Z}\Z_p
\hookrightarrow \End_W(L,F)={\Cal O}_L\subset C'_p.\tag5.13$$

Let $N_L = \End_W(L,F)\cap V'_p$, as in (5.8) above. The collections of
endomorphisms $\j$ induce collections of elements of
$\End_W(L,F)$ and of $V'_p$, which we will denote by the same
letters. Let $M$ be the $\Z_p$-submodule of $N_L$ spanned by the
components $j_1,j_2,j_3,j_4$ of $\j$. We have the following
commutative diagram:
$$\matrix
\{j_1,\dots,j_4\}&\subset& C_\xi&\lra&\End_W(L,F)&\supset& N_L&\supset&M \\
\nass
\downarrow&{}&\downarrow&{}&\downarrow&{}&\downarrow&{}&{}\\
\nass
V'&\subset&C'&\lra& C'_p&\supset& V'_p&{}&{}&.\endmatrix\tag5.14$$
Recall that $T=T_\xi\in {\roman{Sym}}_4(\Z_{(p)})\subset
{\roman{Sym}}_4(\Z_p)$ is the matrix of inner products of the
elements $j_1,\dots,j_4$ with respect to the quadratic form on
$V'_p$. Thus we have the following basic observation:
\proclaim{Lemma 5.5} At a point of intersection
$\xi\in {\Cal Z}(T,\omega)\cap \Cal M^{\text{ss}}(\F)$ with
corresponding lattice $L\in X$, the matrix $T_\xi=T$ is
represented by the lattice $N_L
=\End_W(L,F)\cap V'_p$ in the quadratic space $V'_p$. In fact,
$T$ is the matrix for the restriction of the quadratic form on
$N_L$ to the sublattice $M$ spanned by $j_1,\dots,j_4$.
\endproclaim

Suppose $L\in X\setminus X_0$ with associated distinguished
lattice $\tl$. Recall that ${\Cal O}_L\subset {\Cal O}_{\tilde
L}$ and let $\om$ be the $\Z_p$-subalgebra of
$\otl=\End_W(\tl,F)$ generated by $j_1,\dots,j_4$, i.e., by $M$.
Also let $C(M)$ be the Clifford algebra of $M$. Let $\fn_L
= {\roman{red}}_\tl(N_L)$ and let $\fm=\fm_{\tilde L}={\roman{red}}_\tl(M)$, so
that
$$\matrix M&\subset&N_L&\subset&N_{\tilde L}\\
\nass
\downarrow&{}&\downarrow&{}&\downarrow\\
\nass
\fm_{\tilde L}&\subset&\fn_{ L}&\subset&\tfn.\endmatrix\tag5.15$$

\proclaim{Lemma 5.6} Suppose that $L\in X\setminus X_0$ with associated
distinguished
lattice $\tl$. \hfill\break (i) The natural map $C(M)\rightarrow
\om$ is an isomorphism. \hfill\break (ii) There is a commutative
diagram
$$
\matrix \om&\subset &\otl\\
\nass
\downarrow&{}&\downarrow\\
\nass
{\roman{red}}_\tl(\om)&{}&{\roman{red}}_\tl(\otl)\simeq M_2({\Bbb
F}_{p^2})\\
\nass
||&{}&\cup\\
\nass
C(\fm_\tl)&\hookrightarrow&C(\tfn)={\Bbb F}_{p^2}\cdot 1_2,\endmatrix$$
\endproclaim
\demo{Proof}
The inclusion of $C(\fm_\tl)$ into $C(\tfn)$ is induced by the
inclusion of quadratic spaces $\fm\subset \tfn$. We obviously
have a commutative diagram with surjective vertical arrows,
$$
\matrix C(M)&\longrightarrow &{\Cal O}_M\\
\nass
\downarrow&{}&\downarrow\\
\nass
C(\fm)&\longrightarrow&{\roman{red}_\tl}({\Cal O}_M).\\
\endmatrix\tag5.16$$
But the upper horizontal arrow is surjective since both algebras
are generated by $M$. This proves that the lower horizontal arrow
is surjective. By the statement at the beginning it is also
injective which proves the equality sign at the south-west corner
of the diagram in (ii). The rest of the Lemma follows from Lemma 5.3.
\qed\enddemo

Next let us consider the case when $L\in X_0$. We use somewhat
similar notation: let $\fn_L= {\roman{red}}_L(N_L)$ and
$\fm=\fm_L={\roman{red}}_L({\Cal O}_M)$.

The same arguments yield:

\proclaim{Lemma 5.7} Suppose that $L\in X_0$ is superspecial.
There is a commutative diagram:
$$
\matrix \om&\subset &\ol\\
\nass
\downarrow&{}&\downarrow\\
\nass
{\roman{red}}_L(\om)&\subset&{\roman{red}}_L(\ol)\simeq M_2({\Bbb
F}_{p^2})\\
\nass
||&{}&||\\
\nass
C(\fm_L)&{}&C(\fn_L)&.\endmatrix$$
\endproclaim
Our next task will be to show that the matrix $T\ {\roman{mod}}\,
p$ in $M_4( {\Bbb F}_p)$ controls the size of $\fm$. More
precisely, we now list the possibilities for ${\frak m}$ in the
non-superspecial and the superspecial case separately.
\proclaim{Lemma 5.8}
Let $L\in X\setminus X_0$ with associated $\tilde{L}\in \tilde X$.
The possibilities for ${\frak m}={\frak m}_{\tilde L}$ are the
following:

\roster
\item"{(i)}" If $\dim_{{\Bbb F}_p}{\frak m}=1$, then $T$ has rank 1 modulo
$p$ and does not represent 1.
\item"{(ii)}" If ${\frak m}=0$, then $p\mid T$.
\endroster
\endproclaim

\demo{Proof}
The first alternative corresponds to the case where ${\frak
m}={\frak n}_{\tilde L}$, by Lemma 5.3. The assertion now follows from
Lemma 5.5.
\qed\enddemo

\proclaim{Lemma 5.9}
Let $L\in X_0$. The possibilities for $\fm=\fm_L$ and
$C(\fm)\subset M_2(\fps)$ are the following:
\roster
\item"{(i)}" The rank of $T\ {\roman{mod}}\ p$ is 3, or equivalently
$\dim\fm =3$. Then $C(\fm)\simeq M_2({\Bbb F}_{p^2})$.
\item"{(ii)}" The rank of $T\ {\roman{mod}}\ p$ is 2.
Then $\dim\fm =2$ and $C(\fm)\simeq M_2({\Bbb F}_p)$.
\item"{(iii)}" The rank of $T\ {\roman{mod}}\ p$ is 1 and $\dim\fm =2$.
Then $\fm$ is of the form $\fm =\fm_0+{\frak r}$ where $\fr$ is the
radical and $\dim \fm_0= \dim{\frak r} =1$. In this case
$$C(\fm)\simeq C(\fm_0)^\sim[\e]/(\e^2)$$
where
$$C(\fm_0) \simeq \cases \fp\oplus\fp & \\ \fps&\endcases,$$
and the element $\e$ acts on $C(\fm_0)$ by the nontrivial
automorphism of order $2$.
\item"{(iv)}" The rank of $T\ {\roman{mod}}\ p$ is 1 and $\dim\fm =1$. Then
$$C(\fm) \simeq \cases \fp\oplus\fp &
\\ \fps&\endcases.$$
\item"{(v)}" $T\equiv 0\ {\roman{mod}}\ p$ and $\dim \fm =1$. Then
$$C(\fm)=\wedge(\fm)\simeq {\Bbb F}_p[\e\,] / (\e^2)\ \ .$$
\item"{(vi)}" $\fm=0$. Then $T\equiv 0\ {\roman{mod}}\ p$
and $C(\fm)={\Bbb F}_p$ is in the center of $M_2({\Bbb
F}_{p^2})$.
\endroster
\endproclaim

In cases (iii) resp.\ (iv), $\fm_0$ resp.\ $\fm$ is a
nondegenerate line, so that the quadratic form on it is
isomorphic to either $x^2$ or $a x^2$, with $a\in {\Bbb
F}^\times_p\setminus {\Bbb F}^{\times,2}_p$, yielding a Clifford
algebra $\fp\oplus\fp$ or $\fps$.
\proclaim{Lemma 5.10} In cases (iii) and (iv) above, when an $\fps$ arises in
the Clifford algebra $C(\fm)$, this $\fps$ is not central in $M_2(\fps)$.
\endproclaim
\demo{Proof}
Choose $x\in \fm$ spanning a nondegenerate line for which $C(\fp
x) \simeq \fps$. Then $x$ is an endomorphism of the
$2$-dimensional $\fps$-vector space $\tl_0/F\tl_0$ with
${\roman{tr}}(x)=0$, and with $x^2=q(x)\cdot id$ where
$q(x)\notin {\Bbb F}^{\times,2}_p$. This last condition is
equivalent to our hypothesis on the Clifford algebra. Thus, $x$
has two distinct eigenvalues $\pm\sqrt{q(x)}$ on $\tl_0/F\tl_0$,
and hence does not lie in the center.
\qed\enddemo

We can now describe the intersections of
our special cycle with the supersingular locus.
\proclaim{Theorem 5.11}
Suppose that $\xi\in {\Cal Z}(T,\omega)$ with image in the supersingular locus
$\Cal M^{\text{ss}}(\F)$ with corresponding $L\in X$.\hfill\break
(i) The rank of $T=T_\xi$ modulo $p$ is at most $3$.\hfill\break (ii)
If $T_\xi$ represents 1, then $L\in X_0$ and $\xi$ is a point of
proper intersection.\hfill\break (iii) If $L\in X\setminus X_0$,
with associated distinguished lattice $\tl$, the whole
distinguished $\Bbb P_{\tilde L}$ associated to $\tl$ in the
supersingular locus, and passing through $\xi$, occurs in ${\Cal
Z}(T,\omega)$. In particular, $\xi$ is not a point of proper
intersection.
\endproclaim
\demo{Proof} The reduction of $T=T_\xi$ modulo $p$
is the matrix for the quadratic form on the images of
$j_1,\dots,j_4$ in $\fm =\fm_L$, resp.\ $\fm=\fm_{\tilde L}$, and
$\fm$ has dimension at most $3$. This proves (i). If $L\in
X\setminus X_0$, then $T=T_{\xi}$ does not represent $1$, by Lemma
5.8. Furthermore in this case by Lemma 5.3.\
$$C({\frak m})\subset \F_{p^2}\cdot 1\subset {\roman{red}}_{\tilde
L}({\Cal O}_{L'})\ \ ,$$ for any $L'\in {\Bbb P}_{\tilde L}$,
which is not superspecial. This implies $M\subset \Cal O_M\subset {\Cal
O}_{L'}$.
If now $L'\in {\Bbb P}_{\tilde L}$ is
superspecial it follows that $M\subset {\Cal O}_M\subset
{\Cal O}_{L'}$ by specialization.
\qed\enddemo

It remains to consider the cases where $T$ does not represent $1$. We first
treat the case when $p\mid T$.
\proclaim{Theorem 5.12} Suppose that $p\mid T$ and that
$\xi\in {\Cal Z}(T,\omega)$ has image in $\Cal M^{\text{ss}}(\F)$
with corresponding $L\in X_0$. Then $\xi$ is not a point of proper
intersection. More precisely:\hfill\break (i) If
$\fm={\roman{red}}_L(M)=0$ then each of the $p+1$ distinguished
$\Bbb P^1$'s through $pr(\xi)\in{\Cal M}^{\roman{ss}}$ occurs in
the image of ${\Cal Z}(T,\omega)$, i.e., for every distinguished
$\tl$ with $\tilde L\supset L\supset F\tilde L$, we have $M\subset
\End_W(\tl,F)$; furthermore ${\roman{red}}_\tl(M)=0$.\hfill\break
(ii) If $\fm={\roman{red}}_L(M)$ is a null line in $\fn_L$, then
there is a unique distinguished $\tl$ with $\tilde L\supset
L\supset F\tilde L$ and with $M\subset \End_W(\tl,F)$; furthermore
${\roman{red}}_\tl(M)=0$. Hence there is a unique distinguished
$\Bbb P^1$ passing through $pr(\xi)\in{\Cal M}^{\roman{ss}}$ and
contained in the image of${\Cal Z}(T,\omega)$.
\endproclaim
\demo{Proof} We may assume that $L^\bot =L$. First suppose that $\fm$ is a null
line, and choose
$x_0\in M$ such that $\bar{x}_0={\roman{red}}_L(x_0)$ spans
$\fm={\roman{red}}_L(M)$.
The endomorphism $\bar{x}_0$ of the two dimensional vector space
$L_0/FL_0$ satisfies $\bar{x}_0^2=0$ but $\bar{x}_0\ne 0$. Thus
$\im(\bar{x}_0)$ is a line in $L_0/FL_0$.
\proclaim{Lemma 5.13} Assume that $L=L^\perp$.
The lattice $\tilde L$ defined by $F\tl = x_0(L)+FL$ lies in $\tilde{X}$ and
$L\in {\Bbb P}_{\tilde L}$.
Moreover, $M\subset \End_W(\tl,F)$, and
${\roman{red}}_\tl(M)=0$.
\endproclaim
\demo{Proof} Since $x_0$ commutes with $F$, we clearly have
$F(F\tl)=x_0(FL)+F^2L = x_0(FL)+pL \subset F\tl$. Similarly one
sees that $V(F\tilde L)\subset F\tilde L$, i.e.,  $pF\tilde
L\subset F(F\tilde L)$, hence $F\tilde L$ is admissible.
Note that $(F\tl)^\perp\supset L^\perp=L \supset F\tl$ with
$[(F\tl)^\perp:L]=[L:F\tl]=1$.

To show that $F\tl$ is distinguished,
it will suffice to prove that $F(F\tl)\subset p(F\tl)^\perp$, i.e.,
that $<F^2\tl,F\tl>\, \subset pW$.
But
$$\align
<F^2\tl,F\tl>\  &=\ <x_0(FL) + pL,x_0(L)+FL> \\ {}&\subset\
<x_0(FL),x_0(L)> + pW\\ {}&=\ <FL,x_0^2(L)>+ pW\tag5.17\\
{}&\subset\ <FL,FL>+ pW\\ {}&\subset\ pW.\endalign$$ Here we have
used the fact that $x_0^*=x_0$ and that $\bar{x}_0^2=0$, i.e., that
$x_0^2(L)\subset FL$. We conclude that $F\tilde L\in \tilde X$ and
hence also $\tilde L\in\tilde X$.

Next, we must show that every element of $M$ preserves $\tl$ or, equivalently,
$F\tilde M$.
In fact, we show that
$M\cdot \tl\subset F\tl$, so that ${\roman{red}}_\tl(M)=0$. First consider the
reduction sequence
$$0\lra M_0\lra M\ \overset{{\roman{red}}_L}\to{\lra}\  \fp\cdot \bar{x}_0\lra
0,\tag5.18$$
where
$$M_0=\{ \ y\in M\mid y(L)\subset FL\ \}.\tag5.19$$
It suffices to prove the inclusions $x_0(\tl)\subset F\tl$ and
$y(\tl)\subset F\tl$ for all $y\in M_0$.
Recall that, for $x\in M$, $x^2=q(x)\cdot id$.
Since $p\mid T_\xi$, the resulting quadratic form on $M/pM$ is
identically zero, and so $C(M/pM)=\wedge(M/pM)$. In particular,
for any $x_1$ and $x_2\in M$, $x_1x_2\equiv -x_2x_1 \mod p$, i.e.,
$$x_1x_2(L) \subset x_2x_1(L)+pL.\tag5.20$$
Now, for $y\in M_0$,
$$\align
y(F\tl) &= yx_0(L)+y(FL)\\
{}&\subset x_0y(L) + pL + F(y(L))\tag5.21\\
{}&\subset x_0(FL) + F^2L\\
{}&\subset F(x_0(L)+FL) = F^2\tl.
\endalign$$
Next, observe that $x_0^2 = q(x_0)\cdot id$ and $q(x_0)\equiv 0\mod p$
implies that $x_0^2(L)\subset pL$, not just $FL$. Thus
$$\align
x_0(F\tl) &=x_0^2(L) + F x_0(L)\\
{}&\subset pL + F x_0(L)\tag5.22\\
{}&= F(FL + x_0(L)) = F^2\tl.\endalign$$
This completes the proof of the Lemma.
\qed\enddemo

To finish the proof of (ii), we show that the distinguished lattice
$F\tilde L$ constructed in Lemma~5.13 is unique. Note that
$\ker(\bar{x}_0)=\im(\bar{x}_0)$. If $\tl' = W\cdot u + FL$ is
another distinguished lattice, whose image $\ell' =\tl'/FL$ is
distinct from $\ker(\bar{x}_0)$, then
$$\bar{x}_0(\ell') = \im(\bar{x_0}) \ne \ell',\tag5.23$$
so that $\tl'$ is not preserved by $x_0$.

Now suppose that ${\roman{red}}_L(M)=0$, i.e., that $M\cdot L\subset FL$. Let
$F\tl = W\cdot u+FL$ be any distinguished lattice with $L\supset F\tl\supset
FL$.
We want to show that, for any $x\in M$, $x(\tl)\subset F\tl = p\tl^\perp$ or,
equivalently $x(F\tilde
L)\subset F^2\tilde L=p\cdot (F\tilde L)^{\bot}$.
But now
$$\align
<x(F\tl),F\tl>\ &=\ <W x(u)+Fx(L),W u+FL>\tag5.24\\
{}&\subset\ W<x(u),u> + pW.
\endalign$$
But now, since $x^*=x$,
$$<x(u),u>\  =\  <u,x(u)>\  =\  -<x(u),u>\tag5.25$$
so that $<x(u),u>=0$. Thus $<x(F\tl),F\tl>\ \subset pW$, i.e.,
$x(F\tl)\subset pF(\tilde L)^\bot = F^2\tl$, as required.
This concludes the proof of Theorem~5.12.
\qed\qed\enddemo

We now turn to the case when $p\nmid T$ but $T$ does not represent $1$.

\proclaim{Theorem 5.14} Suppose that $p\nmid T$ and that $T$ does not represent
$1$.
Let $\xi\in{\Cal Z}(T,\omega)$ with
$pr(\xi)\in {\Cal M}^{ss}({\Bbb F})$ and with corresponding $L\in X_0$.  Then
$\xi$ is
not a point of proper intersection. More precisely:
\hfill\break
(i) If $\dim_{{\Bbb F}_p}{\frak m}=1$ and ${\frak m}$
does not represent $1$, then precisely two of the $p+1$
distinguished ${\Bbb P}_{\tilde L}$'s through $pr(\xi)$ occur in the
image of $\Cal Z(T,\o)$. These are the only two distinguished lattices
$\tilde L_1$ and $\tilde L_2$ with $\tilde L_i\supset L\supset
F\tilde L_i$ and with $M\subset {\roman {End}}_W(\tilde L_i, F)$
$(i=1,2)$. Furthermore ${\roman {red}}_{\tilde L_i}(M)\ne (0)$,
$i=1,2$.
\hfill\break
(ii) If $\dim_{{\Bbb F}_p} {\frak m}=2$ and ${\frak m}= {\frak
m}_0+\fr$ as in Lemma 5.9, (iii), where ${\frak m}$ does not
represent $1$, then precisely one of the $p+1$ distinguished ${\Bbb
P}_{\tilde L}$'s  through $pr(\xi)$ occurs in $pr(\Cal Z(T,\o))$.
It is the only distinguished lattice $\tilde L_1$ with $\tilde
L_1\supset L\supset F\tilde L_1$ and with $M\subset
{\roman{End}}_W(\tilde L_1, F)$. Furthermore ${\roman{red}}_{\tilde
L_1}(M)\ne (0)$.
\endproclaim
\demo{Proof}
We may again assume $L^{\bot}=L$. We first consider case (i). Choose
$x_0\in M$ such that $\overline x_0={\roman{red}}_L(x_0)$ spans ${\frak m}
={\roman{red}}_L(M)$. Then $\overline x_0$ is an automorphism of $L_0/FL_0$
with
$$\overline x_0^2=\varepsilon\cdot 1\ \ ,$$
where $\varepsilon\in {\Bbb F}_p^\times\setminus {\Bbb F}_p^{\times,2}$.
Furthermore,
by Lemma 5.8., $\overline x_0$ is not central. Hence $\overline x_0$ has two
distinct eigenvalues $\varepsilon_1= \sqrt{\varepsilon}$ and
$\varepsilon_2=-\sqrt{\varepsilon}$ in ${\Bbb F}_{p^2}$. Let $E_1$ and $E_2$
be the corresponding eigenspaces and let $F\tilde L_1$ and $F\tilde L_2$ be
the corresponding lattices in ${\Cal L}$,
$$FL\subset F\tilde L_i\subset L\ \ , \ \ i=1,2\ \ .$$
Since $\overline x_0$ commutes with $F$ and $V$, the eigenspaces $E_1$ and
$E_2$ are preserved by $F$ and $V$, hence $F\tilde L_1$ and $F\tilde L_2$
are admissible. To see that $F\tilde L_1$ and $F\tilde L_2$ are
distinguished, it suffices to see that $F\tilde L_i\subset p\cdot \tilde
L_i^\bot$, i.e., that
$$\langle F\tilde L_i, \tilde L_i\rangle\subset p\cdot W\ \ ,\ \ i=1,2\ \ .$$
Equivalently we have to see that the eigenspaces $E_1$ and $E_2$ are
isotropic with respect to the antihermitian form (5.11) on $L_0/FL_0$.
If $v\in E_i$, then $\overline x_0v=\varepsilon_i
\cdot v$ and
$$
\varepsilon\cdot (v,v)
=(\varepsilon v,v)=(\overline x_0^2v,v)= (\overline x_0v, \overline x_0v)
=(\varepsilon_iv, \varepsilon_iv)
=-\varepsilon\cdot (v,v).
$$
It follows that $F\tilde L_1$ and $F\tilde L_2$ are distinguished.

The
lattices $\tilde L_i=F^{-1}(F\tilde L_i)\in \tilde X$ for $i=1$ and $2$ are
the distinguished lattices appearing in the statement of (i). We have
${\roman{red}}_{\tilde L_i}(M)\ne 0$ since $\overline x_0$ induces an
automorphism of the eigenspace $E_i$. On the other hand any $y\in M$ with
${\roman{red}}_L(y)=0$, i.e., with $y(L)\subset FL$, also satisfies
$y(\tilde L_i)\subset L\subset \tilde L_i$. It follows that $M\subset
{\roman{End}}_W(\tilde L_i, F)$, hence (i).

Now we consider case (ii). Let $\overline x_0\in{\frak m}_0$ with
$\overline x_0^2=\varepsilon\cdot 1$, for  $\varepsilon\in {\Bbb
F}_p^\times\setminus {\Bbb F}_p^{\times,2}$,  and let $\overline y_0$
be a generator of the
radical ${\frak r}$. Then $\overline x_0\overline y_0= -\overline y_0
\overline x_0$.
Therefore $\overline y_0$ maps the eigenspace $E_1$ of
$\overline{x}_0$ in $L_0/FL_0$ to the eigenspace $E_2$ and the
eigenspace $E_2$ to $E_1$. Since $\overline y_0^2=0$, but
$\overline y_0\ne 0$, precisely one of the two eigenspaces is
annihilated by $\overline y_0$. The corresponding lattice is
distinguished and yields as in case (i) the lattice $\tilde L_1$
appearing in the statement of (ii).
\qed\enddemo

\proclaim{Corollary 5.15}  Let
$\xi\in \Cal Z(T,\o)\subset
{\Cal Z}(d_1,\omega_1)\times_{\Cal M}\dots\times_{\Cal M}{\Cal Z}(d_r,
\omega_r),$
where $d_i\in {\roman{Sym}}_{n_i}({\Bbb Q})_{>0}$
with $n_1+\dots +n_r=4$, cf.\ (3.1) and (3.6).
Then $\xi$ is a point of
proper intersection if and only if its fundamental matrix
$T=T_{\xi}$ is non-singular and represents 1 over ${\Bbb Z}_p$. In
this case $\xi$ is supersingular and superspecial.
\endproclaim

A topic we have not touched upon in the present paper is to
describe the shape of the intersection of our cycles in the case of
improper intersection, or, equivalently, to describe, for $T\in
{\roman{Sym}}_4({\Bbb Q})_{>0}$, the cycle ${\Cal Z}(T,\omega)$
when its dimension is positive. We refer to the companion paper
\cite{\krhb} to the present one for more information on this topic.

\vfill\eject
\subheading{\Sec6. Intersection multiplicities}

In this section we consider the intersection multiplicity at a
point of proper intersection. More precisely we return to the setup
of the third section, i.e., we fix a decomposition $4=n_1+\dots
+n_r$, where $n_i\ge 1$ for all $i$, elements $d_i\in
{\roman{Sym}}_{n_i} ({\Bbb Q})_{>0}$ and
$\omega_i\subset V({\Bbb A}_f^p)^{n_i}$
giving rise to special cycles ${\Cal Z}(d_1,\omega_1),
\ldots, {\Cal Z}(d_r,\omega_r)$. We fix a point $\xi\in {\Cal Z}(d_1,\omega_1)
\times_{\Cal M}\ldots \times_{\Cal M} {\Cal Z}(d_r,\omega_r)$ with
${\roman{det}}(T_{\xi})\ne 0$ and where $T=T_{\xi}$ represents 1 over
${\Bbb Z}_p$. Let $\xi$ correspond to $(A,\iota, \lambda,
\overline{\eta}^p; {\j}_1,\ldots, \j_r)$. Since
${\roman{det}}(T)\ne 0$ and since $T$ represents 1 over
${\Bbb Z}_p$, the associated Dieudonn\'e module $L$ is superspecial
and corresponds to a formal group ${\Cal A}$ of dimension 2 and
height 4, with a collection of endomorphisms $\j =(j_1, \ldots,j_4)$
spanning a ${\Bbb Z}_p$-submodule $M$ of rank 4 in
${\roman{End}}_W(L,F)$. By changing the trivialization of the
rational Dieudonn\'e module we may assume that $L=L^\bot$, i.e.,
that ${\Cal A}$ is equipped with a principal quasi-polarization
$\lambda_{\Cal A}$. An immediate application of the theorem of
Serre and Tate shows that the formal completions at $\xi$ of ${\Cal
M}$ and of its closed subschemes may be interpreted as versal
deformation spaces, i.e., with obvious notation,
$$\align
\hat{\Cal M}_{\xi}
&={\roman{Def}}({\Cal A}, \lambda_{\Cal A})\tag6.1\\
\hat{\Cal Z}(d_i,\omega_i)_{\xi}
&={\roman{Def}}({\Cal A}, \lambda_{\Cal A}; \j_i)\tag6.2\\
\hat{\Cal Z}(T,\o)_\xi&=\big({\Cal Z}(d_1,\omega_1)\times_{\Cal
M}\dots\times_{\Cal M} {\Cal
Z}(d_r,\omega_r)\big)_{\xi}\hskip-10pt\hat{\phantom{M}}\tag6.3\\
{}&={\roman{Def}}({\Cal A}, \lambda_{\Cal A}; \j )=
{\roman{Def}}({\Cal A}, \lambda_{\Cal A}; M).
\endalign$$
Here $\o=\o_1\times\cdots\times\o_r$. Recall that any $x\in M$
satisfies $x^*=x$ and $tr(x)=0$, that the quadratic form on $M$ is
given by $x^2=q(x)\cdot id$, and that $T$ is the matrix for that
quadratic form with respect to the basis $j_1,\dots,j_4$. Since $T$
represents 1 over ${\Bbb Z}_p$, there exists an $x_0\in M$ such
that $x_0^2 = id$. The quadratic lattice $M$ can be written as
$M=\Z_p\cdot x_0 + M_0$, where $M_0 = x_0^\perp$. Moreover, if
$x\in M_0$, then $x x_0 = -x_0 x$, since the subalgebra of
$\End_W(L,F)$ generated by $M$ is the image of $C(M)$, the Clifford
algebra of $M$.
\par
The idempotents $e_1=\frac12(1+x_0)$, $e_2=\frac12(1-x_0)$ -- recall that $p\ne
2$ --
give a splitting $\Cal A \simeq \Cal A_1\times\Cal A_2$, with $\Cal A_1=e_1\Cal
A$
and $\Cal A_2= e_2\Cal A$ of dimension $1$ and
height $2$. If $x\in M_0$,  $x e_1 = e_2 x$, and so $M_0$ can be viewed as a
submodule of $Hom(\Cal A_1,\Cal A_2)$. Let $L_i= e_iL$ be the Dieudonn\'e
module of ${\Cal A}_i$. Then
$L=L_1\oplus L_2$. Furthermore $L_1$ and $L_2$ are paired trivially under
the symplectic pairing on $L$.
Indeed, if $v_1\in L_1$ and $v_2\in L_2$, we have
$$< v_1,v_2>\ =\ < e_1v_1, e_2v_2>\ =\ < v_1, e_1e_2v_2>\ =0\ \ ,$$
since $e_1^\ast =e_1$. It follows that $<\ ,\ >$ induces a perfect symplectic
pairing
$<\ ,\ >_i$  on $L_i$,  i.e., ${\Cal A}_1$ and ${\Cal A}_2$ are equipped with
principal
quasi-polarizations. Since a principal quasi-polarization on a $p$-divisible
formal group of dimension 1 and height
2 deforms automatically we obtain a natural identification
$${\roman{Def}}({\Cal A}, \lambda_{\Cal A};M)={\roman{Def}} ({\Cal A}_1, {\Cal
A}_2;M_0)\ \ .$$
The length $e(\xi)$ of the local Artin ring appearing on the right was
determined
by
Gross and Keating in section 5 of \cite{\grosskeating}. Since we have assumed
in all of the above that $p\ne 2$, we may as well continue to make this
assumption,
although Gross and Keating do not.
Choose a basis $\psi_1,\ \psi_2,\ \psi_3$ for $M_0$ such that
$$q(u_1\psi_1+u_2\psi_2+u_3\psi_3) = \e_1p^{a_1}u_1^2
+\e_2p^{a_2}u_2^2+\e_3p^{a_3}u_3^2,$$
with $0\le a_1\le a_2\le a_3$. Thus, over $\Z_p$, $T$ is equivalent to
the diagonal matrix $\text{diag}(1,\e_1p^{a_1},\e_2p^{a_2},\e_3p^{a_3})$.

\proclaim{Proposition 6.1} {\rm(Gross, Keating, \cite{\grosskeating},
Proposition~5.4)}
If $a_1+ a_2$ is even,
then
$e(\xi)=e_p(T_\xi)$ is equal to:
$$\multline
\sum\limits_{i=0}^{a_1-1} (i+1)(a_1+a_2+a_3 -3i) p^i
+\sum\limits_{i=a_1}^{(a_1+a_2-2)/2}
(a_1+1)(2a_1+a_2+a_3-4i)p^i\\
+\frac12 (a_1+1)(a_3-a_2+1)p^{(a_1+a_2)/2}.\qquad\qquad\qquad
\endmultline$$
If $a_1+ a_2$ is odd,
then
$e(\xi)=e_p(T_\xi)$ is equal to:
$$
\sum\limits_{i=0}^{a_1-1} (i+1)(a_1+a_2+a_3 -3i) p^i
+\sum\limits_{i=a_1}^{(a_1+a_2-1)/2}
(a_1+1)(2a_1+a_2+a_3-4i)p^i.$$
\endproclaim
\proclaim{Corollary 6.2}
The cycles ${\Cal Z}(d_1,\omega_1), \ldots, {\Cal Z}(d_r,\omega_r)$
intersect transversally at the point $\xi$ if and only if
${\roman{ord}} ({\roman{det}}\ T_{\xi})=1$.
\endproclaim
\demo{Proof} The above formulas show that $e(\xi)=1$ if and only if
${\roman{ord}}({\roman{det}}\, T_{\xi})=1$. If this is the case,
the cycles ${\Cal Z}(d_i, \omega_i)$ have to be irreducible and
reduced locally at $\xi$, and the intersection multiplicity in the
sense of Serre, which is bounded by the length, is equal to 1. In
this case the cycles ${\Cal Z}(d_i, \omega_i)$ are all regular at
$\xi$ and their tangent spaces give a direct sum decomposition of
the tangent space of ${\Cal M}$ at $\xi$. For all this, cf.\
\cite{\fulton}, Prop.\ 8.2, and Example 8.2.1.
\enddemo
At this point we have completely answered the question a) at the
end of section 3. What is not clear is whether the length $e(\xi)$
is indeed the intersection multiplicity of ${\Cal Z}(d_1,\omega_1),
\ldots, {\Cal Z}(d_r,\omega_r)$ at $\xi$. This is the content of
the question b) of section 3.
\proclaim{Conjecture 6.3}
Let $\xi$ be an isolated intersection point of
${\Cal Z}(d_1,\omega_1),\ldots, {\Cal Z}(d_r,\omega_r)$.
Then
$$({\Cal O}_{{\Cal Z}(d_1,\omega_1)} \mathop{\otimes}\limits^{\Bbb L}\ldots
\mathop{\otimes}\limits^{\Bbb L}{\Cal O}_{{\Cal Z}(d_r,\omega_r)})_{\xi} =
({\Cal O}_{{\Cal Z}(d_1,\omega_1)} \mathop{\otimes}
\ldots\otimes {\Cal O}_{{\Cal Z}(d_r,\omega_r)})_{\xi}\ \ ,$$
hence $e(\xi)$ is the intersection multiplicity of
${\Cal Z}(d_1,\omega_1), \ldots, {\Cal Z}(d_r,\omega_r)$ at $\xi$.
\endproclaim
We stress that this conjecture is reasonable only because ${\Cal
M}$ is smooth over ${\roman{Spec}}\, {\Bbb Z}_{(p)}$. Indeed,
Genestier \cite{\genestier} has shown that in the
Drinfeld-Cherednik situation of bad reduction the analogues of the
special cycles considered here may have embedded components. On the
other hand, assume in our situation that ${\Cal Z}(d_i,\omega_i)$
is an intersection of $n_i$ divisors in ${\Cal M}$. Then if $\xi$
is an isolated intersection point of ${\Cal
Z}(d_1,\omega_1),\ldots, {\Cal Z}(d_r,\omega_r)$ it follows that
each partial intersection ${\Cal Z}(\underline
d_{i_1},\omega_{i_1})\cap\ldots\cap {\Cal Z}(d_{i_s},\omega_{i_s})$
($1\le i_1\le\ldots\le i_s\le r$) is locally at $\xi$ a complete
intersection. Hence it also follows that the length $e(\xi)$ is the
intersection multiplicity of ${\Cal Z}(d_1,\omega_1),\ldots, {\Cal
Z}(d_r,\omega_r)$ at $\xi$, and the above conjecture holds true.

{\bf Remark 6.4.}
Assume that $\xi\in {\Cal Z}(d_1,\omega_1)\cap\dots\cap {\Cal
Z}(d_r, \omega_r)$ is a point with fundamental matrix $T=T_{\xi}$
which is non-singular and represents over ${\Bbb Z}_p$ a unit
$\varepsilon\in {\Bbb Z}^\times_p \setminus {\Bbb Z}^{\times,
2}_p$. Therefore there exists $x_1\in M$ such that
$x_1^2=\varepsilon \cdot {\roman{id}}$. Hence we obtain an action
of ${\Bbb Z}_{p^2} ={\Bbb Z}_p[\sqrt{\varepsilon}]$ on ${\Cal
A}$,
$$\alpha: {\Bbb Z}_{p^2}\longrightarrow {\roman{End}}({\Cal A})\ \
.$$
We may write $M$ in the form $M={\Bbb Z}_p\cdot x_1+M_1$, where
$M_1=x_1^\bot$. For $x\in M_1$ we have $xx_1=-x_1x$. Comparing with
the definitions in the companion paper to this one, we see that
$({\Cal A}, \alpha)$ is precisely one of the formal groups with
${\Bbb Z}_{p^2}$-action considered there \cite{} and that the
elements of $M_1$ are {\it special endomorphisms}\/ in the sense of
that paper. In particular, the formal completion of ${\Cal
Z}(d_1,\omega_1)\cap\dots\cap {\Cal Z}(d_r, \omega_r)$ at $\xi$
coincides with the formal completion of the corresponding
subvariety of the Hilbert-Blumenthal surface considered in
\cite{\krhb}.

\vfill\eject
\subheading{\Sec7. The total contribution of isolated points}

In this section we will consider the total contribution of the points
of proper intersection of our special
cycles. Using our previous results and a counting argument, we are able to give
an
explicit formula.

We return to the global situation of sections 1 and 2 and fix data
as follows. We assume as always that $p\nmid 2 D(B)$ and that
$K=K^p\cdot K_p$ where $K_p$ is the standard maximal compact
subgroup ( see the end of section 4),  and where $K^p$ is neat. We
then have the moduli scheme ${\Cal M}={\Cal M}_{K^p}$ which is
smooth over ${\roman{Spec}}\, {\Bbb Z}_{(p)}$. As in section 3, we
fix $n_1,\ldots, n_r$ with $1\leq n_i\leq 4$ and with $n_1+\ldots
+n_r=4$. For $i=1, \ldots, r$, choose positive definite matrices
$d_i\in {\roman{Sym}}_{n_i}({\Bbb Z}_{(p)})_{>0}$ and $K^p$-invariant
open compact subsets $\omega_i\subset V({\Bbb A}_f^p)^{n_i}$. We
then have the cycles
${\Cal Z}(d_i,\omega_i)$, $i=1,\ldots, r$.
We then define the  contribution of the points of proper
intersection to the intersection number of ${\Cal
Z}(d_1,\omega_1),\ldots, {\Cal Z}(d_r,\omega_r)$ to be
$$<{\Cal Z}(d_1,\omega_1),\ldots, {\Cal Z}(d_r,
\omega_r)>_p^{\roman{proper}} \ :=\sum_{\xi} e(\xi)\ \ .\tag7.1$$
Here the sum runs over the points of proper intersection $\xi$ in
${\Cal Z}(d_1,\o_1)\times_{\Cal M}\ldots\times_{\Cal M} {\Cal Z}(d_r,\o_r)$,
and $e(\xi)$ denotes
the length of the local ring at $\xi$, as described in section 6.
Note that, if Conjecture~6.3 were
known to hold, this is also the local intersection multiplicity at
$\xi$.

In the special case $r=1$, we let $d_1=T$, and we have the cycle
${\Cal Z}(T,\omega)$, whose image in $\Cal M$ lies in the
supersingular locus ${\Cal M}^{ss}$. Then ${\Cal Z}(T,\omega)$ is a
collection of isolated points if and only if $T$ represents 1 over
${\Bbb Z}_p$ (Corollary~5.15). In this case we use the notation
$$<{\Cal Z}(T,\omega)>_p\ = \sum_{\xi\in {\Cal Z}(T,\omega)} e(\xi)\ .\tag7.2$$

In general, by (3.6) and the analysis of the previous sections, we may write
$$
<{\Cal Z}(d_1,\omega_1),\ldots, {\Cal Z}(d_r, \omega_r)>_p^{\roman{proper}}
{}= \sum_T <{\Cal Z}(T,\omega)>_p,
\tag7.3$$
where the summation is over $T\in {\roman{Sym}}_4({\Bbb Z}_{(p)})_{>0}$
which are nonsingular,  represent 1 over ${\Bbb Z}_p$, and have
diagonal blocks $d_1,\ldots, d_r$:
$$T=\pmatrix d_1&{}&\dots&{}\\
{}&d_2&\dots&{}\\
\vdots&\vdots&\ddots&{}\\
{}&{}&\dots&d_r\endpmatrix.$$

We will now give more explicit expressions for the above entities. For
this it will suffice to give an expression for (7.2). But the
results of section 6 show that the intersection multiplicities
$e(\xi)$ in the sum of (7.2) only depend on $T$ and even only on
its ${\Bbb Z}_p$-equivalence class. As in Proposition~6.1,
we denote this integer by
$e_p(T)$ and thus may write
$$<{\Cal Z}(T,\omega)>_p\ = e_p(T)\cdot \vert {\Cal Z}(T, \omega)({\Bbb
F})\vert\ .\tag7.4$$
It remains to determine the cardinality of ${\Cal Z}(T,\omega)(\F)$.

As before, let $B'$ be the definite quaternion algebra with
discriminant $D(B)p$, let $C'=M_2(B')$, and let $V'=\{x\in C'\mid
x'=x\text{ and } tr(x)=0\ \}$. Let $G'$ be as in (4.5). Recall that
we also have fixed an isomorphism $G'(\A_f^p)\simeq G(\A_f^p)$, and
a base point $\xi_o= (A_o,\iota_o,\lambda_o,\overline{\eta}_o^p)\in
\Cal M^{\text{ss}}(\F)$ such that the associated Dieudonn\'e module
$L_o\in X$ is superspecial, with stabilizer $K'_p$ in $G'(\Q_p)$.
Then, under the parametrization (4.7), the set of superspecial
points in $\Cal M^{\text{ss}}(\F)$ corresponds to the double coset
space
$$G'(\Q)\back \bigg( G'(\Q_p)/K'_p\times G(\A_f^p)/K^p\bigg),\tag7.5$$
cf.\ Corollary~4.15. For a superspecial point
$(A,\iota,\lambda,\bar{\eta}^p)$ of $\Cal M^{\text{ss}}(\F)$, the
choice an isogeny $\gamma:(A,\iota)\rightarrow (A_o,\iota_o)$
compatible with the polarizations determines a pair $(g^p,g_p)\in
G'(\Q_p)/K'_p\times G(\A_f^p)/K^p$, and the passage to
$G'(\Q)$-orbits removes the dependence on the choice of $\gamma$.

The choice of an isogeny $\gamma$ also yields an identification of
the space $\End^0(A,\iota)^{\text{op}}$ with
$\End^0(A_o,\iota_o)^{\text{op}}= C'$, and of the space of special
endomorphisms of $(A,\iota,\lambda)$ with $V'(\Q)$. Let
$\Omega'_T({\Bbb Q})\subset V'(\Q)^4$ be the fibre over $T$ of the
map defined by the quadratic form on $V'({\Bbb Q})$,
$$V'({\Bbb Q})^4\longrightarrow {\roman{Sym}}_4({\Bbb Q})\ .\tag7.6$$

Returning to the set ${\Cal Z}(T,\omega)(\F)$,
we consider the map
$${\Cal Z}(T,\omega)(\F)\hookrightarrow
G'(\Q)\back \bigg( \Omega'_T({\Bbb Q})\times G'(\Q_p)/K'_p\times
G(\A_f^p)/K^p\bigg)\tag7.7$$ defined as follows. To a point
$\xi=(A,\iota, \lambda, \overline{\eta}^p;\j)
\in {\Cal Z}(T,\omega)(\F)$, and a choice of
isogeny $\gamma:(A,\iota)\rightarrow (A_o,\iota_o)$, there is an
associated triple  $(\gamma_*\j,g^p,g_p)$, where $\gamma_*\j\in
V'(\Q)^4$ is the 4-tuple of endomorphisms determined by $\j$ and
$\gamma$. Again, the passage to $G'(\Q)$-orbits removes the
dependence on the choice of $\gamma$.

It is not difficult to describe the image of ${\Cal Z}(T,\omega)(\F)$.
For $\y\in \Omega'_T({\Bbb Q})$, the triple $(\y,g^p,g_p)$ lies in the image
if and only if
\roster
\item"(i)" The images of the components of $\y$ under the inclusion
$V'\hookrightarrow \End_W({\Cal L}, F)$ preserve the lattice
$g_p\overline{L}$, and
\item"(ii)" The image of $\y$ under $\eta^p$ lies in $g^p\cdot\omega$.
\endroster
We note that the condition (i) is equivalent to the assertion that
the components of the 4-tuple $g_p^{-1}\y$ lie in
$$V'({\Bbb Z}_p) =V'({\Bbb Q}_p)\cap \End_W(L,F)=N_L\ .\tag7.8$$
We let $\varphi'_p$ be the characteristic function of
$V'({\Bbb Z}_p)^4$, let $\varphi_f^p={\roman{char}}(\omega)$ be the
characteristic function of $\o$,
and set $\varphi_f'=\varphi'_p\otimes \varphi_f^p$. Then
$\varphi_f'\in S(V'({\Bbb A}_f)^4)^{K'}$.
Conditions (i) and (ii) can then be summarized as follows.
\proclaim{Lemma 7.1}
The $G'(\Q)$-orbit of the triple $(\y,g^p,g_p)$ lies in the image
of ${\Cal Z}(T,\omega)(\F)$ if and only if $\varphi_f'(g^{-1}\y)\ne
0$, where $g = (g^p,g_p)\in G'(\A_f)$.
\endproclaim

Note that the function $(\y,g)\mapsto \ph'_f(g^{-1}\y)$ is
invariant under the diagonal action of $G'(\Q)$ on the left and
under the action of $K'=K'_pK^p$ and of $Z'(\A_f)$ on the right.
The total contribution of the superspecial points may be expressed
as an integral.

\proclaim{Theorem 7.2}
Let $T\in {\roman{Sym}}_4({\Bbb Z}_{(p)})_{>0}$ be non-singular and such
that $T$ represents 1 over ${\Bbb Z}_p$. Let $\omega\subset V({\Bbb
A}_f^p)^4$ be $K^p$-invariant open and compact,
and let $K'=K'_pK^p\subset G'(\A_f)$. Let $\pr(K')$ be the image of $K'$
in $Z'(\A_f)\back G'(\A_f) \simeq SO(V')(\A_f)$. Then
$$<{\Cal Z}(T,\omega)>_p\ =e_p(T)\cdot {\roman{vol}}(\pr(K'))^{-1}
\cdot I_{T,f}(\varphi_f')\ .$$
Here $\varphi_f'=\varphi'_p\otimes \varphi_f^p\in S(V({\Bbb
A}_f)^4)$ as above, and $I_{T,f}(\varphi_f')$ denotes the theta integral
$$I_{T,f}(\varphi_f')= \int\limits_{G'({\Bbb Q})Z'(\A_f)\back
G'({\Bbb A}_f)} \sum_{\y\in \Omega'_T({\Bbb Q})}
\varphi'_f(g^{-1}\y)\,dg\ .$$
\endproclaim
The measure $dg$ is induced by an arbitrary Haar measure on
$Z'(\A_f)\back G'({\Bbb A}_f)$ and the atomic measure on
$Z'(\Q)\back G'({\Bbb Q})$. The coefficient $e_p(T)$ is given by
the formulas in Proposition~6.1. The identity of the Theorem
remains valid if $T$ is nonsingular but not positive definite,
since, in that case, $T$ is not represented by $V'$, and hence both
sides of the identity vanish.

\demo{Proof} By Lemma~7.1, we see that
$$|{\Cal Z}(T,\omega)(\F)| = \sum_{ G'(\Q)\back\big( \Omega'_T(\Q) \times
G'(\A_f)/K'Z'(\A_f)\big)}
\ph'_f(g^{-1}\y).\tag7.9$$
On the other hand, since $\pr(K')$ is neat,
the stabilizer in $Z'(\Q)\back G'(\Q)$ of a coset $gK'Z'(\A_f)/Z'(\A_f)$ is
trivial.
Thus, we have
$$|{\Cal Z}(T,\omega)(\F)|
= \vol(\text{pr}(K'))^{-1} \int\limits_{G'(\Q)Z'(\A_f)\back G'(\A_f)}
\sum_{\y\in \O_T'(\Q)} \ph'_f(g^{-1}\y)\, dg,\tag 7.10$$
for a measure as described in the Theorem.
In combination with (7.4),
this gives the claimed expression.
\qed\enddemo

\proclaim{Corollary 7.3}
In the situation of the beginning of this section,
$$<{\Cal Z}(d_1,\omega_1),\ldots, {\Cal Z}(d_r, \omega_r)>_p^{\roman{proper}}\
=
\sum_T e_p(T)\ {\roman{vol}}(\pr(K'))^{-1}\cdot I_{T,f}(\varphi_f')\ .$$
where $T$ runs over all $T\in {\roman{Sym}}_4({\Bbb Z}_{(p)})_{>0}$
which represent 1 over ${\Bbb Z}_p$ and have
diagonal blocks $d_1,\ldots, d_r$. The function
$\varphi_f'=\varphi'_p\otimes \varphi_f^{p}\in S(V'({\Bbb A}_f)^4)$ is defined
by
$$\align
\varphi'_p
&
={\roman{char}}\ V'({\Bbb Z}_p)^4\\
\varphi_f^{p}
&
={\roman{char}}(\omega_1\times \ldots\times \omega_r).
\endalign$$
\endproclaim

\noindent{\bf Remark 7.4.} Formula (7.10) expresses the quantity $|{\Cal
Z}(T,\omega)(\F)| $ as a product of orbital integrals. More
precisely, note that the components of $\y\in \O_T'(\Q)$ span a
$4$-dimensional subspace of the $5$-dimensional space $V'$. Since
$G'$ acts on $V'$ via its projection to $SO(V')$, the stabilizer of
$\y$ in $G'(\Q)$ is precisely $Z'(\Q)$, the kernel of this
projection. Since $G'(\Q)$ acts transitively on $\O_T'(\Q)$, we can
unfold to obtain:
$$\align
|{\Cal Z}(T,\omega)(\F)| &= \vol(K')^{-1} \int_{Z'(\Q)\back G'(\A_f)}
\ph'_f(g^{-1}\y)\ dg\tag7.11\\
\nass
{}&=\vol(K')^{-1} \vol(Z'(\Q)\back Z'(\A_f))\, O_{\y}(\ph^p)\,O_{\y}(\ph_p'),
\endalign$$
for orbital integrals
$$O_{\y}(\ph^p) = \int_{Z(\A_f^p)\back G(\A_f^p)} \ph_f^p(g^{-1}\y)\,
dg,\tag7.12$$
and
$$O_{\y}(\ph_p') = \int_{Z'(\Q_p)\back G'(\Q_p)} \ph'_p(g^{-1}\y)\,
dg.\tag7.13$$

In our main theorem (in section 9), we will identify the right hand
sides of the formulas of Theorem~7.2 and Corollary~7.3 as special values of
derivatives of
Fourier coefficients of certain Eisenstein series. In the next section we will
explain more
precisely the Eisenstein series in question.

\vfill\eject
\define\stuff{\prime}
\define\Lt{\theta_T}
\define\Ltf{\theta_{T,f}}
\define\ev{\text{\rm ev}}

\subheading{\Sec8. Fourier coefficients of Siegel Eisenstein series}

In this section, we recall, from \cite{\ked}, the construction of
certain incoherent Siegel Eisenstein series and the structure of
the Fourier coefficients of their derivative at $s=0$, the center
of symmetry. To be more precise, these Eisenstein series occur on
the metaplectic cover of the symplectic group of rank $4$ over
$\Q$, and have an odd functional equation. Their Fourier
coefficients are parameterized by rational symmetric matrices $T\in
Sym_4(\Q)$. In \cite{\ked}, a formula was given for the derivative
at $s=0$ of such a coefficient, when $\det(T)\ne0$.

We retain the notation of section 1, and we refer to sections 1 --
6 of \cite{\ked} for more details. Thus $B$ is an indefinite
quaternion algebra over $\Q$ of discriminant $D(B)$, $C = M_2(B)$,
$V$ is given by (1.1), and $G$ is given by (1.3), etc.. In
particular,  $V$ is a five-dimensional quadratic space over $\Q$
with signature $(3,2)$. Let $\chi=\chi_V$ be the quadratic
character of $\A^\times/\Q^\times$ attached to $V$: $\chi(x) =
(x,\det(V))_\A$, where $(\ ,\ )_\A$ is the global Hilbert symbol.
Note that $\chi_\infty(-1)=1$.

Let $W$ be a symplectic vector space of dimension $8$ over $\Q$,
with a fixed symplectic basis $e_1,\dots,e_4,e'_1,\ldots, e'_4$,
and let $H_\A$ be the metaplectic extension of $Sp(W_\A)$, with
Siegel parabolic $P_\A$. For $s\in \C$ and for $\chi$ as above, let
$I_4(s,\chi)$ be the global degenerate principal series
representation of $H_\A$. As explained in \cite{\ked}, (2.9), the
representation $I_4(0,\chi)$ has a direct sum decomposition into
two types of irreducible representations. One of these types are
the irreducible summands, like $\Pi_4(V)$, associated to
five-dimensional quadratic spaces with character $\chi_V$. The
other type are the irreducible summands associated to {\bf
incoherent collections,} in the sense of section 2 of \cite{\ked}.
One such summand is $\Pi_4(\Cal C)$, associated to the incoherent
collection $\Cal C$, defined as follows. For any finite prime
$\ell$, $\Cal C_\ell=V_\ell$, while $\Cal C_\infty = V'_\infty$,
where $V'_\infty$ is the quadratic space over $\R$ of signature
$(5,0)$. There is a surjective map
$$\lambda_f: S((\Cal C_{\A_f})^4) = S(V(\A_f)^4) \lra \Pi_4(\Cal C)_f\subset
I_4(0,\chi)_f.\tag8.1$$

A section $\P(s)\in I_4(s,\chi)$ is {\bf standard} if its
restriction to the standard maximal compact subgroup $K_H$ in
$H_\A$ is independent of $s$. For $\ph_f\in S(V(\A_f)^4)$, let
$\P_f(s)$ be the standard section of $I_4(s,\chi)_f$ such that
$\P_f(0)=\lambda(\ph_f)$. Let $\P(s) =
\P_\infty^{\frac52}(s)\tt\P_f(s)$, where $\P_\infty^{\frac52}(s)$
is the standard section of $I_4(s,\chi)_\infty$ whose restriction
to $K_{H_\infty}$ is the character $\det^{\frac52}$. Then $\P(s)$
is an incoherent section with $\P(0)\in\Pi_4(\Cal C)$. The
incoherent Eisenstein series
$$E(h,s,\P)=\sum_{\gamma\in P_\Q\back H_\Q} \P(\gamma h,s)\tag8.2$$
converges for $\text{Re}(s)> \frac52$, and its analytic continuation
vanishes at the point $s=0$, \cite{\ked}. There is a Fourier expansion
$$E(h,s,\P)=\sum_{T\in Sym_4(\Q)}  E_T(h,s,\P),\tag8.3$$
with respect to the unipotent radical of $P$.
When $\P(s)=\tt_{\ell} \P_\ell(s)$ is a factorizable
section, and when $\det(T)\ne 0$, there is a product formula
$$E_T(h,s,\P) =\prod_{\ell\le\infty} W_{T,\ell}(h_\ell,s,\P_\ell),\tag8.4$$
where $W_{T,\ell}(h_\ell,s,\P_\ell)$ is the local generalized Whittaker
integral,
cf. section 4 of \cite{\ked}. For fixed $h$, $T$, and $\P$, there is a finite
set of places
$S$ such that, \cite{\ked}, Proposition~4.1,
$$\prod_{\ell\notin S} W_{T,\ell}(h_\ell,s,\P_\ell) =
\zeta^S(2s+4)^{-1}\zeta^S(2s+2)^{-1},\tag8.5$$
and hence
$$E_T(h,s,\P) =\zeta^S(2s+4)^{-1}\zeta^S(2s+2)^{-1}
\cdot \prod_{\ell\in S} W_{T,\ell}(h_\ell,s,\P_\ell).\tag8.6$$
Since $\det(T)\ne 0$, the factors $W_{T,\ell}(h_\ell,s,\P_\ell)$ have an entire
analytic
continuation.

Fix $T$ with ${\roman{det}}(T)\ne 0$. Since $E_T(h,0,\P)=0$, at least one of
the factors in the product formula (8.6) vanishes at $s=0$. In
particular, by Proposition~1.4 of \cite{\ked}, the factor at $\ell$
vanishes whenever the five-dimensional quadratic space $\Cal
C_\ell$ does not represent $T$. Let $\text{Diff}(T,\Cal C)_f$ be
the set of finite places at which $\Cal C_\ell$ fails to represent
$T$, and let
$$
\text{Diff}(T,\Cal C) = \cases \text{Diff}(T,\Cal C)_f \cup \{\infty\}&
\text{ if $sig(T) = (3,1)$ or $(1,3)$}\\
\text{Diff}(T,\Cal C)_f &\text{ otherwise.}
\endcases\tag8.7
$$
By Corollary~5.3 of \cite{\ked},  $|\text{Diff}(T,\Cal C)|$ is odd;
and, by Corollary~5.4 of loc.\ cit.,
$$\ord{s=0}E_T(h,s,\P) \ge |\text{Diff}(T,\Cal C)|.\tag8.8$$
Thus, the only nonsingular $T$ for which $E'_T(h,0,\P)$ can be
nonzero are those for which $|\text{Diff}(T,\Cal C)|=1$. We will
relate the value $E'_T(h, 0, \Phi)$ for ${\roman{Diff}}(T, {\Cal
C})= \{ p\}$ to the numbers $<{\Cal Z}(T,\omega)>_p$ in the previous
section.

\comment
To do something similar for $<{\Cal Z}(d_1,
\omega_1),\ldots, {\Cal Z}(d_r, \omega_r)>_p^{\roman{proper}}$ we will need
more notation.

For $n_1,\ldots, n_r$ with $1\le n_i\le 4$ and with $n_1+\ldots+n_r=4$, let
$$W=W_1+\ldots+W_r\tag8.9$$
be a decomposition of $W$ into symplectic subspaces of
dimensions $2n_i$, compatible with the fixed symplectic basis, and let
$$\iota: H_{1,\A}\times\ldots\times H_{r,\A}\lra H_\A\tag8.10$$
be the corresponding homomorphism of metaplectic groups, covering the
embedding
$$\iotai:Sp(W_{1,\A})\times\ldots\times Sp(W_{r,\A})\hookrightarrow
Sp(W_\A).\tag8.11$$

Let
$$F(h_1,\ldots,h_r,\P) =
\frac{\partial\phantom{s}}{\partial
s}\bigg\{E(\iota(h_1,\ldots,h_r),s,\P)\bigg\}
\bigg|_{s=0}\tag8.12$$
be the restriction to $H_{1,\A}\times\ldots\times H_{r,\A}$ of the
derivative of $E(h,s,\P)$ at $s=0$. By Corollary~6.3 of \cite{\ked},
the Fourier coefficients of this function with respect to the
unipotent radicals of the Siegel parabolics of $H_1,\ldots,H_r$ are
indexed by $d_1\in Sym_{n_1}(\Q),\ldots, d_r\in Sym_{n_r}(\Q)$, and
have the following form:
$$\align
& F_{d_1,\ldots,d_r}(h_1,\ldots,h_r,\P)\\
&
= \sum_{p\le\infty}
F_{d_1,\ldots,d_r}(h_1,\ldots,h_r,\P)_p
+F_{d_1,\ldots,d_r}(h_1,\ldots,h_r,\P)_{\text{sing}},\tag8.13\\
\endalign$$
where
$$F_{d_1,\ldots,d_r}(h_1,\ldots,h_r,\P)_p= \sum\limits_{\matrix
\nass
T\in Sym_4(\Q), \ \det(T)\ne 0\\
\nass
{\roman{diag}}(T)=(d_1,\ldots,d_r)\\
\nass
\text{\rm Diff}(T,\Cal C)=\{p\}\endmatrix}
E'_T(\iota(h_1,\ldots,h_r),0,\P),\tag8.14$$ and
$$F_{d_1,\ldots,d_r}(h_1,\ldots,h_r,\P)_{\text{sing}}=\sum\limits_{
\matrix
\nass
T\in Sym_4(\Q), \ \det(T)= 0\\
\nass
{\roman{diag}}(T)= (d_1,\ldots, d_r)
\endmatrix} E'_T(\iota(h_1,\dots, h_r),0,\P).\tag8.15$$
\endcomment

Let us fix a finite prime $p$. We wish to give a formula for
$E'_T(h,0,\Phi)$ if $T\in {\roman{Sym}}_4({\Bbb Q})$ is
nonsingular with ${\roman{Diff}}(T, {\Cal C})=\{ p\}$. Let
$B^{\stuff}$ be the definite quaternion algebra over $\Q$ which is
ramified at $p$ and whose invariants coincide with those of $B$ at
all finite primes other than $p$. Let $C^{\stuff}
= M_2(B^{\stuff})$, and let
$$V^{\stuff} =\{\ x\in M_2(B^{\stuff})\mid x'=x \text{ and } tr(x)=0\
\},\tag8.9$$
with quadratic form defined by squaring, as in section 1. Let
$G^{\stuff}=GSpin(V^{\stuff})$ be defined by the analogue of (1.3). Note that
there is an
exact sequence
$$1\lra Z'\lra G^{\stuff}\lra SO(V^{\stuff})\lra1\tag8.10$$
of algebraic groups over $\Q$, where $Z'$ is the center of $G^{\stuff}$.

We fix identifications $B^{\stuff}(\A_f^p) = B(\A_f^p)$, and hence
$V^{\stuff}(\A_f^p)=V(\A_f^p)$, and $G^{\stuff}(\A_f^p)=G(\A_f^p)$.  We
also assume that $\ph_f\in S(V(\A_f)^4)$ is factorizable, so that
$\ph_f=\ph_p\tt\ph_f^p$, and we can view $\ph_f^p$ as a Schwartz
function on $V^{\stuff}(\A_f^p)^4$. Recall that there is a surjective
map
$$\lambda^{\stuff}_f:S(V^{\stuff}(\A_f)^4)\lra \Pi_4(V^{\stuff})_f\subset
I_4(0,\chi)_f.\tag8.11$$

Recall, \cite{\sweet}, \cite{\ked}, that the local degenerate
principal series representation $I_{4,p}(0,\chi_p)$ has a direct
sum decomposition with irreducible factors
$$I_{4,p}(0,\chi_p) = R_4(V_p)\oplus R_4(V_p^{\stuff}).\tag8.12$$
Let $T\in{\roman{Sym}}_4({\Bbb Q})$ be nonsingular with $\text{\rm
Diff}(T,\Cal C)=\{p\}$. Then the linear functional
$$W_{T,p}(h,0,\cdot): I_{4,p}(0,\chi_p) \lra \C\tag8.13$$
vanishes identically on $ R_4(V_p)= R_4(\Cal C_p)$, and does not
vanish identically on the summand $R_4(V_p^{\stuff})$, \cite{\ked},
Proposition~1.4. We choose a standard section $\P_p'(s)$, with
$\P'_p(0)\in R_4(V_p^{\stuff})$, and such that for a given $h\in
H_{{\Bbb Q}_p}$,
$$W_{T,p}(h,0,\P'_p)\ne 0.\tag8.14$$
Let $\ph'_p\in S((V^{\stuff}_p)^4)$ be a Schwartz function whose image
$\lambda_p(\ph'_p)$ in $I_{4,p}(0,\chi_p)$ is $\P'_p(0)$. Note that
$V^{\stuff}$ is positive definite, and let $\ph'_\infty\in
S((V^{\stuff}_\infty)^4)$ be the Gaussian, $\ph'_\infty(x) =
\exp{(-\pi tr(q(x)))}$. Finally, let $\ph_f^{\stuff}=\ph'_p\tt\ph_f^p$ so
that
$$\varphi^{\stuff}=\ph'_\infty\tt\ph_f^{\stuff} =
\ph'_\infty\tt\ph'_p\tt\ph_f^p
\in S(V^{\stuff}(\A)^4).\tag8.15$$

Recall that the metaplectic group $H_\A$ acts on the space
 $S(V^{\stuff}(\A)^4)$ via the
Weil representation $\o=\o_\psi$, defined using our fixed additive character
$\psi$ of
$\A/\Q$.
For $g\in G^{\stuff}(\A)$ and $h\in H_\A$, let
$$\theta(g,h,\ph^{\stuff}) = \sum_{\y\in V^{\stuff}(\Q)^4}
\big(\o(h)\ph^{\stuff}\big)(g^{-1}\y)\tag8.16$$
be the theta function attached to $\ph^{\stuff}$, and
let
$$I(h,\ph^{\stuff})= \frac12\int_{G^{\stuff}(\Q)Z(\A)\back G^{\stuff}(\A)}
\theta(g,h,\ev(\ph^{\stuff}))\, dg,\tag8.17$$
for the Tamagawa measure $dg$ on $Z(\A)\back G^{\stuff}(\A)$, and
where $\ev(\ph^{\stuff})$ denotes the projection of $\ph^{\stuff}$
to the subspace of functions all of whose local components are even,
cf. \cite{\ked}, (7.19). Note that $\theta(g,h,\ph^{\stuff})$ can be defined
by the same formula for $g\in O(V^{\stuff})(\A)$, and that
$$I(h,\ph^{\stuff})= \int_{O(V^{\stuff})(\Q)\back O(V^{\stuff})(\A)}
\theta(g,h,\ph^{\stuff})\, dg,\tag8.18$$
where $\text{vol}(O(V^{\stuff})(\Q)\back O(V^{\stuff})(\A),dg)=1$.

For $g\in G'(\A)$ and $h\in H_\infty$, let
$$\Lt(g,h;\ph') = \sum_{\y\in \O_T'(\Q)}
\big(\o(h)\ph'\big)(g^{-1}\y),\tag8.19$$
and
$$\Ltf(g,\ph'_f) = \sum_{\y\in \O_T'(\Q)} \ph'_f(g^{-1}\y),\tag8.20$$
where $\Ltf(g,\ph'_f)$ depends only on $g_f$.
Recall that $\ph'=\ph'_\infty\tt\ph'_f$, with
$$\ph'_\infty(x) = e^{-\pi tr((x,x))},\tag8.21$$
the Gaussian attached to $V'$. This function is invariant under
$G'(\R)$ and, for $\y\in \O_T'(\Q)$, it has the value
$$\ph'_\infty(\y) = e^{-2\pi tr(T)}.\tag8.22$$
Therefore, for $h\in H_\infty$, we have
$$\align
\theta_T(h,g;\ph') &=  \sum_{\y\in \O_T'(\Q)}
\big(\o(h)\ph'\big)(g^{-1}\y)\tag8.23\\
\nass
{}&= \big(\o(h)\ph'_\infty\big)(\y_0)\cdot \Ltf(g,\ph'_f),\endalign
$$
where $\y_0$ is any fixed element of $\O_T'(\Q)$.

For $h\in H_\infty$, and for
$\y_0\in \O_T'(\Q)$, set
$$ W_T^{\frac52}(h):=\big(\o(h)\ph'_\infty\big)(\y_0).\tag8.24$$
More explicitly, as in (11.74)
of \cite{\ked}, if $h$ has Iwasawa decomposition
$h = (n(b) m(a) k, t) \in Sp_4(\R)\times \C^1 \simeq Mp_4(W_\infty)$,
for $b\in Sym_4(\R)$, $a\in GL_4(\R)^+$, and $k\in K_{H_\infty}$, then
$$\align
W^{\frac52}_T(h) &= t\cdot\det(a)^{\frac52}\, e(tr(Tb))\, e^{-\pi
tr({}^taTa)}\,
\det(k)^{\frac52}\tag8.25\\
{}&= t \cdot\det(a)^{\frac52}\, e(tr(T\tau))\, \det(k)^{\frac52},
\endalign$$
where $\tau = b + i a {}^ta$.

Recalling that $Z'(\R)\back G'(\R)\simeq SO(V')(\R)$ is compact, we have the
following formula for the $T$-th Fourier coefficient of the theta integral:
$$
\align
{}&2\,I_T(h,\ph')\tag8.26\\
{}&= \int_{G^{\stuff}(\Q)Z'(\A)\back G^{\stuff}(\A)}
\theta_T(g,h;\ev(\ph'))\, dg\\
\nass
{}&=W^{\frac52}_T(h)\cdot
\int_{G^{\stuff}(\Q)Z'(\A)\back G^{\stuff}(\A)} \Ltf(g,\ev(\ph'_f))\, dg\\
\nass
{}&=
W^{\frac52}_T(h)\cdot \vol(SO(V')(\R),d_\infty g)\cdot
\int_{G^{\stuff}(\Q)Z'(\A_f)\back G^{\stuff}(\A_f)} \Ltf(g,\ev(\ph'_f))\,
d_fg,\endalign
$$
where $d_fg$ is the measure arising from the counting measure on
$Z'(\Q)\back G'(\Q)$ and the Haar measure on $Z'(\A_f)\back
G^{\stuff}(\A_f)\simeq SO(V')(\A_f)$ coming from some choice of a
gauge form $\mu$ on $SO(V')$. Also $d_\infty g$ is the Haar measure
on $SO(V')(\R)$ induced by $\mu$.

With the notation just described, and for $h\in H_\infty$,
Corollary~6.3 of \cite{\ked} specializes to
$$
E'_T(h,0, \Phi)
= \frac{W'_{T,p}(e,0,\Phi_p)}{W_{T,p} (e,0,\Phi'_p)}
\cdot 2 I_T(h,\varphi^{\stuff}),\tag8.27
$$
if $T$ is nonsingular with ${\roman{Diff}}(T,{\Cal C})=\{ p\}$.

Substituting the expression (8.26) for the Fourier coefficient of the theta
integral found
above, we obtain:
\proclaim{Proposition 8.1} Suppose that
$\P(s)=\P_\infty^{\frac52}(s)\tt\P_f(s)$ with $\P_f(0) =
\lambda_f(\ph_f)$, is an incoherent standard section. For $h\in
H_{\infty}$, and for each $T\in Sym_4(\Q)$ with $\det(T)\ne0$ and
$\text{Diff}(T,\Cal C)=\{p\}$, choose $\ph'_p$ and $\P'_p(s)$, such
that $W_{T,p}(e,0,\Phi'_p)\ne 0$. Then
$$E'_T(h,0,\Phi)= \vol(SO(V')(\R))\cdot W_T^{5/2}(h)\cdot
\frac{W'_{T,p}(e,0,\Phi_p)}{W_{T,p}
(e,0,\Phi'_p)} \cdot I_{T,f}(\varphi_f^{\stuff}).$$
Here
$$\align
I_{T,f}(\ph_f^{\stuff})&=\int_{G^{\stuff}(\Q)Z'(\A_f)\back G^{\stuff}(\A_f)}
\sum_{\y\in \O_T^{\stuff}(\Q)}
\ev(\ph_f^{\stuff})(g^{-1}\y)\ d_fg\\
\nass
{}&=\int_{G^{\stuff}(\Q)Z'(\A_f)\back G^{\stuff}(\A_f)}
\Ltf(g,\ev(\ph'_f))\ d_fg,\endalign
$$
and the measures are as described after (8.26) above.
\endproclaim

If the function $\ph'_f$ is locally even, then the integral
$$I_{T,f}(\ph_f') = \int_{G'(\Q)Z'(\A_f)\back G'(\A_f)} \Ltf(g,\ph'_f)\
d_fg\tag8.28$$ occurs in Theorem~7.2, where the measure arises from
an arbitrary Haar measure on $Z'(\A_f)\back G'(\A_f)$, and the
quantity
$$\vol(\pr(K'))^{-1} I_{T,f}(\ph_f')\tag8.29$$
is independent of the choice. Therefore, we can obtain the expression
$$\multline
E'_T(h,0,\Phi) = \vol(SO(V')(\R)\pr(K'))
\cdot W_T^{5/2}(h)\cdot\\
\cdot\frac{W'_{T,p}(e,0,\Phi_p)}{W_{T,p}
(e,0,\Phi'_p)} \cdot
\vol(\pr(K'))^{-1}I_{T,f}(\varphi_f'),\endmultline\tag8.30$$
where the factor $\vol(SO(V')(\R)\pr(K'))$ is computed using the Tamagawa
measure on $SO(V')(\A)$. Hence, since $\pr(K')$ is neat,
$$\vol(SO(V')(\R)\pr(K')) = 2 |SO(V')(\A): SO(V')(\Q)SO(V')(\R)\pr(K')|^{-1},$$
and the quantities in (8.30) separated by a dot do not depend on
any choice of measure.

\vfill\eject
\subheading{\Sec9. The main theorem}

In this section we assemble the results of previous sections and state our main
results.

We begin by further specializing the formula of Proposition~8.1. Specifically,
we need more information about the factor
$$\frac{W'_{T,p}(e,0,\P_p)}{W_{T,p}(e,0,\P_p')}.\tag9.1$$
Fix the prime $p$ with $p\nmid 2D(B)$, and assume that $\ph_p$ is the
characteristic function of $V(\Z_p)^4$. Recall that $\P_p(s)$ is
the standard section with $\P_p(0) =\lambda_p(\ph_p)$. Also, let
$\ph'_p$ be the characteristic function of the lattice
$V'(\Z_p)^4$, and let $\P'_p(s)$ be the standard
section with $\P'_p(0) = \lambda'_p(\ph'_p)$.

Recall that a nonsingular $T\in \text{\rm Sym}_4(\Q_p)$ is represented by
precisely one of the quadratic spaces $V(\Q_p)$ and $V'(\Q_p)$, \cite{\ked},
Proposition~1.3.
\proclaim{Proposition 9.1} Suppose
that $\varphi_p, \varphi'_p,
\Phi_p, \Phi'_p$ are as above, and that $T\in \text{\rm Sym}_4(\Q_p)$ with
$\det(T)\ne 0$.\hfill\break
(i) If $W'_{T,p}(e,0,\Phi_p)\ne 0$, then $T\in{\roman{Sym}}_4({\Bbb Z}_p)$.
\hfill\break
(ii) If $T\in {\roman{Sym}}_4({\Bbb Z}_p)$
and if $T$ is represented by $V'({\Bbb Q}_p)$, then
$W_{T,p}(e,0,\Phi'_p)\ne 0\ \ $.\hfill\break
(iii) If $T\in {\roman{Sym}}_4({\Bbb Z}_p)$ is represented by $V'(\Q_p)$, and
if
$T$ represents $1$, then
$${W'_{T,p}(e, 0,\Phi_p)\over W_{T,p}(e, 0, \Phi'_p)} ={1\over 2}
{\roman{log}}\,
p\cdot (p^2+1)(p-1)\cdot e_p(T),$$
where $e_p(T)$ is the local intersection multiplicity given in Proposition~6.1.
\endproclaim
The proof will be given in section~10.

A subset $\o\subset V(\A_f^p)^n$ is said to be {\bf locally
centrally symmetric} if it is invariant under the action of the
group $\mu_2(\A_f^p)$. The characteristic function $\ph_\o\in
S(V(\A_f^p)^n)$ of such a set is locally even, as in (8.17),\ i.e.
$\ph_\o={\roman{ev}}(\ph_\o)$. The function
$\ph'_f=\ph'_p\tt\ph_\o\in S(V'(\A_f)^n)$ is then locally even as
well, so that the expression (8.30) holds for the derivative of the
Fourier coefficients of the associated Eisenstein series.

Our first main result is the following.
\proclaim{Theorem 9.2}
Assume that $p\nmid 2 D(B)$ and that $\varphi_p, \varphi'_p,
\Phi_p, \Phi'_p$ are as above. Let $\omega \subset V({\Bbb A}_f^p)^4$ be a
locally centrally symmetric $K^p$-invariant compact open subset.
Let $\Phi(s)=\Phi_\infty(s)\otimes \Phi_p(s)\otimes
\Phi_f^p(s)$ be the standard section corresponding to
$\ph=\varphi_\infty\otimes \varphi_p\otimes
\varphi_f^p\in S(V^{(p)}(\A)^4)$ with $\varphi_f^p={\roman{char}} (\o)$, cf
Lemma~7.1.
Suppose that $T\in{\roman{Sym}}_4({\Bbb Q})$ with
${\roman{det}}(T)\ne 0$ and with $\text{\rm Diff}(T,\Cal C)=\{p\}$.
\hfill\break
(i) If $T\not\in {\roman{Sym}}_4({\Bbb Z}_{(p)})_{>0}$,
then ${\Cal Z}(T,\omega)=\emptyset$, $<{\Cal Z}(T, \omega)>_p=0$,  and
$$E'_T(h,0,\Phi)=0.$$
(ii) If $T\in {\roman{Sym}}_4({\Bbb Z}_{(p)})_{>0}$ represents 1 over ${\Bbb
 Z}_p$,
then ${\Cal Z}(T,\omega)$ is zero dimensional, and, for $h\in H_{\infty}$,
$$E'_T(h,0,\Phi) = \frac12\vol(SO(V')(\R))\cdot
W_T^{5/2}(h)\cdot{\roman{vol}}(\pr(K))
\cdot {\roman{log}}\, p <{\Cal Z}(T, \omega)>_p\ \ .$$
\endproclaim

Note that, if $T\in {\roman{Sym}}_4({\Bbb Z}_{(p)})_{>0}$ does not
represent $1$, then $\Cal Z(T,\o)$ contains components of the
supersingular locus (Corollary~5.15 and Theorems~5.12 and~5.14). In
this case, we do not have a formula for $<{\Cal Z}(T, \omega)>_p$.

In Theorem~9.2, the chosen gauge form $\mu$ on $SO(V') = Z'\back G'$ determines
the Haar measure on $SO(V')(\R)$ used to compute $\vol(SO(V')(\R))$. The
corresponding
gauge form on the inner twist $SO(V) = Z\back G$ determines the measure on
$Z'(\A_f)\back G'(\A_f)$ used to compute $\vol(\pr(K))$. Note that
the product  $\vol(SO(V')(\R)\vol(\pr(K))$ is independent of the choice of
$\mu$.

\demo{Proof of Theorem~9.2} Beginning with formula (8.30),
and using (iii) of Proposition~9.1 and Theorem~7.2, we have
$$\align
E'_T(h,0,\P) &= \vol(SO(V')(\R)\pr(K'))\cdot W_T^{5/2}(h)\cdot\\
\nass
{}&\qquad\cdot \frac{W'_{T,p}(e,0,\Phi_p)}{W_{T,p} (e,0,\Phi'_p)}
\cdot \vol(\pr(K'))^{-1}I_{T,f}(\varphi_f')\tag9.2\\
\nass
{}&= \vol(SO(V')(\R)\pr(K'))\cdot W_T^{5/2}(h)\cdot\\
\nass
{}&\qquad\cdot {1\over 2} {\roman{log}}\, p\cdot (p^2+1)(p-1)\cdot
e_p(T) \cdot \vol(\pr(K'))^{-1}I_{T,f}(\varphi_f')\\
\nass
{}&= \vol(SO(V')(\R)\pr(K'))\cdot W_T^{5/2}(h)\cdot\\
\nass
{}&\qquad\cdot {1\over 2} {\roman{log}}\, p\cdot (p^2+1)(p-1)\cdot
<{\Cal Z}(T, \omega)>_p.\endalign
$$
To finish the proof, we simply note the following relation between
volumes.
\proclaim{Lemma 9.3} Recall that $K_p = GL_2(\Cal O_{B_p})\cap G(\Q_p)$
and $K'_p = GL_2(\Cal O_{B'_p})\cap G'(\Q_p)$. Then, for the
Haar measures on $Z'(\A_f)\back G'(\A_f)$, $Z'(\Q_p)\back G'(\Q_p)$,
$Z(\A_f)\back G(\A_f)$, and $Z(\Q_p)\back G(\Q_p)$
determined by the fixed gauge form $\mu$ and the corresponding form on the
inner twist, \hfill\break
$$\frac{\vol(\pr(K))}{\vol(\pr(K'))} = \frac{\vol(\pr(K_p))}{\vol(\pr(K'_p))}
= (p^2+1)(p-1).$$
\endproclaim
\noindent This finishes the proof of Theorem~9.2
\qed\enddemo

\demo{Proof of Lemma~9.3, following Kottwitz \cite{\kottwitztam}}
We may replace $G/Z$ and $G'/Z'$ by their simply connected
coverings $\tilde G$ resp.\ $\tilde G'$ and ${\roman{pr}}(K_p)$ and
${\roman{pr}}(K'_p)$ by their inverse images $\tilde K_p$ resp.\
$\tilde K'_p$. We use on $\tilde G(\Q_p)$ resp.\ $\tilde G'(\Q_p)$
the Haar measure induced by a top differential form on the
$\Z_p$-form of $\tilde G$ resp.\ $\tilde G'$ corresponding to an
Iwahori subgroup $\tilde I_p\subset \tilde K_p$ resp.\ $\tilde
I'_p\subset \tilde K'_p$. These measures are compatible, cf.\
\cite{\kottwitztam}, p.\ 632. The volumes of $\tilde I_p$ and
$\tilde I'_p$ are related as follows. Choose as in
\cite{\kottwitztam} a maximal split torus $S$ in $\tilde G$ and a
maximal torus $S_1$ containing $S$ which splits over an unramified
extension. We also denote by $S_1$ the canonical $\Z_p$-form of
$S_1$. Choose $S', S'_1$ of the same sort for $\tilde G'$. Then
$${\vol(\tilde I_p)\over \vol(\tilde I'_p)} = {S_1(\F_p)\over S'_1(\F_p)}
= {(p-1)^2\over p^2-1}\ \ ,$$
since in the case at hand $S_1\cong {\Bbb G}_m^2$ and $S'_1\cong
{\roman{Res}}_{\Q_{p^2}/\Q_p}{\Bbb G}_m$. The result follows since
$$\vert\tilde K_p/\tilde I_p\vert = 1+2p+2p^2+2p^3+p^4\ \ ,
\ \ \vert \tilde K'_p/\tilde I'_p\vert = p+1\ \ ,$$
hence
$${\vol({\roman{pr}}(K_p))\over
\vol({\roman{pr}}(K'_p))} = {\vol(\tilde I_p)\over \vol(\tilde I'_p)}
\cdot {\vert \tilde K_p/\tilde I_p\vert\over \vert \tilde
K'_p/\tilde I'_p\vert} = {(p-1)^2\over p^2-1} \cdot {p^4-1\over
p-1}\ \ .\qquad\qed$$
\enddemo

We next formulate the corresponding result for the intersection of
special cycles.

For $n_1,\ldots, n_r$ with $1\le n_i\le 4$ and with
$n_1+\ldots+n_r=4$, let $d_i\in {\roman{Sym}}_{n_i}({\Bbb
Z}_{(p)})_{>0}$ and fix locally centrally symmetric $K^p$-invariant
open compact subsets $\omega_i\subset V(\A_f^p)^{n_i}$. Let
$$W=W_1+\ldots+W_r\tag9.3$$
be a decomposition of $W$ into symplectic subspaces of
dimensions $2n_i$, compatible with the fixed symplectic basis, and let
$$\iota: H_{1,\A}\times\ldots\times H_{r,\A}\lra H_\A\tag9.4$$
be the corresponding homomorphism of metaplectic groups, covering the
embedding
$$\iota:Sp(W_{1,\A})\times\ldots\times Sp(W_{r,\A})\hookrightarrow
Sp(W_\A).\tag9.5$$
Restricting to the archimedean place, for $(h_1,\dots,h_n)\in
H_{1,\infty}\times\ldots\times H_{r,\infty}$, we have
$$W_T^{\frac52}(\iota(h_1,\dots,h_r)) = W_{d_1}^{\frac52}(h_1)\dots
W_{d_r}^{\frac52}(h_r),
\tag9.6$$
where $T$ has diagonal blocks $d_1,\dots,d_r$.
Thus, by (7.3), we obtain:
\proclaim{Corollary 9.4} With the above notations,
$$\multline
\sum_T E'(\iota(h_1,\dots,h_r),0,\P) =
\frac12\vol(SO(V')(\R))\cdot
W_{d_1}^{\frac52}(h_1)\dots W_{d_r}^{\frac52}(h_r)\\
\times{\roman{vol}}(\pr(K))
\cdot {\roman{log}}\, p
<{\Cal Z}(d_1, \omega_1),\dots,{\Cal Z}(d_r, \omega_r)>^{\text{proper}}_p,
\endmultline$$
where the intersection number on the right side is defined by
(7.3), and the summation runs over $T\in
\text{Sym}_4(\Z_{(p)})_{>0}$ such that $\text{Diff}(T,\Cal
C)=\{p\}$, $\text{diag}(T)=(d_1,\dots,d_r)$, and $T$ represents $1$
over $\Z_p$. Also, $\P$ is determined as in Theorem~9.2 with
$\o=\o_1\times
\dots\times\o_r$.
\endproclaim

Of course, the left side of the expression of Corollary~9.4 is part
of the $(d_1,\dots,d_r)$-th Fourier coefficient of the pullback
$$F(h_1,\dots,h_r;\P):= E'(\iota(h_1,\dots,h_r),0,\P),\tag9.7$$
cf. \cite{\ked}, (6.13).
This result gives an analogue of the results of \cite{\ked}.

\vfill\eject
\subheading{\Sec10. Representation densities}

In this section, we give the proof of Propositions~9.1,
which is based on a formula of Kitaoka, \cite{\kitaokatwo}, for representation
densities. We then describe a conjectural generalization of Kitaoka's
formula. In this section, for $x\in \Q_p^\times$, $\chi(x) = (x,p)_p$.

We begin by recalling the well known relation between the values of
the function $W_{T,p}(e,s,\P_p)$, at integer values of $s$ and
classical representation densities.

For a suitable choice of basis for $V(\Z_p)$ the quadratic form $q$
has matrix
$$S=S_0= \pmatrix 1&{}&{}\\
{}&{}&\frac12\cdot1_2\\ {}&\frac12\cdot1_2&{}\endpmatrix.\tag10.1$$
For $r\ge 0$, let
$$S_r  =  \pmatrix S_0&{}&{}\\
{}&{}&\frac12\cdot1_r\\ {}&\frac12\cdot1_r&{}\endpmatrix.\tag10.2$$
For nonsingular matrix $T\in Sym_4(\Z_p)$, let
$$\a_p(S_r,T)=
\lim_{t\rightarrow\infty}
p^{-t(10+8r)}\#\{\ x\in M_{5+2r,4}(\Z/p^t\Z)\ \mid \ S_r[x]-T\in
p^t Sym_4(\Z_p)\ \}\tag10.3$$
be the classical representation
density \cite{\kitaokabook}, p.98. This quantity depends only on the
$GL_4(\Z_p)$-equivalence class of $T$, so we assume that
$$T=\text{diag}(\e_0 p^{a_0},\e_1 p^{a_1},\e_2 p^{a_2}, \e_3
p^{a_3}),\tag10.4$$
with $\e_i\in \Z_p^\times$ and $0\le a_0\le a_1\le a_2\le a_3$.
Then, as explained in Corollary~A.1.5 of \cite{\ked},
$W_{T,p}(e,r,\Phi_p)=0$ if
$T\in {\roman{Sym}}_4({\Bbb Q}_p)\setminus {\roman{Sym}}_4({\Bbb Z}_p)$,
and
$$W_{T,p}(e,r,\P_p) = \a_p(S_r,T)\tag10.5$$
if $T\in{\roman{Sym}}_4({\Bbb Z}_p)$, since the factor
$\gamma_p(V_p)$ in loc.cit. is $1$ in our present case. Recall --
see \cite{\kitaokatwo}, Lemma~9 and the discussion on pp. 450--453,
for example --  that $\a_p(S_r,T)$ is a rational function of
$X=p^{-r}$, i.e.\ there is a rational function $A_{S,T}(X)$ such
that
$$\alpha_p(S_r, T)= A_{S, T}(p^{-r})\ \ .\tag10.6$$
We therefore have
$$W'_{T,p}(e,0,\P_p) = -\log(p)\cdot
\frac{\partial}{\partial X}\{A_{S,T}(X)\}\big\vert_{X=1}.\tag10.7$$

At this point we have proved part (i) of Proposition 9.1.

Similarly, let $\ph'_p$ be the characteristic function of the
lattice $V'(\Z_p)^4=V^{(p)}(\Z_p)^4$ and let $\P'_p(s)$ be the
corresponding standard section. Again, for  a suitable choice of
basis for $V'(\Z_p)$, the quadratic form on $V'(\Z_p)$ has matrix
$$S'=S_0'=\text{diag}(1,1,\b,p,-p\b),\tag10.8$$
where $\beta\in {\Bbb Z}^{\times}_p\setminus {\Bbb Z}^{\times,
2}_p$.
Again, the factor $\gamma_p(V'_p)=1$, and so
$$W_{T,p}(e,0,\P'_p) = p^{-4}\cdot \a_p(S_0',T).\tag10.9$$

The following two results imply parts (ii) and (iii) of Proposition~9.1.
\proclaim{Proposition 10.1} Suppose that $T\in\text{\rm Sym}_4(\Z_p)$ is not
represented
by $V(\Q_p)$ and that $T$ represents $1$.
Let $e_p(T)$ be the local intersection
multiplicity, given by the formulas of Proposition~6.1. Then,
$$\align
W'_{T,p}(e,0,\P_p) &= -\log(p)\cdot
\frac{\partial}{\partial X}\{A_{S,T}(X)\}\bigg\vert_{X=1}\\
\nass
{}&=\log\, p\cdot (1-p^{-4})(1-p^{-2}) \cdot e_p(T).\endalign$$
\endproclaim

\proclaim{Proposition 10.2} Suppose that $T\in\text{\rm Sym}_4(\Z_p)$
with $\det(T)\ne0$ represents $1$.
Then
$$W_{T,p}(e,0,\P'_p) = p^{-4}\cdot\a_p(S_0',T) =
\cases p^{-4}(1-p^{-2})2(p+1)&\text{if $V'(\Q_p)$ represents $T$}\\
\nass
0&\text{ otherwise.}
\endcases$$
\endproclaim

Of course, we would like to have analogous information about
$W'_{T,p}(e,0,\P_p)$ and $W_{T,p}(e,0,\P'_p)$ for all $T$. At
first, we simply restrict to the case where $p\nmid T$, so that, we
may assume that $a_0=0$, i.e.,
$$T=\text{diag}(\e_0,\e_1 p^{a_1},\e_2 p^{a_2},\e_3 p^{a_3}).\tag10.10$$
Note that $S\simeq 1_5$.
Then, by the standard reduction formula, \cite{\kitaokanote}, p.149,
%$$\alpha_p(S_1\bot U,T_1\bot U)= \alpha_p(S_1\bot U,U)\cdot\alpha_p(S_1,
%T_1)$$
%if $U$ is a unimodular symmetric matrix, cf. also \cite{\kitaokabook}
%Corollary~5.6.1.)
$$\a_p(S_r,T) = \a_p(S_r,\e_0)\a_p(\tilde S_{r},\tilde T),\tag10.11$$
where $\tilde S_{r}$ is obtained by adding a split space of
dimension $2r$ to
$$\tilde S=\text{diag}(1,1,1,\e_0)\tag10.12$$
and
$$\tilde T = \text{diag}(\e_1 p^{a_1},\e_2 p^{a_2},\e_3 p^{a_3}).\tag10.13$$
Note that
$$\a_p(S_r,\e_0) = (1+ \chi(\e_0)p^{-2-r}) = (1+\chi(\e_0)p^{-2}X),\tag10.14$$
where $X=p^{-r}$, \cite{\tonghai}.

Now suppose that $\chi(\e_0)=1$, i.e., that $T$ represents $1$.
Let $H_{2m}$ be the split quadratic form of
rank $2m$ over $\Z_p$, so that
$$H_{2m}=\pmatrix {}&1_m\\1_m&{}\endpmatrix.\tag10.15$$
Then $\tilde S_r$ is isomorphic to the split space $H_{2r+4}$, and
Kitaoka gives an explicit formula for
the representation density $\a_p(H_{2m},\tilde T)$ for any ternary form
$\tilde T$, \cite{\kitaokatwo}. His formulas, in the cases $a_1-a_2$ even and
$a_1-a_2$ odd, are given as a sum of five double sums!
These can be simplified to yield the following expressions:
\proclaim{Proposition 10.3}{\rm (Kitaoka, \cite{\kitaokatwo})}
Let $X= p^{-r}$, and let
$$\tilde T = \text{diag}(\e_1 p^{a_1},\e_2 p^{a_2},\e_3 p^{a_3}),$$
with $0\le a_1\le a_2\le a_3$.

Let
$$\chi(\tilde T) = \cases 1&\text{ if $a_1\equiv a_2\equiv a_3 \mod(2)$,}\\
\chi(-\e_1\e_2)&\text{ if $a_1\equiv a_2\not\equiv a_3 \mod(2)$,}\\
\chi(-\e_1\e_3)&\text{ if $a_1\not\equiv a_2\not\equiv a_3 \mod(2)$,}\\
\chi(-\e_2\e_3)&\text{ if $a_1\not\equiv a_2\equiv a_3 \mod(2)$.}
\endcases$$
(i) If $a_1\equiv a_2\mod(2)$, then
$$\align
\nass
\frac{\a_p(H_{2r+4},\tilde T)}{(1-p^{-2}X)(1-p^{-2}X^2)} &={}\\
\nass
\sum\limits_{\ell=0}^{\frac{a_1+a_2}2-1}&
p^\ell\bigg( \sum\limits_{k=0}^{\text{min}(a_1,\ell)}
X^{2\ell-k} +\chi(\tilde T) X^{a_1+a_2+a_3 +k-2\ell}\bigg)\\
&{}+p^{\frac{a_1+a_2}2}X^{a_2}
\bigg(\sum\limits_{k=0}^{a_1} X^k\bigg) \bigg(\sum\limits_{j=0}^{a_3-a_2}
(\e X)^j\bigg),\\
\noalign{\noindent
where $\e=\chi(-\e_1\e_2)$.\hfill\break
(ii) If $a_1\not\equiv a_2\mod(2)$, then}
\nass
\nass
\frac{\a_p(H_{2r+4},\tilde T)}{(1-p^{-2}X)(1-p^{-2}X^2)}&{}=\\
\nass
\sum\limits_{\ell=0}^{\frac{a_1+a_2-1}2}&
p^\ell\bigg( \sum\limits_{k=0}^{\text{min}(a_1,\ell)}
X^{2\ell-k} +\chi(\tilde T) X^{a_1+a_2+a_3 +k-2\ell}\bigg).\endalign$$
\endproclaim

Note that these expressions exhibit the functional equation of the
local degenerate Whittaker function under $X\mapsto X^{-1}$.
Evaluating at $X=1$ and taking (10.11) and (10.14) into account, we
obtain:
\proclaim{Corollary 10.4} Suppose that $T$ represents $1$.\hfill\break
(i) If $a_1\equiv a_2\mod(2)$, then
$$\align
\frac{\a_p(S,T)}{(1-p^{-2})(1-p^{-4})}&=
(1 +\chi(\tilde T))\sum\limits_{\ell=0}^{\frac{a_1+a_2}2-1}
\big(\text{min}(a_1,\ell)+1\big)\,p^\ell\\
{}&\qquad +p^{\frac{a_1+a_2}2}
\big(a_1+1\big) \bigg(\sum\limits_{j=0}^{a_3-a_2}
\e^j\bigg),\\
\noalign{\noindent
where $\e=\chi(-\e_1\e_2)$.\hfill\break
(ii) If $a_1\not\equiv a_2\mod(2)$, then}
\frac{\a_p(S,T)}{(1-p^{-2})(1-p^{-4})}&=
(1 +\chi(\tilde T)) \sum\limits_{\ell=0}^{\frac{a_1+a_2-1}2}
\big(\text{min}(a_1,\ell)+1\big)p^\ell.\endalign$$
In case (ii), this quantity vanishes if and only if $\chi(\tilde T)=-1$.
In case (i), if $a_2\equiv a_3\mod(2)$, then $\chi(\tilde T)=1$ and there are
an odd number of terms in the last sum, so that the whole expression is
nonzero.
If $a_2\not\equiv a_3\mod(2)$, then $\chi(\tilde T)=\chi(-\e_1\e_2)=\e$, so
that
the whole expression vanishes if and only if $\chi(\tilde T)=-1$.
\endproclaim

\proclaim{Proposition 10.5} Suppose that $T$ represents $1$.
Also suppose that $\chi(\tilde T)=-1$, so that $T$ is not represented by $S$,
i.e., by $V(\Q_p)$ \hfill\break
(i) If $a_1\equiv a_2\mod(2)$, then
$$\align
\frac{\partial}{\partial X}\bigg\{&
\frac{A_{S,T}(X)}{(1-p^{-2}X^2)(1-p^{-4}X^2)}\bigg\}\bigg\vert_{X=1}\\
\nass
{}&\qquad\qquad= -\sum\limits_{\ell=0}^{\frac{a_1+a_2}2-1}
p^\ell\bigg( \sum\limits_{k=0}^{\text{min}(a_1,\ell)}
(a_1+a_2+a_3 +2k-4\ell)\bigg)\\
&\qquad\qquad\qquad\quad -p^{\frac{a_1+a_2}2}
\big(a_1+1\big) \big(\frac{a_3-a_2+1}2\big).\\
\nass
\noalign{\noindent (ii) If $a_1\not\equiv a_2\mod(2)$, then}
\nass
\frac{\partial}{\partial X}\bigg\{&
\frac{A_{S,T}(X)}{(1-p^{-2}X)(1-p^{-4}X^2)}\bigg\}\bigg\vert_{X=1}\\
\nass
{}&\qquad\qquad=-\sum\limits_{\ell=0}^{\frac{a_1+a_2-1}2}
p^\ell\sum\limits_{k=0}^{\text{min}(a_1,\ell)}
(a_1+a_2+a_3 +2k-4\ell).\endalign$$
\endproclaim

After a short manipulation, these expressions coincide, up to sign, with those
given in Proposition~6.1 for the local intersection multiplicity $e_p(T)$!
\proclaim{Corollary 10.6} Suppose that $p\nmid T$ and $\e_0$ is a square,
i.e., that $T$ represents $1$ over $\Z_p$. Also
suppose that $T$ is not represented by $S$. Then
$$\frac{\partial}{\partial X}\big\{
A_{S,T}(X)\big\}\bigg\vert_{X=1}
=-(1-p^{-2})(1-p^{-4})\,e_p(T),$$
where $e_p(T)$ is as in Proposition~6.1.
\endproclaim

This completes the proof of Proposition~10.1.

\demo{Proof of Proposition~10.2} We apply the reduction formula
to obtain:
$$\a_p(S'_r,T) = \a_p(S'_r,\e_0)\a_p(\tilde S'_{r},\tilde T),\tag10.15$$
where $\tilde S_{r}$ is obtained by adding a split space of
dimension $2r$ to
$$\tilde S'=\text{diag}(1,\e_0\beta,p,-p\beta)\tag10.16$$
and $\tilde T$ is as in (10.12).

If $\e_0$ is a square, then
$$\a_p(\tilde S',\e_0) = (1-\chi(-1)p^{-1}),$$
\cite{\tonghai}. On the other hand,
$\tilde S'$ is just the norm form on the maximal order of the
division quaternion algebra over $\Q_p$.
\proclaim{Lemma 10.7}
$$\a_p(\tilde S',\tilde T) = 2 (1+\chi(-1)p^{-1})(p+1).$$
\endproclaim
\demo{Proof} Let $\Bbb B$ be the division quaternion algebra over
$\Q_p$, and let $R$ be its maximal order. Then, for a suitable
$\Z_p$-basis, $\tilde S'$ is the matrix for the quadratic form $Q$
given by the reduced norm on $R$. Let
$$A_{p^r}(T) = \#\{ x\in (R/p^r R)^3\mid Q[x] \equiv \tilde T\mod p^r\},$$
so that
$$\a_p(\tilde S',\tilde T) = \lim_{r\to\infty} p^{-6r}A_{p^r}(T).$$
Choose a uniformizer $\pi\in R$ such that $\pi^2=-p$, and hence
$Q[\pi x] = p Q[x]$. Note that $x\in R$ if and only if $Q[x]\in \Z_p$.
Thus there is a bijection
$$\multline
\{ x\in (R/p^r R)^3\mid Q[x] \equiv p\tilde T\mod p^r\}
\isoarrow\\
\{ y\in (R/p^{r-1}\pi R)^3\mid Q[y]\equiv\tilde T\mod p^{r-1}\},
\endmultline
$$
given by $x\mapsto \pi^{-1}x$. Since $|R/\pi R|= p^2$, we have
$$A_{p^r}(p\tilde T) = p^{6} A_{p^{r-1}}(\tilde T),$$
and hence
$$\a_p(\tilde S',p\tilde T) = \a_p(\tilde S',\tilde T).$$
Thus, we may replace $\tilde T$ by $T'=\diag(\e_1,\e_2 p^{a_2-a_1},\e_3
p^{a_3-a_1})$.
Here $\e_1$ can be taken to be equal to either $1$ or $\beta$.
Using reduction, we have
$$\a_p(\tilde S',T') = \a_p(\tilde S',\e_1)\a_p(S'',T''),$$
where
$$S'' = \diag(\e_1\beta,p,-\beta p),\qquad\text{ and }\qquad
T''=\diag(\e_2 p^{a_2-a_1},\e_3 p^{a_3-a_1}).$$
By Theorem~3.1  of \cite{\tonghai},
$$\a_p(\tilde S',\e_1) = (1+\chi(-1)p^{-1}).$$
If $\e_1=1$, the form $S''$ is just the norm form on the
trace zero elements in $R$, while, if $\e_1=\beta$, then $S''$
is isomorphic to $\beta$ times this norm form. Since $\a_p(\beta S'',\beta T'')
=\a_p(S'',T'')$, Proposition~8.6 of \cite{\ked} yields
$$\a_p(S'',T'') = \cases 2(p+1)&\text{ if $T''$ is anisotropic,}\\
0&\text{ otherwise.}\endcases$$

\qed\enddemo
\medskip\noindent
{\bf Remark 10.8.} The value $\alpha_p(\tilde S', \tilde T)$ is
given erroneously by Gross and Keating as $2\cdot (1+p^{-1})\cdot
(p+1)$, \cite{\grosskeating}, Proposition 6.10.

Thus
$$\align
\a_p(S',T) &= \a_p(S',1)\a_p(\tilde S',\tilde T)\\
\nass
{}&= 2(1- p^{-2}) (p+1),\endalign
$$
as claimed in Proposition~10.2.\qed\qed\enddemo

\bigskip
\centerline{\bf Notes on Clifford algebras}
\noindent
{\bf A.1.} Let $(V,q)$ be a non-degenerate quadratic space of
dimension 5 over a field $F$ of characteristic not 2. Let $C(V)$ be
its Clifford algebra, with its 2-grading
$$C(V)= C^+(V)\oplus C^-(V)\ \ .$$
The {\bf Clifford involution} $c\mapsto c'$ of $C(V)$ is the unique
involution which acts by the identity map on $V\subset C^-(V)$.
Thus
$$(v_1 \cdots v_r)'= v'_r\cdots v'_1\ \ .$$
If $v_1,\ldots, v_5$ is a basis for $V$, then the element $\delta
=v_1\cdots v_5$ lies in the center of $C(V)$ and satisfies
$$\delta'=\delta\ \ .$$
Let
$$G=G\, {\roman{Spin}} (V)= \{ g\in C^+(V)^{\times}\mid gVg^{-1}
=V,\ \ \hbox{and}\ gg'=\nu(g)\}$$
which may be considered as an algebraic group over
${\roman{Spec}}\, F$.
\par\noindent
{\bf A.2.} In this section suppose that $F$ is algebraically closed
and choose a Witt decomposition of the quadratic space $V$,
$$V=V_+\oplus V_0\oplus V_-$$
where ${\roman{dim}}\, V_{\pm}=2$ and $V_{\pm}$ are maximal
isotropic subspaces of $V$. Let $v_0\in V_0$ be a basis vector with
$q(v_0)=1$. We recall the Spin representation of $G$. We use the
identifications of representations of $C(V)$,
$$\align
C(V)/C(V)C(V_-)_{>0} &
=C(V_+\oplus V_0)=\\
&
=C(V_+)(1+v_0)\oplus C(V_+)(1-v_0)\ \ .
\endalign$$
As $C(V)^+$-modules the last two modules are isomorphic. Either one
of them defines the Spin representation $W$ of $G$. Its dimension
is 4.

Fix an isomorphism $\Lambda^2 V_+=F$ and let
$$\lambda: W\to F$$
be the linear functional obtained by composing this isomorphism
with the projection of $C(V_+)=\Lambda (V_+)$ onto $\Lambda^2V_+$.
We obtain an alternating $F$-form on $W$ by
$$<x,y> =\lambda (x'y)\ \ .$$

\proclaim{Lemma} For $c\in C(V)$, and for $x$ and $y\in W$,
$$<\sigma(c)x,y>=<x,\sigma(c')y>.$$
In particular, for $g\in G=GSpin(V)$,
$$<\sigma(g)x,\sigma(g)y> = \nu(g)<x,y>.$$
Here $\sigma(g)$ denotes the spin representation action of $g$ on $W$,
and $\nu:G\lra F^\times$,  $\nu(g)=gg'$ is the restriction to $G$ of the
spinor norm on $C(V)$.
\endproclaim
\demo{Proof} Choose a basis $e_0,\ e_1,\ v_0,\ f_0,\ f_1$ for
$V$ such that the matrix for the quadratic form is
$$\pmatrix {}&{}&{}&1&0\\
{}&{}&{}&0&1\\ {}&{}&1&{}&{}\\ 1&0&{}&{}&{}\\
0&1&{}&{}&{}\endpmatrix.$$ In $C(V)$, $v_0^2=1$, $e_0f_0+f_0e_0 =
1$, $e_1f_1+f_1e_1=1$, $e_0^2=0$, $v_0(1+v_0)=(1+v_0)$, etc. The
spin representation $W = C(V_+)(1+v_0)$ has basis  $(1+v_0)$,
$e_0(1+v_0)$, $e_0e_1(1+v_0)$, and $e_1(1+v_0)$. We take $\lambda$
to be the coefficient of $e_0e_1(1+v_0)$ and the symplectic form
has matrix
$$J=\pmatrix {}&1_2\\-1_2&{}\endpmatrix.$$

It is easy to check that
$$\sigma(e_0) = \pmatrix
0&0&0&0\\
1&0&0&0\\
0&0&0&1\\
0&0&0&0\endpmatrix$$
$$\sigma(e_1) = \pmatrix
0&0&0&0\\
0&0&0&0\\
0&-1&0&0\\
1&0&0&0\endpmatrix$$
$$\sigma(v_0) = \pmatrix
1&0&0&0\\
0&-1&0&0\\
0&0&1&0\\
0&0&0&-1\endpmatrix$$
$$\sigma(f_0) = \pmatrix
0&1&0&0\\
0&0&0&0\\
0&0&0&0\\
0&0&1&0\endpmatrix$$
and
$$\sigma(f_1) = \pmatrix
0&0&0&1\\
0&0&-1&0\\
0&0&0&0\\
0&0&0&0\endpmatrix.$$
If $\sigma(c)$ is any of these matrices, then
$J{}^t\sigma(c)J^{-1} = \sigma(c)$, and hence, for any $c\in C(V)$,
$J{}^t\sigma(c)J^{-1} = \sigma(c')$, as claimed.
\qed\enddemo

\proclaim{Corollary} $\sigma: G=GSpin(V)\isoarrow GSp(W).$
\endproclaim

{\bf A.3.} In this section $F$ is again arbitrary, of
characteristic not 2.

\proclaim{Lemma}
Let $(V,q)$ be a non-degenerate quadratic space of dimension 5. The
subspace $\delta\cdot V\subset C^+(V)$ is characterized as:
$$\delta\cdot V=\{\ x\in C^+(V)\mid x'=x\text{ and } tr(x)=0\ \}.$$
\endproclaim

\demo{Proof}
Recall that $\delta\in C^-(V)$ is central in $C(V)$ and satisfies
$\delta'=\delta$. It is, thus, clear that $x=\delta v$ satisfies
$x'=x$. On the other hand, $x^2 = q(v)\delta^2=a$ lies in $F$, the
center of $C^+(V)$. In addition, if $x\ne 0$, then $x$ cannot lie
in the center of $C^+(V)$, since, if it did, then $v=\delta^{-1}x$
would lie in the center of $C(V)$, and this is not the case. If
$a=0$, so that $x^2=0$, the condition $tr(x)=0$ is immediate. If
$x^2=a\ne 0$, choose $u\in V$ with $q(u)\ne 0$ but with $(u,v)=0$,
and set $y=\delta u$. Then $xy=-yx$, and so, over an algebraic
closure of $F$, left multiplication  by $y$ gives an isomorphism
between the $\pm\sqrt{a}$ eigenspaces of $x$, and thus these spaces
have the same dimension and $tr(x)=0$. This proves that $\delta V$
is contained in the space on the right hand side. The converse
inclusion will be proved further down.
\qed\enddemo

Let $B$ be a quaternion algebra over $F$ with main involution
$\iota$, and let $C=M_2(B)$ with involution $x\mapsto x'
= {}^tx^\iota$. Let
$$\align
V_B&=\{\ x\in C\ \mid\ x'=x \text{ and } tr(x)=0\ \}\\
{}&=\{\ x=\pmatrix a&b\\b^\iota&-a\endpmatrix\ \mid \ a\in F, \ b\in B\ \}.
\endalign$$
Note that
$$xx'= x^2 = \pmatrix a^2+\nu(b)&{}\\{}&a^2+\nu(b)\endpmatrix,$$
so that the inclusion $V_B\hookrightarrow M_2(B)$ induces
a homomorphism
$$C(V_B,q_B) \lra M_2(B),$$
where the quadratic form on $V_B$ is $q_B(x) = xx'$ . The diagram
$$
\matrix
C(V_B,q_B)&\lra&M_2(B)\\
\nass
\downarrow{\prime}&{}&\downarrow{\prime}\\
\nass
C(V_B,q_B)&\lra&M_2(B)\endmatrix$$ commutes, and induces an
isomorphism $C^+(V_B,q_B) \isoarrow M_2(B)$, compatible with the
involutions.

Conversely, let $V$ be a nondegenerate quadratic space of dimension
5. The Clifford involution induces an isomorphism $C^+(V)\simeq
C^+(V)^{\roman{op}}$, hence $C^+(V)$ is of the form
$$C^+(V)\simeq M_2(B)\ \ ,$$
for a quaternion algebra $B$ over $F$. We may choose the
isomorphism compatible with the involutions $x\mapsto x'$. This map
then carries $\delta V$ into $V_B$. For dimension reasons we obtain
an isometry,
$$(V,\delta^2\cdot q_V)\simeq (V_B, q_B)\ \ .$$
This also concludes the proof of the lemma above.

\proclaim{Corollary}
$$G=\{ g\in C^+(V)^{\times}\mid gg'= \nu (g)\}\ \ .$$
\endproclaim

{\bf A.4.} Any involution of the central simple algebra $C=M_n(B)$,
has the form $x\mapsto h x' h^{-1}$ where $x'={}^tx^\iota$, where
$h\in GL_n(B)$ with $h'=\pm h$.  If $h'=h$, we say that the
involution is of {\it main type}, while, if $h'=-h$, we say that
$h$ is of {\it nebentype}. As observed above, the Clifford
involution on $M_2(B)$ is of main type.

Let $E$ be a central simple algebra over $F$, with $\dim_FE=16^2$,
and with a nontrivial involution $x\mapsto x^\eta$ whose
restriction to $F$ is trivial. Then there is a quaternion algebra
$B$ over $F$ and an isomorphism $E\simeq M_8(B)$. For a quaternion
algebra $B_1$ over $F$, let $C_1=M_2(B_1)$ and let $x\mapsto
x^{\eta_1}$ be an involution of $C_1$ whose restriction to $F$ is
trivial. Suppose that there is a (unitary) homomorphism
$$i_1:C_1=M_2(B_1)\hookrightarrow E=M_8(B)$$
such that
$$i_1(c)^\eta = i_1(c^{\eta_1}).$$
Let
$$C_2 = Cent_E(i_1(C_1))$$
be the centralizer of the image of $C_1$ and let $i_2:C_2\hookrightarrow E$
be the natural inclusion. Then $C_2\simeq M_2(B_2)$, where
$$B_1\tt B_2\simeq M_2(B),$$
and we have an isomorphism
$$i=i_1\tt i_2:C_1\tt C_2 \isoarrow E$$
such that
$$i(c_1\tt c_2)^\eta = i(c_1^{\eta_1}\tt c_2^{\eta_2}),$$
for an involution $\eta_2$ of $C_2$.

\proclaim{Proposition} The types of the involutions
$\eta = \eta_1\tt\eta_2$ are:
$$\text{main} = \cases \text{ main }\tt\text{ neben }\\
\text{ neben }\tt\text{ main }\endcases$$
and
$$\text{neben} = \cases \text{ main }\tt\text{ main }\\
\text{ neben }\tt\text{ neben. }\endcases$$
\endproclaim
\demo{Proof} We can assume that $F$ is algebraically closed.
Then, on $B=B_1=B_2=M_2(F)$,
$$\pmatrix a&b\\c&d\endpmatrix^\iota =\pmatrix d&-b\\-c&a\endpmatrix,$$
as usual, and a main involution on
$C_1=C_2=M_2(B)$ or on $E=M_8(B)$,
is given by $x\mapsto {}^tx^\iota$. On $M_n(B)\simeq M_{2n}(F)$,
this amounts to applying transpose on the matrix of the $2\times 2$ blocks
and then applying $\iota$ blockwise. We denote this type of transpose by
$x\mapsto {}^tx$ and write $x\mapsto {}^Tx$ for the usual transpose on
$M_{2n}(F)$.
Let $\tau=\pmatrix {}&1\\-1&{}\endpmatrix\in B= M_2(F)$, so that
$\tau^\iota = -\tau$, and, for $x\in B$,
$$\tau x^\iota \tau^{-1} = {}^Tx.$$
Setting
$$h=h_n=\text{diag}(\tau,\dots,\tau)$$
in $M_n(B)$, we have involutions of nebentype
$$x\mapsto h {}^tx^\iota h^{-1}$$
which, on $M_n(B)\simeq M_{2n}(F)$ are just given by
$x\mapsto {}^Tx$, the {\it usual transpose}, rather than
the blockwise transpose.

Now consider the explicit isomorphism
$$\align
i:M_4(F)\tt M_2(B) &\isoarrow M_8(B)\\
(a_{ij})\tt y &\mapsto (a_{ij}y).\endalign$$
Applying the involution of main type on $M_8(B)$, we have
$${}^ti(x\tt y)^\iota = i({}^Tx\tt {}^ty^\iota).$$
Similarly, applying the involution of nebentype on $M_8(B)$, we have
$${}^Ti(x\tt y) = i({}^Tx\tt h {}^ty^\iota h^{-1}) = i({}^Tx\tt {}^Ty).$$

Every involution on $E$ compatible with the isomorphism
$i:C_1\tt C_2\isoarrow E$ is conjugate to one of these two by
an element of the form $g=i(g_1\tt g_2)$, with ${}^tg^\iota =\pm g$. Note that
$${}^tg^\iota = g \iff \cases {}^Tg_1=g_1\text{ and } {}^tg_2^\iota = g_2&\\
{}^Tg_1=-g_1\text{ and } {}^tg_2^\iota =- g_2,&\endcases$$
and
$${}^tg^\iota = -g \iff \cases {}^Tg_1=g_1\text{ and } {}^tg_2^\iota = -g_2&\\
{}^Tg_1=-g_1\text{ and } {}^tg_2^\iota = g_2.&\endcases$$
Also observe that the involution
$$x\mapsto g_1{}^Tx g_1^{-1} = g_1h{}^tx^\iota h^{-1} g_1^{-1}$$
is of main type if ${}^Tg_1 = -g_1$ and of nebentype if ${}^Tg_1=g_1$, since
$${}^Tg_1 = h{}^tg_1h^{-1}=\pm g_1 \iff \pm {}^t(g_1h)^\iota = -g_1h.$$
\qed\enddemo

{\bf A.5.} Let $B$ be a quaternion algebra over $\R$.
\proclaim{Lemma} For $\tau\in B^\times$ with
$\tau^\iota=\pm\tau$, the involution $x\mapsto x^* = \tau x^\iota \tau^{-1}$ on
$B$
is positive if and only if:
$$\cases \tau^\iota = -\tau \text{ and $\tau^2<0$} &\text{ if $B=M_2(\R)$ }\\
\tau^\iota=\tau &\text{ if $B=\H$ is division.}\endcases$$
In particular, if $B=\H$, then $x^*=x^\iota$ is the unique positive involution
on $B$.
\endproclaim
\demo{Proof} Take $\tau\in B^\times$ such that $\tau^\iota=-\tau$
and $\tau^2<0$. Note that the condition on $\tau^2$ is automatic
when $B=\H$. Choose $\eta\in B^\times$ such that $\eta\tau=-\tau\eta$
and $\eta^\iota=-\eta$. Then every element $x\in B$ can be written uniquely
in the form $x= a+ b\eta$ with $a$ and $b\in \R(\tau)\simeq \C$. Then
$$x^* = \tau( a+b\eta)^\iota \tau^{-1} = a^\iota -\eta^\iota b^\iota$$
and
$$\align
tr(xx^*)&=tr((a+b\eta)(a^\iota -\eta^\iota b^\iota))\\
{}&=tr(aa^\iota +b\eta a^\iota - a\eta^\iota b^\iota - b \eta \eta^\iota
b^\iota)\\
{}&=2(aa^\iota + bb^\iota \eta^2).\endalign$$
If $B=M_2(\R)$, then $\eta^2>0$, and this quantity is positive, while,
if $B=\H$, then $\eta^2<0$, and this quantity can be negative. Note that,
when $B=M_2(\R)$, then an involution defined by a $\tau$ with $\tau^2>0$
cannot be positive.
\qed\enddemo

{\bf A.6.} Let $B=M_2({\Bbb R})$ and let $C=M_2(B)$ with involution
$x'={}^tx^{\iota}$ as above and let
$$V=\{ x\in C;\ x'=x\ \text{and}\ {\roman{tr}}(x)= 0\}\ \ .$$
Then the signature of $V$ for the form $q_B$ of A.3 is (3,2).
%\bye
\comment

{\it The following material was commented out in the original printout}

Let $B$ be an indefinite quaternion algebra over $\Q$, let
$C=M_2(B)$, with involution $x'={}^tx^\iota$, as above, and let
$$V=\{\ x\in C\ \mid \ x'=x \text{ and } tr(x)=0\ \}.$$
Note that the signature of $V$ for the form $q_B$ of 1.4 is $(3,2)$.

Fix a maximal order $\Cal O_B$ in $B$ such that $\Cal O_B^\iota = \Cal O_B$,
and let $\Cal O_C=M_2(\Cal O_B)$. Let $D(B)$ be the product of the
primes at which $B_p$ is division, and choose $\tau\in B^\times$
such that $\tau^\iota=-\tau$, $\tau^2=-D(B)$, and
$\tau \Cal O_B \tau^{-1}=\Cal O_B$.  By 1.7, the map
$x\mapsto x^*=\tau x^\iota \tau^{-1}$ is a positive involution of $B$
preserving $\Cal O_B$.
Also, for
$$\a=\pmatrix \tau&{}\\{}&\tau\endpmatrix\in M_2(B)$$,
$\a'=-\a$ and $x^*=\a x' \a^{-1} = \a^{-1}x'\a$ is a positive involution of
$C$,
preserving $\Cal O_C$.

Let
$$G=GSpin(V)=\{g\in C^\times\mid gg'=\nu(g)\ \},$$
so that $G$ is a twisted form over $\Q$ of $GSp_4$.

Let $U=\Cal O_C$, viewed as a module for $\Cal O_C$ under both left and right
multiplication. Define an alternating form:
$$<\ ,\ >\ :U\times U\lra \Z$$
by
$$<x,y>\ = tr(y'\a^{-1} x).$$
Then
$$<cx,y>\ = tr(y'\a^{-1}cx) = tr(y'\a^{-1}c\a\a^{-1} y) =\ <x,c^*y>,$$
and
$$<xc,y>\ = tr(y'\a^{-1}xc) = tr(cy'\a^{-1} x) =\ <x,yc'>.$$

Let $D$ be the space of oriented negative $2$-planes in $V(\R)$. This space has
two components and the group $G(\R)$ acts transitively on it, via its action
on $V(\R)$. For $z\in D$, let $z_1$, $z_2\in z$ be a properly oriented basis
such that the restriction of the quadratic form $q_B$ from $V(\R)$ to $z$
has matrix $-1_2$ for the basis $z_1$, $z_2$. Let
$j_z = z_1z_2\in C(\R)$. Viewing $j_z$ as the image of the element $z_1z_2\in
C(V(\R))$,
the Clifford algebra of $V(\R)$ and recalling the commutative diagram of 1.4,
we see that $j_z'=-j_z$ and that $j_z^2=-z_1^2z_2^2=-1$.  Hence, $j_zj_z'=1$
and so,
$j_z\in G(\R)$, and
$$<x j_z,y j_z>\ =\ <x,y>.$$

In fact, there is an isomorphism
$$\C\isoarrow C^+(z)\qquad i\mapsto z_1z_2,$$
where $C^+(z)$ is the even Clifford algebra of the real $2$-plane $z$.
The composition of this map with the map
$$C^+(z)\subset C^+(V(\R)) \isoarrow M_2(B(\R))$$
of 1.4 induces a morphism, defined over $\R$,
$h_z:\Bbb S \lra G$.
Note that $h_z(i) = j_z$.
The space $D$ can thus be viewed as the space of conjugacy
classes of such maps under the action of the group $G(\R)$.

For $\tau\in B^\times$, as above, let
$$\tau_0= D(B)^{-\frac12} \tau\in B^\times(\R),$$
so that $\tau_0^2=-1$.
Choose $\eta\in B^\times$ such that
$$\eta\tau=-\tau\eta,\qquad\text{ and }\qquad \eta^\iota=-\eta.$$
Since $B$ is indefinite, $\eta^2>0,$ and we can set
$$\eta_0=(\eta^2)^{-\frac12}\eta\in B^\times(\R),$$
so that $\eta_0^2=1$. The vectors
$$\pmatrix {}&\eta_0\\ -\eta_0&{}\endpmatrix, \qquad\text{ and } \qquad
\pmatrix {}&\tau_0\eta_0\\-\tau_0\eta_0&{}\endpmatrix \in V(\R)$$
form a standard basis of an oriented negative $2$-plane $z_0\in D$,
and
$$j_{z_0}=\pmatrix {}&\eta_0\\ -\eta_0&{}\endpmatrix
\pmatrix {}&\tau_0\eta_0\\-\tau_0\eta_0&{}\endpmatrix = \pmatrix
\tau_0&{}\\{}&\tau_0
\endpmatrix
=D(B)^{-\frac12} \a =:\a_0.$$

\proclaim{Lemma} For any $z\in D$,
$$<xj_z,y>\ =\ < y j_z,x>,$$
and, for $x\in U(\R)$, $x\ne 0$,
$$<x j_z,x>\ >0,$$
if $z$ lies in the same connected component of $D$ as $z_0$, and
$$<x j_z,x>\ <0,$$
if $z$ and $z_0$ lie in opposite components.
\endproclaim
\demo{Proof} For the first assertion:
$$<x j_z,y>\ =\ -<x,y j_z>\ =\ <y j_z,x>.$$
For the second, write $z= g z_0$ for $g\in G(\R)$, so that
$$j_z= g j_{z_0} g^{-1}= g \a_0 g^{-1}.$$
Then, we have
$$\align
<x j_z,x>&=\ <x g \a_0 g^{-1}, x>\\
{}&=\nu(g)^{-1} <x g \a_0,x g > \\
{}&= \nu(g)^{-1}tr((xg)'\a^{-1} xg\a_0)\\
{}&= \nu(g)^{-1} D(B)^{-\frac12}\, tr(\a (xg)' \a^{-1} (xg))\\
{}&= \nu(g)^{-1} D(B)^{-\frac12}\,tr((xg)^*(xg)).
\endalign$$
Since $x\mapsto x^*$ is a positive involution, this gives the claim.
\qed\enddemo

Let $D^+$ be the connected component of $D$ containing $z_0$. Then,
for any $z\in D^+$, we obtain a (principally) polarized abelian variety
$$A_z=(U(\R), j_z, U(\Z), <\ ,\ >)$$
with $\dim A_z=8$ and with an action
$$\iota: \Cal O_C\hookrightarrow End(A_z).$$
If
$$\gamma\in \Gamma =\{\ g\in G(\Q)\mid U(\Z) g = U(\Z)\ \},$$
then, right multiplication by $\gamma^{-1}$ induces an isomorphism
$$A_z\isoarrow A_{\gamma z}.$$

{\bf 1.} Let $B$ be a quaternion algebra over $F$ with main involution $\iota$,
and let
$$V_0=\{\ x\in B\mid tr(x)=0\}.$$
For $x\in V_0$, $x^2=-\nu(x)$ so that $(V_0,-\nu)$ is a quadratic space.
The inclusion $V_0\hookrightarrow B$ induces an algebra homomorphism
$$\ph_0:C(V_0,-\nu)\lra B$$
of the Clifford algebra $C(V_0,-\nu)$. Writing $B$ as a cyclic algebra $(a,b)$
with basis $1,\ i,\ j,\ ij=k$,
$$C(V_0,-\nu)=C(V_0)\simeq C^+(V_0)\tt F(\delta_0)$$
where $\delta_0=ijk$ (product in $C(V_0)$). Then
$$\ph_0(\delta_0) = ijk = k^2 =-ab$$
lies in $F$ and so, we obtain an isomorphism
$$\ph_0:C^+(V_0,-\nu) \simeq B.$$
The image of $V_0$ and the element $\delta_0$ lie in the odd part
$C^-(V_0)$ so that $\delta_0V_0\subset C^+(V_0)$. Note that
$\delta_0^\iota=-\delta_0$, so that
$$\delta_0V_0=\{\ x\in C^+(V_0)\mid tr(x)=0\ \},$$
and  $GSpin(V_0)\simeq B^\times$.

\endcomment

\hfuzz=20pt
\widestnumber\key{666}

\vfill\eject
\subheading{References}

\ref\key{\boutotcarayol} \by J.-F. Boutot and H. Carayol   %\boutotcarayol
\paper Uniformisation p-adique des courbes de Shimura: les th\'eor\`emes
de Cerednik et de Drinfeld
\inbook Courbes modulaires et courbes de Shimura
\publ Ast\'erisque
%\vol 196--197
{\bf 196-197}
\yr 1991
\pages 45--158
\endref

\ref\key{\deligne}
\by P. Deligne
\paper Travaux de Shimura
\jour Sem.\ Bourbaki 389, Springer Lecture Notes
%\vol 244
{\bf 244}
\publaddr Berlin
\yr 1971
\endref

\ref\key{\delignetwo}
\by P. Deligne
\paper La conjecture de Weil pour les surfaces $K 3$
\jour Invent.\ math.\
%\vol 244
{\bf 15}
\yr 1972
\pages 206-226
\endref

\ref\key{\eichler}
\by M. Eichler
\book Quadratische Formen und orthogonale Gruppen
\publ Springer-Verlag
\publaddr Berlin
\yr 1974
\endref

\ref\key{\fulton}
\by W.\ Fulton
\book Intersection theory
\publ Springer-Verlag
\publaddr Berlin
\yr 1984
\endref

\ref\key{\genestier}
\by A. Genestier
\paper Letter to M. Rapoport
\pages August 12, 1996
\endref

\ref\key{\grosskeating}
\by B. H. Gross and K. Keating
\paper  On the intersection of modular correspondences
\jour Invent.\ math.
%\vol 112
{\bf 112}
\yr 1993
\pages 225--245
\endref

\ref\key{\hashimoto}
\by K. Hashimoto
\paper Class numbers of poitive definite ternary quaternion hermitian forms
\jour Proc. Japan Acad.
%\vol 59
{\bf 59}
\yr 1983
\pages 490--493
\endref

\ref\key{\hashibuki}
\by K. Hashimoto and T. Ibukiyama
\paper On class numbers of positive definite binary quaternion hermitian forms
\jour J. Fac. Sci. Univ. Tokyo, Sect IA
%\vol 27   %28  %% to appear
{\bf 27}
\yr 1980  %1981  %%??
\pages 549--601 % 695--699
\endref

\ref\key{\ibukatoort}
\by T. Ibukiyama, T. Katsura, and F. Oort
\paper Supersingular curves of genus two and class numbers
\jour Compositio Math.
%\vol 57
{\bf 57}
\yr 1986
\pages 127--152
\endref

\ref\key{\kaiser}
\by C. Kaiser
\paper Ein getwistetes fundamentales Lemma f\"ur die $\text{\rm GSp}_4$
\jour Dissertation, Bonn
\yr 1997
\endref

\ref\key{\katoort}
\by T. Katsura and F. Oort
\paper Families of supersingular abelian surfaces
\jour Comp.\ Math.
%\vol 62
{\bf 62}
\yr 1987
\pages 107--167
\endref

\ref\key{\kitaokanote}
\by Y. Kitaoka
\paper A note on local representation densities of quadratic forms
\jour Nagoya Math. J.
%\vol 92
{\bf 92}
\yr 1983
\pages 145--152
\endref

\ref
\key{\kitaokatwo}
\by Y. Kitaoka
\paper Fourier coefficients of Eisenstein series of degree 3
\jour Proc. of Japan Akad.
%\vol 60
{\bf 60}
\yr 1984
\pages 259--261
\endref

\ref
\key{\kitaokabook}
\by Y. Kitaoka
\book Arithmetic of Quadratic Forms
\publ Cambridge Univ. Press
\yr 1993
\bookinfo Cambridge Tracts in Mathematics {\bf 106}
\publaddr Cambridge, U.K.
\endref

\ref\key{\kottwitztam}
\by R. Kottwitz
\paper Tamagawa numbers
\jour Ann.\ of Math.
%\vol 127
{\bf 127}
\yr 1988
\pages 629--646
\endref

\ref\key{\kottwitz}
\by R. Kottwitz
\paper Points on some Shimura varieties over finite fields
\jour J.\ AMS
%\vol 5
{\bf 5}
\yr 1992
\pages 373--444
\endref

\ref\key{\cycles}   % Duke paper
\by S. Kudla
\paper Algebraic cycles on Shimura varieties of orthogonal type
\jour Duke Math. J.
%\vol 86
{\bf 86}
\yr 1997
\pages 39--78
\endref

\ref\key{\ked}   % Annals paper
\by S. Kudla
\paper Central derivatives of Eisenstein series and height pairings
\jour Ann.\ of Math.
\yr to appear
\endref

\ref\key{\kmihes}
\by S. Kudla and J. Millson
\paper Intersection numbers of cycles on locally symmetric spaces and Fourier
coefficients of holomorphic modular forms in several complex variables
\jour Publ.\ math.\ IHES
%\vol 71
{\bf 71}
\yr 1990
\pages 121--172
\endref

\ref\key{\krhb}
\by S. Kudla and M. Rapoport
\paper Arithmetic Hirzebruch Zagier cycles
\jour in preparation
\endref

\ref\key{\krdrin}
\by S. Kudla and M. Rapoport
\paper Height parings on Shimura curves and $p$-adic uniformization
\jour in preparation
\endref

\ref\key{\moretbailly}
\by L. Moret-Bailly
\paper  Familles de courbes et de vari\'et\`es ab\`eliennes sur
${\Bbb P}^1$, II.\ exemple
\inbook L.\ Szpiro (ed.): S\'eminaire sur les pinceaux de courbes de
genre au moins deux
\publ Ast\'erisque
%\vol 86
{\bf 86}
\yr 1981
\pages 125-140
\endref

\ref\key{\oort}
\by F. Oort
\paper Which abelian surfaces are products of elliptic curves?
\jour Math. Ann.
%\vol 214
{\bf 214}
\yr 1975
\pages 35--47
\endref

\ref\key{\ribetone}
\by K. Ribet
\paper On modular representations of $\text{\rm Gal}(\bar{\Bbb Q})/\Bbb Q)$
arising from modular forms
\jour Invent.\ Math.
%\vol 100
{\bf 100}
\yr 1990
\pages 431--476
\endref

\ref\key{\ribettwo}
\by K. Ribet
\paper Bimodules and abelian surfaces
\jour Advanced Studies in Pure Math.
%\vol 17
{\bf 17}
\yr 1989
\pages 359--407
\endref

\ref\key{\serre}
\by J.-P.\ Serre
\paper Alg\`ebre locale. Multiplicit\' es.
\jour Springer Lecture Notes
%\vol 11
{\bf 11}
\publaddr Berlin
\yr 1965
\endref

\ref
\key{\shih}
\by K.-Y. Shih
\paper Existence of certain canonical models
\jour Duke Math. J.
%\vol 45
{\bf 45}
\yr 1978
\pages 63--66
\endref

\ref
\key{\shimura}
\by G.\ Shimura
\paper Arithmetic of alternating forms and quaternion hermitian forms
\jour J.\ Math.\ Soc.\ Japan
%\vol 15
{\bf 15}
\yr 1963
\pages 63--65
\endref

\ref\key{\stamm}
\by H. Stamm
\paper On the reduction of the Hilbert-Blumenthal-moduli scheme with
$\Gamma_0(p)$-level structure
\jour Forum Math.
%\vol 9
{\bf 9}
\yr 1997
\pages 405--455
\endref

\ref\key{\sweet}
\by W. J. Sweet
\paper A computation of the gamma matrix of
a family of $p$-adic zeta integrals
\jour J. Number Theory
%\vol
{\bf 55}
\yr 1995
\pages 222-260
\endref

\ref\key{\tonghai}
\by Tonghai Yang
\paper An explicit formula for local densities of quadratic forms
\jour preprint MPI Bonn
\yr 1997
\endref
\vskip3 cm
\noindent
\line{Stephen S.\ Kudla\hfill Michael Rapoport}
\line{Department of Mathematics\hfill Mathematisches Institut}
\line{University of Maryland\hfill der Universit\"at zu K\"oln}
\line{College Park, MD 20742\hfill Weyertal 86-90}
\line{{}\hfill D -- 50931 K\"oln}
\line{USA\hfill Germany}

\bye